\documentclass[3p,times]{elsarticle}
\usepackage[latin9]{inputenc}

\usepackage{kbordermatrix}

\usepackage{psfrag,
	makecell,
	epsfig,
	epstopdf,
	amsthm,
	amssymb,
	url,
	subfigure,
	algorithm,
	balance,
	enumerate,
	color,
	setspace,
	tikz,
	pgfplots,
	amsmath,	
	float
}

\usepackage{verbatim}
\usepackage{lmodern}
\usepackage{graphicx}
\usepackage{amsmath,mathrsfs,amssymb}
\usepackage{amsfonts}
\usepackage{caption}
\usepackage{array}
\usepackage{longtable}
\usepackage{algpseudocode}
\usepackage{amsmath}
\usepackage{mathpazo}
\usepackage{float}
\usepackage{algorithm}
\graphicspath{ {images/} }
\usepackage{epstopdf}
\usepackage{enumitem}
\usepackage{booktabs}
\usepackage{todonotes}
\usepackage{tabulary}
\usepackage{multirow}
\usepackage{textcomp}
\usepackage{longtable}
\usepackage{lipsum}
\usepackage{amsbsy}
\usepackage{amstext}
\usepackage{amsthm}
\usepackage{amssymb}
\usepackage{mathdots}
\DeclareTextSymbolDefault{\textquotedbl}{T1}
\usepackage[unicode=true,
bookmarks=false,
breaklinks=false,pdfborder={0 0 1},colorlinks=false]
{hyperref}
\hypersetup{
	colorlinks,bookmarksopen,bookmarksnumbered,citecolor=red,urlcolor=red}

\makeatletter

\pdfpageheight\paperheight
\pdfpagewidth\paperwidth
\usepackage[english]{babel}
\usepackage{lineno,hyperref}
\modulolinenumbers[5]
\providecommand{\tabularnewline}{\\}
\usepackage{babel}

\newtheorem{lem}{Lemma}

\begin{document}

\begin{frontmatter}

\title{Special Orthogonal Group SO(3), Euler Angles, Angle-axis, Rodriguez Vector and Unit-Quaternion: Overview, Mapping and Challenges}

\author[HAH]{Hashim A.~Hashim}
\ead{hhashim@tru.ca}

\address[HAH]{Department of Engineering and Applied Science, Thompson Rivers University, BC, Canada, V2C-0C8.}
\vspace{9cm}
\begin{abstract}
The attitude of a rigid-body in the three dimensional space has a
unique and global definition on the Special Orthogonal Group $\mathbb{SO}\left(3\right)$.
This paper gives an overview of the rotation matrix, attitude dynamics
and parameterization. The four most frequently used methods of attitude
representations are discussed with detailed derivations, namely Euler
angles, angle-axis parameterization, Rodriguez vector, and unit-quaternion.
The mapping from one representation to others including $\mathbb{SO}\left(3\right)$
is given. Also, important results which could be useful for the process
of filter and/or control design are given. The main weaknesses of
attitude parameterization using Euler angles, angle-axis parameterization,
Rodriguez vector, and unit-quaternion are illustrated.\\
\noindent\rule{16.5cm}{1pt}
\vspace{0.2cm}
\end{abstract}
\begin{keyword}
Special Orthogonal Group 3, Euler angles, Angle-axis,
Rodriguez Vector, Unit-quaternion, $SO\left(3\right)$, Mapping, Parameterization,
Attitude, Control, Filter, Observer, Estimator, Rotation, Rotational
matrix, Transformation matrix, Orientation, Transformation, Roll,
Pitch, Yaw. 
\end{keyword}
\end{frontmatter}{}
\vspace{2.4cm}

\noindent\rule{16.5cm}{2pt}

\textbf{\textcolor{red}{Bibtex formatted citatio}}\textcolor{red}{n}:\\

@article\{hashim2019special, 

title=\{Special Orthogonal Group SO(3),
Euler Angles, Angle-axis, Rodriguez Vector and Unit-Quaternion: Overview,
Mapping and Challenges\}, 

author=\{Hashim A. Hashim\}, 

journal=\{ArXiv preprint ArXiv:1909.06669\}, 

year=\{2019\} 

\}

\noindent\rule{16.5cm}{2pt}
\newpage
\tableofcontents{}
\newpage

\listoftables
\noindent\rule{16.5cm}{2pt}

\section{Introduction}

Automated or semi-automated systems have been proven to be indispensable
in the majority of contemporary engineering applications. Rigid-bodies
rotating and/or moving in space constitute a fundamental part of these
applications, which include unmanned aerial vehicles (UAVs), satellites,
spacecrafts, rotating radars, underwater vehicles, ground vehicles,
robotic systems, etc. \textbf{Rotational matrix}, also termed \textquotedbl\textbf{attitude}\textquotedbl ,
describes the orientation of a rigid body rotating in the three dimensional
(3D) space. Thus, in this article the words \textbf{attitude}, \textbf{orientation}
or \textbf{rotational matrix} will be used interchangeably. The estimation
and/or control process of a rigid-body moving and rotating in space
is defined by pose (attitude and position) problem. Obviously, attitude
is a major part of the pose problem. The successful estimation and/or
control process of a rigid-body rotating and/or moving in space depends
primarily on the precise representation of the attitude of the rigid-body.
There are several approaches to attitude parameterization of a rigid-body
in the 3D space. Some of those approaches use $3\times3$ parameterization,
while others use three or four components. The attitude can be parameterized
through a $3\times3$ orthogonal matrix, which is known to follow
the Special Orthogonal Group $\mathbb{SO}\left(3\right)$. Attitude
parameterization could also be achieved by three components through
making use of Rodriguez vector or Euler angles (roll, pitch, and yaw).
Additionally, the attitude can be parameterized using four components,
for instance, unit-quaternion and angle-axis parameterization. The
mapping of the three components (Euler angles or Rodriguez vector)
or four components (unit-quaternion or angle-axis) to $\mathbb{SO}\left(3\right)$
must have orthogonal configuration with the determinant equal to one.
The main problem of the $3\times3$ configuration is that it is neither
Euclidean nor has a vector form. Accordingly, the majority of researchers
prefer to consider the vector form in the process of estimation and/or
control of a rigid-body rotating and/or moving in space. The Euclidean
parameterization includes Rodriguez vector and Euler angles which
lie in $\mathbb{R}^{3}$. Unit-quaternion vector is used extensively
in the estimation and/or control process of a rigid-body rotating
and/or moving in space, in spite of the fact that it is non-Euclidean
and lies in the three-sphere $\mathbb{S}^{3}$. 

A more detailed discussion of the above-mentioned methods of attitude
parameterization and mapping from one method to the other have been
discussed in several journal articles, for instance \cite{hashim2020AtiitudeSurvey},  \cite{stuelpnagel1964parametrization},
\cite{wen1991attitude}, \cite{shuster1993survey}, and \cite{chaturvedi2011rigid}.
The main advantage of representing attitude directly on $\mathbb{SO}\left(3\right)$
is that a $3\times3$ set creates a global and unique representation
of the attitude problem. This implies that any orientation of a rigid-body
in 3D space has a unique rotational matrix. Rodriguez vector, Euler
angles and angle-axis parameterization fail to capture the attitude
at certain configurations and could be limited to singularity, called
``unstable equilibria''. On the other hand, unit-quaternion does
not have the problem of singularity, but it nonetheless suffers from
non-uniqueness in representation of the attitude. The above-named
deficiencies of the attitude representation methods have been considered
in \cite{stuelpnagel1964parametrization} and \cite{shuster1993survey}.
These shortcomings will be elaborated on in the subsequent sections.
The control and/or the estimation process of a rigid-body rotating
and/or moving in space could possibly fail if the vector associated
with the parameterization of the rigid body started at certain configurations. 

Tracing control of any rigid-body is mainly based on attitude. Over
the last few decades, tracking control of a rigid-body has gained
popularity and aroused interest of researchers in the control community.
These applications included spacecrafts (\cite{kang2009linear}, \cite{mohamed2014improved},
and \cite{stevens2015aircraft}), satellites \cite{yeh2002nonlinear},
mobile robots \cite{mohamed2014improved,chwa2010tracking}, etc. Attitude
and pose filtering problem is an essential task in robotics applications.
This filtering problem requires a set of measurements made with respect
to the inertial frame and obtained from the sensors attached to the
rigid-body. In the past, this problem used to be tackled through vehicles
equipped with high quality sensors \cite{hashim2018SO3Stochastic,markley2003attitude,hashim2020AtiitudeSurvey,crassidis2007survey,mahony2008nonlinear}.
High quality sensors are quite expensive and may not be an optimal
fit for small scale systems. With the introduction of micro elector-mechanical
systems (MEMSs) a range of inertial measurement units (IMUs) was proposed.
These units are small, inexpensive, and light-weight, allowing them
to be easily attached to small drones, mini UAVs, spacecrafts, radar,
satellites, mobile robots and other applications. Despite the advantages
that IMUs offer, their quality is fairly low and the measurements
they provide are prone to bias and noise.

The environment is unpredictable and the uncertainties present in
measurements are characterized by irregular behavior \cite{hashim2018SO3Stochastic,mahony2008nonlinear,hashim2019SE3Det,mohamed2019filters,hashim2020AtiitudeSurvey}
and additionally complicated by the fact that the initial conditions
may not be accurately known \cite{hashim2017neuro,hashim2017adaptive}.
Therefore, the attitude and pose filter should be particularly robust
against uncertainties in measurements and any initialized error, and
should be able to converge to the desired solution. Over the past
few decades, several successful filters have been proposed providing
high quality estimation of the true attitude or pose. The attitude
estimation problem used to be addressed either by a Gaussian filter
or a nonlinear deterministic filter. The majority of Gaussian filters
consider the unit-quaternion in the problem formulation and filter
structure, for instance \cite{lefferts1982kalman}, \cite{goddard1998pose},
\cite{choukroun2006novel}, \cite{hashim2020AtiitudeSurvey}, \cite{cheon2007unscented},
and \cite{barrau2015intrinsic}. Nonlinear deterministic filters are
directly developed on $\mathbb{SO}\left(3\right)$ such as, \cite{mahony2008nonlinear},
and \cite{hashim2019SO3Det}. Other filters were evolved directly
on $\mathbb{SO}\left(3\right)$ and $\mathbb{SE}\left(3\right)$,
considering Rodriguez vector and angle-axis parameterization in stability
analysis such as \cite{hashim2018Conf1,hashim2018SO3Stochastic,hashim2019SO3Det,mohamed2019filters,hashim2018SE3Stochastic,hashim2019SE3Det,hashim2019Conf1}.
The filters in \cite{hashim2018Conf1,hashim2018SO3Stochastic} are
nonlinear stochastic filters.

According to literature, the crucial factor of success of any control
and/or estimation process of a rigid-body rotating and/or moving in
space is primarily attributed to attitude parameterization. Therefore,
the first and main focus of this article is to give an overview of
different methods of attitude parameterization and ways to transition
between those methods. In fact, transitioning between methods could
help to avoid the complexities of the design process. The second purpose
is to outline important results associated with $\mathbb{SO}\left(3\right)$,
angle-axis parameterization, Rodriguez vector, and unit-quaternion
which could be very beneficial in the control and/or the estimation
process of a rigid-body rotating and/or moving in space. These results
could significantly simplify the process of filter and/or control
design.\textcolor{red}{}%

The remainder of the paper is organized in the following manner: Section
\ref{sec:Preliminaries-and-Math} presents an overview of all mathematical
and attitude notation, preliminaries to $\mathbb{SO}\left(3\right)$,
angle-axis, Rodriguez vector, and unit-quaternion parameterization,
and some useful identities. Section \ref{sec:OVERVIEW_SO3} gives
an overview of $\mathbb{SO}\left(3\right)$, attitude dynamics and
some important results. Section \ref{sec:Euler-Parameterization}
shows the mapping from/to Euler angles to/from other methods of attitude
representation, Euler rate, angular velocity transformation and the
problem of attitude parameterization using Euler angles. Section \ref{sec:OVERVIEW_Ang_Axis}
illustrates the mapping from/to angle-axis parameterization to/from
other methods of attitude representation, some important results and
the problem of attitude parameterization using angle-axis parameterization.
Section \ref{sec:OVERVIEW_ROD} parameterizes the attitude using Rodriguez
vector and shows the mapping from/to Rodriguez vector to/from other
methods of attitude representation, and Rodriguez vector dynamics.
Also, Section \ref{sec:OVERVIEW_ROD} presents some important results
and the problem of attitude parameterization using Rodriguez vector.
Section \ref{sec:Unit-quaternion} gives an overview of unit-quaternion,
shows the mapping from/to unit-quaternion to/from other attitude representations,
attitude dynamics, and the problem of attitude parameterization
using unit-quaternion. Section \ref{sec:Conclusion} concludes the
work by summarizing attitude parameterization and mapping from one
attitude representation to other representations in Table \ref{tab:SO3},
\ref{tab:Angle-axis}, \ref{tab:Rodriguez-vector}, and \ref{tab:Unit-Quaternion}. 

\noindent\makebox[1\linewidth]{%
	\rule{0.6\textwidth}{1.4pt}%
}
\newpage
\section{Preliminaries and Math Notation \label{sec:Preliminaries-and-Math}}

Table \ref{tab:Table-of-notation1} presents the important math notation
used throughout the paper. 

\begin{table}[H]
	\centering{}\caption{\label{tab:Table-of-notation1}Mathematical Notation}
	\begin{tabular}{lll}
		\toprule 
		\addlinespace
		$\mathbb{N}$  & :  & The set of integer numbers\tabularnewline
		\addlinespace
		$\mathbb{R}_{+}$  & :  & The set of nonnegative real numbers\tabularnewline
		\addlinespace
		$\mathbb{R}^{n}$  & :  & Real $n$-dimensional vector\tabularnewline
		\addlinespace
		$\mathbb{R}^{n\times m}$  & :  & Real $n$-by-$m$ dimensional matrix\tabularnewline
		\addlinespace
		$\left\Vert \cdot\right\Vert $  & :  & Euclidean norm, for $x\in\mathbb{R}^{n}$, $\left\Vert x\right\Vert =\sqrt{x^{\top}x}$\tabularnewline
		\addlinespace
		$\mathbb{S}^{1}$  & :  & Unit-circle, $\mathbb{S}^{1}=\left\{ \left.x=\left[x_{1},x_{2}\right]^{\top}\in\mathbb{R}^{2}\right|\sqrt{x_{1}^{2}+x_{2}^{2}}=1\right\} $ \tabularnewline
		\addlinespace
		$\mathbb{S}^{2}$  & :  & Two-sphere, $\mathbb{S}^{2}=\left\{ \left.x=\left[x_{1},x_{2},x_{3}\right]^{\top}\in\mathbb{R}^{3}\right|\left\Vert x\right\Vert =1\right\} $ \tabularnewline
		\addlinespace
		$\mathbb{S}^{n}$  & :  & $n$-sphere, $\mathbb{S}^{n}=\left\{ \left.x\in\mathbb{R}^{n+1}\right|\left\Vert x\right\Vert =1\right\} $ \tabularnewline
		\addlinespace
		$^{\top}$  & :  & Transpose of a component\tabularnewline
		\addlinespace
		$\times$ & : & Cross multiplication\tabularnewline
		\addlinespace
		$\mathbf{I}_{n}$  & :  & Identity matrix with dimension $n$-by-$n$, $\mathbf{I}_{n}\in\mathbb{R}^{n\times n}$ \tabularnewline
		\addlinespace
		$\mathbf{0}_{n}$  & :  & Zero matrix with dimension $n$-by-$n$, $\mathbf{0}_{n}\in\mathbb{R}^{n\times n}$ \tabularnewline
		\addlinespace
		${\rm det}\left(\cdot\right)$  & :  & Determinant of a component\tabularnewline
		\addlinespace
		${\rm Tr}\left\{ \cdot\right\} $  & :  & Trace of a component\tabularnewline
		\bottomrule
	\end{tabular}
\end{table}

Table \ref{tab:Table-of-notation2} provides some important definitions
and notation related to attitude. 

\begin{table}[H]
	\centering{}\caption{\label{tab:Table-of-notation2}Attitude Notation}
	\begin{tabular}{lll}
		\toprule 
		\addlinespace
		$\left\{ \mathcal{I}\right\} $ & :  & Inertial-frame of reference\tabularnewline
		\addlinespace
		$\left\{ \mathcal{B}\right\} $ & :  & Body-frame of reference\tabularnewline
		\addlinespace
		$\mathbb{GL}\left(3\right)$  & :  & The 3 dimensional General Linear Group\tabularnewline
		\addlinespace
		$\mathbb{O}\left(3\right)$  & :  & Orthogonal Group\tabularnewline
		\addlinespace
		$\mathbb{SO}\left(3\right)$  & :  & Special Orthogonal Group\tabularnewline
		\addlinespace
		$\mathfrak{so}\left(3\right)$  & :  & The space of $3\times3$ skew-symmetric matrices, and Lie-algebra
		of $\mathbb{SO}\left(3\right)$\tabularnewline
		\addlinespace
		$\left[\cdot\right]_{\times}$  & :  & Skew-symmetric of a matrix\tabularnewline
		\addlinespace
		$\boldsymbol{\mathcal{P}}_{a}$  & :  & Anti-symmetric projection operator\tabularnewline
		\addlinespace
		$R$  & :  & Attitude/Rotational matrix/Orientation of a rigid-body, $R\in\mathbb{SO}\left(3\right)$\tabularnewline
		\addlinespace
		$\Omega$  & :  & Angular velocity vector with $\Omega=\left[\Omega_{x},\Omega_{y},\Omega_{z}\right]^{\top}\in\mathbb{R}^{3}$\tabularnewline
		\addlinespace
		${\rm v}^{\mathcal{I}}$  & :  & Vector in the inertial-frame, ${\rm v}^{\mathcal{I}}\in\mathbb{R}^{3}$\tabularnewline
		\addlinespace
		${\rm v}^{\mathcal{B}}$  & :  & Vector in the body-frame, ${\rm v}^{\mathcal{B}}\in\mathbb{R}^{3}$\tabularnewline
		\addlinespace
		$\mathcal{R}_{\xi}$  & :  & Attitude representation obtained using Euler angles $\left(\phi,\theta,\psi\right)$,
		$\mathcal{R}_{\xi}\in\mathbb{SO}\left(3\right)$\tabularnewline
		\addlinespace
		$\xi$  & :  & Euler angle vector with, $\xi=\left[\phi,\theta,\psi\right]^{\top}\in\mathbb{R}^{3}$\tabularnewline
		\addlinespace
		$\mathcal{J}$  & :  & Transformation matrix between Euler rate and angular velocity vector,
		$\mathcal{J}\in\mathbb{R}^{3\times3}$\tabularnewline
		\addlinespace
		$\mathcal{R}_{\alpha}$  & :  & Attitude representation obtained using angle-axis parameterization,
		$\mathcal{R}_{\alpha}\in\mathbb{SO}\left(3\right)$\tabularnewline
		\addlinespace
		$\alpha$  & :  & Angle of rotation, $\alpha\in\mathbb{R}$\tabularnewline
		\addlinespace
		$u$  & :  & Unit vector (axis of parameterization), $u=\left[u_{1},u_{2},u_{3}\right]^{\top}\in\mathbb{S}^{2}\subset\mathbb{R}^{3}$\tabularnewline
		\addlinespace
		$\mathcal{R}_{\rho}$  & :  & Attitude representation obtained using Rodriguez vector, $\mathcal{R}_{\rho}\in\mathbb{SO}\left(3\right)$\tabularnewline
		\addlinespace
		$\rho$  & :  & Rodriguez vector, $\rho=\left[\rho_{1},\rho_{2},\rho_{3}\right]^{\top}\in\mathbb{R}^{3}$\tabularnewline
		\addlinespace
		$\mathcal{R}_{Q}$  & :  & Attitude representation obtained using unit-quaternion vector, $\mathcal{R}_{Q}\in\mathbb{SO}\left(3\right)$\tabularnewline
		\addlinespace
		$Q$  & :  & Unit-quaternion vector, $Q=\left[q_{0},q^{\top}\right]^{\top}=\left[q_{0},q_{1},q_{2},q_{3}\right]^{\top}\in\mathbb{S}^{3}\subset\mathbb{R}^{4}$\tabularnewline
		\addlinespace
		$Q^{*}$  & :  & Complex conjugate of unit-quaternion, $Q^{*}\in\mathbb{S}^{3}\subset\mathbb{R}^{4}$\tabularnewline
		\addlinespace
		$\odot$  & :  & Multiplication operator of two unit-quaternion vectors\tabularnewline
		\bottomrule
	\end{tabular}
\end{table}

Let $x\in\mathbb{R}^{3}$ and $M\in\mathbb{R}^{3\times3}$, and consider
$f\left(x,M\right)$ function such as
\[
f\left(x,M\right)=x^{\top}Mx\in\mathbb{R}
\]
where $^{\top}$ is the transpose of a component. Hence, the arbitrary
transformation through mapping is 
\[
f:\mathbb{R}^{3}\times\mathbb{R}^{3\times3}\rightarrow\mathbb{R}
\]

\subsection{$\mathbb{SO}\left(3\right)$ Preliminaries\label{sec:SO3PPF_Math-notation}}

Let $\mathbb{GL}\left(3\right)$ stand for the 3 dimensional general
linear group that is a Lie group characterized with smooth multiplication
and inversion. $\mathbb{O}\left(3\right)$ refers to the Orthogonal
Group and is a subgroup of the general linear group 
\begin{equation}
	\mathbb{O}\left(3\right):=\left\{ \left.A\in\mathbb{R}^{3\times3}\right|A^{\top}A=AA^{\top}=\mathbf{I}_{3}\right\} \label{eq:OVERVIEW_PER_SO3_1}
\end{equation}
$\mathbb{SO}\left(3\right)$ denotes the Special Orthogonal Group
and is a subgroup of the Orthogonal Group and the general linear group.
The attitude of a rigid body is defined as a rotational matrix $R$:
\begin{equation}
	\mathbb{SO}\left(3\right):=\left\{ \left.R\in\mathbb{R}^{3\times3}\right|R^{\top}R=RR^{\top}=\mathbf{I}_{3}\text{, }{\rm det}\left(R\right)=+1\right\} \label{eq:OVERVIEW_PER_SO3_2}
\end{equation}
where $\mathbf{I}_{3}$ is the identity matrix with $3$-dimensions
and ${\rm det\left(\cdot\right)}$ is the determinant of the associated
matrix. The associated Lie-algebra of $\mathbb{SO}\left(3\right)$
is termed $\mathfrak{so}\left(3\right)$ and is defined by 
\begin{equation}
	\mathfrak{so}\left(3\right):=\left\{ \left.\mathcal{A}=\left[\begin{array}{ccc}
		0 & -\alpha_{3} & \alpha_{2}\\
		\alpha_{3} & 0 & -\alpha_{1}\\
		-\alpha_{2} & \alpha_{1} & 0
	\end{array}\right]\right|\mathcal{A}^{\top}=-\mathcal{A}\right\} \label{eq:OVERVIEW_PER_SO3_3}
\end{equation}
with $\mathcal{A}\in\mathbb{R}^{3\times3}$ being the space of skew-symmetric
matrices. Define the map $\left[\cdot\right]_{\times}:\mathbb{R}^{3}\rightarrow\mathfrak{so}\left(3\right)$
such that 
\begin{equation}
	\mathcal{A}=\left[\alpha\right]_{\times}=\left[\begin{array}{ccc}
		0 & -\alpha_{3} & \alpha_{2}\\
		\alpha_{3} & 0 & -\alpha_{1}\\
		-\alpha_{2} & \alpha_{1} & 0
	\end{array}\right]\in\mathfrak{so}\left(3\right),\hspace{1em}\alpha=\left[\begin{array}{c}
		\alpha_{1}\\
		\alpha_{2}\\
		\alpha_{3}
	\end{array}\right]\in\mathbb{R}^{3}\label{eq:OVERVIEW_PER_SO3_4}
\end{equation}
For all $\alpha,\beta\in\mathbb{R}^{3}$, we have 
\[
\left[\alpha\right]_{\times}\beta=\alpha\times\beta
\]
where $\times$ is the cross product between the two vectors. Let
the vex operator be the inverse of a skew-symmetric matrix $\left[\cdot\right]_{\times}$
to vector form, denoted by $\mathbf{vex}:\mathfrak{so}\left(3\right)\rightarrow\mathbb{R}^{3}$
such that 
\begin{equation}
	\mathbf{vex}\left(\mathcal{A}\right)=\alpha\in\mathbb{R}^{3}\label{eq:OVERVIEW_PER_SO3_5}
\end{equation}
for all $\alpha\in\mathbb{R}^{3}$ and $\mathcal{A}\in\mathfrak{so}\left(3\right)$
as given in \eqref{eq:OVERVIEW_PER_SO3_4}. Let $\boldsymbol{\mathcal{P}}_{a}$
denote the anti-symmetric projection operator on the Lie-algebra $\mathfrak{so}\left(3\right)$,
defined by $\boldsymbol{\mathcal{P}}_{a}:\mathbb{R}^{3\times3}\rightarrow\mathfrak{so}\left(3\right)$
such that 
\begin{equation}
	\boldsymbol{\mathcal{P}}_{a}\left(\mathcal{B}\right)=\frac{1}{2}\left(\mathcal{B}-\mathcal{B}^{\top}\right)\in\mathfrak{so}\left(3\right)\label{eq:OVERVIEW_PER_SO3_PA_6}
\end{equation}
for all $\mathcal{B}\in\mathbb{R}^{3\times3}$. Hence, one obtains
\begin{align}
	\mathbf{vex}\left(\boldsymbol{\mathcal{P}}_{a}\left(\mathcal{B}\right)\right) & =\beta\in\mathbb{R}^{3}\label{eq:-36}\\
	\left[\beta\right]_{\times} & =\boldsymbol{\mathcal{P}}_{a}\left(\mathcal{B}\right)=\frac{1}{2}\left(\mathcal{B}-\mathcal{B}^{\top}\right)\label{eq:OVERVIEW_PER_SO3_VEX_7}
\end{align}

\subsection{Angle-axis parameterization preliminaries}

The attitude of a rigid body can be established given angle of rotation
$\alpha\in\mathbb{R}$ and a unit-axis $u=\left[u_{1},u_{2},u_{3}\right]^{\top}\in\mathbb{S}^{2}$.
This method of attitude reconstruction is termed angle-axis parameterization
\cite{shuster1993survey}. The mapping of angle-axis parameterization
to $\mathbb{SO}\left(3\right)$ is governed by $\mathcal{R}_{\alpha}:\mathbb{R}\times\mathbb{S}^{2}\rightarrow\mathbb{SO}\left(3\right)$
such that 
\begin{align}
	\mathcal{R}_{\alpha}\left(\alpha,u\right) & ={\rm exp}\left(\alpha\left[u\right]_{\times}\right)\nonumber \\
	& =\mathbf{I}_{3}+\sin\left(\alpha\right)\left[u\right]_{\times}+\left(1-\cos\left(\alpha\right)\right)\left[u\right]_{\times}^{2}\label{eq:SO3PPF_PER_AX_ANG}
\end{align}

\subsection{Rodriguez vector preliminaries}

The attitude can be obtained through Rodriguez parameters vector $\rho=\left[\rho_{1},\rho_{2},\rho_{3}\right]^{\top}\in\mathbb{R}^{3}$.
The related map to $\mathbb{SO}\left(3\right)$ is given by $\mathcal{R}_{\rho}:\mathbb{R}^{3}\rightarrow\mathbb{SO}\left(3\right)$
with 
\begin{align}
	\mathcal{R}_{\rho}= & \frac{1}{1+\left\Vert \rho\right\Vert ^{2}}\left(\left(1-\left\Vert \rho\right\Vert ^{2}\right)\mathbf{I}_{3}+2\rho\rho^{\top}+2\left[\rho\right]_{\times}\right)\label{eq:OVERVIEW_PER_ROD}
\end{align}

\subsection{Unit-quaternion preliminaries}

The unit-quaternion is defined by $Q=\left[q_{0},q^{\top}\right]^{\top}\in\mathbb{S}^{3}$
where $q_{0}\in\mathbb{R}$ and $q=\left[q_{1},q_{2},q_{3}\right]^{\top}\in\mathbb{R}^{3}$
such that 
\begin{equation}
	\mathbb{S}^{3}=\left\{ \left.Q\in\mathbb{R}^{4}\right|\left\Vert Q\right\Vert =1\right\} \label{eq:OVERVIEW_PER_Q_1}
\end{equation}
Let $Q=\left[q_{0},q^{\top}\right]^{\top}\in\mathbb{S}^{3}$, hence,
$Q^{*}=Q^{-1}\in\mathbb{S}^{3}$ can be defined as follows 
\begin{equation}
	Q^{*}=Q^{-1}=\left[\begin{array}{c}
		q_{0}\\
		-q
	\end{array}\right]\in\mathbb{S}^{3}\label{eq:OVERVIEW_PER_Q_2}
\end{equation}
where $Q^{*}$ and $Q^{-1}$ are complex conjugate and inverse of
the unit-quaternion, respectively.

\subsection{Useful math identities}

This subsection presents a list of identities which are going to prove
subsequently useful in the article. For any Lie bracket $\left[\mathcal{A},\mathcal{B}\right]$,
we have 
\begin{align}
	\left[\mathcal{A},\mathcal{B}\right] & =\mathcal{A}\mathcal{B}-\mathcal{B}\mathcal{A},\quad\mathcal{A},\mathcal{B}\in{\rm \mathbb{R}}^{3\times3}\label{eq:OVERVIEW_PER_Ident1}\\
	\left[\mathcal{A}\upsilon,\mathcal{B}\upsilon\right] & =\left(\mathcal{A}\mathcal{B}-\mathcal{B}\mathcal{A}\right)\upsilon,\quad\mathcal{A},\mathcal{B}\in{\rm \mathbb{R}}^{3\times3},\upsilon\in\mathbb{R}^{3}\label{eq:OVERVIEW_PER_Ident2}\\
	{\rm Tr}\left\{ \left[\mathcal{A},\mathcal{B}\right]\right\}  & =0,\quad\mathcal{A},\mathcal{B}\in{\rm \mathbb{R}}^{3\times3}\label{eq:OVERVIEW_PER_Ident3}
\end{align}
where ${\rm Tr}\left\{ \cdot\right\} $ refers to the trace of a matrix.
Let $R\in\mathbb{SO}\left(3\right)$ and $\upsilon,\omega\in{\rm \mathbb{R}}^{3}$,
then the following identities hold 
\begin{align}
	\left[R\upsilon\right]_{\times} & =R\left[\upsilon\right]_{\times}R^{\top}\label{eq:OVERVIEW_PER_Identity_used1}\\
	-\left[\omega\right]_{\times}\omega & =\left[0,0,0\right]^{\top}\label{eq:OVERVIEW_PER_Identity_used2}\\
	\left[\omega\right]_{\times}\upsilon & =-\left[\upsilon\right]_{\times}\omega\label{eq:OVERVIEW_PER_Identity_used3}\\
	-\left[\omega\right]_{\times}\left[\upsilon\right]_{\times} & =\left(\omega^{\top}\upsilon\right)\mathbf{I}_{3}-\upsilon\omega^{\top}\label{eq:OVERVIEW_PER_Identity_used4}\\
	{\rm Tr}\left\{ \mathcal{A}\left[\upsilon\right]_{\times}\right\}  & =0,\quad\mathcal{A}=\mathcal{A}^{\top}\in{\rm \mathbb{R}}^{3\times3}\label{eq:OVERVIEW_PER_Identity_used5}\\
	A^{\top}\left[\upsilon\right]_{\times}+\left[\upsilon\right]_{\times}A & =\left[\left({\rm Tr}\left\{ A\right\} \mathbf{I}_{3}-A\right)\upsilon\right]_{\times},\hspace{1em}A\in\mathbb{R}^{3\times3}\label{eq:OVERVIEW_PER_Identity_used6}\\
	{\rm Tr}\left\{ A\left[\upsilon\right]_{\times}\right\}  & ={\rm Tr}\left\{ \boldsymbol{\mathcal{P}}_{a}\left(A\right)\left[\upsilon\right]_{\times}\right\} \nonumber \\
	& =-2\mathbf{vex}\left(\boldsymbol{\mathcal{P}}_{a}\left(A\right)\right)^{\top}\upsilon,\quad A\in\mathbb{R}^{3\times3}\label{eq:OVERVIEW_PER_Identity_used7}
\end{align}
Before we proceed further, any $R$ is a rotational matrix and it
is important to note that 
\[
R=\mathcal{R}_{\xi}=\mathcal{R}_{\alpha}=\mathcal{R}_{\rho}=\mathcal{R}_{Q}\in\mathbb{SO}\left(3\right)
\]
with respect to the method of attitude parameterization.

\noindent %
\noindent\makebox[1\linewidth]{%
	\rule{0.6\textwidth}{1.4pt}%
}

\section{Special Orthogonal Group $\mathbb{SO}\left(3\right)$ \label{sec:OVERVIEW_SO3}}

The matrix $R\in\mathbb{R}^{3\times3}$ is said to represent the attitude
of a rigid-body if and only if $R\in\mathbb{SO}\left(3\right)$ which,
in turn, is true when the following two conditions are satisfied: 
\begin{enumerate}
	\item ${\rm det}\left(R\right)=+1$. 
	\item $R^{\top}R=RR^{\top}=\mathbf{I}_{3}$, which means that $R^{-1}=R^{\top}$. 
\end{enumerate}
$R$ is called rotation by inversion if ${\rm det}\left(R\right)=-1$,
which does not belong to $\mathbb{SO}\left(3\right)$ and therefore,
will not be considered in our analysis. The focus of this analysis
is $R\in\mathbb{SO}\left(3\right)$, which satisfies both of the above-mentioned
conditions. The main advantage of the $R\in\mathbb{SO}\left(3\right)$
representation is that the attitude $R$ is global and unique, implying
that each physical orientation of a rigid-body corresponds to a unique
rotational matrix.

\subsection{Attitude dynamics \label{subsec:OVERVIEW_SO3_Kinematics}}

\begin{figure}[h]
	\centering{}\includegraphics[scale=0.55]{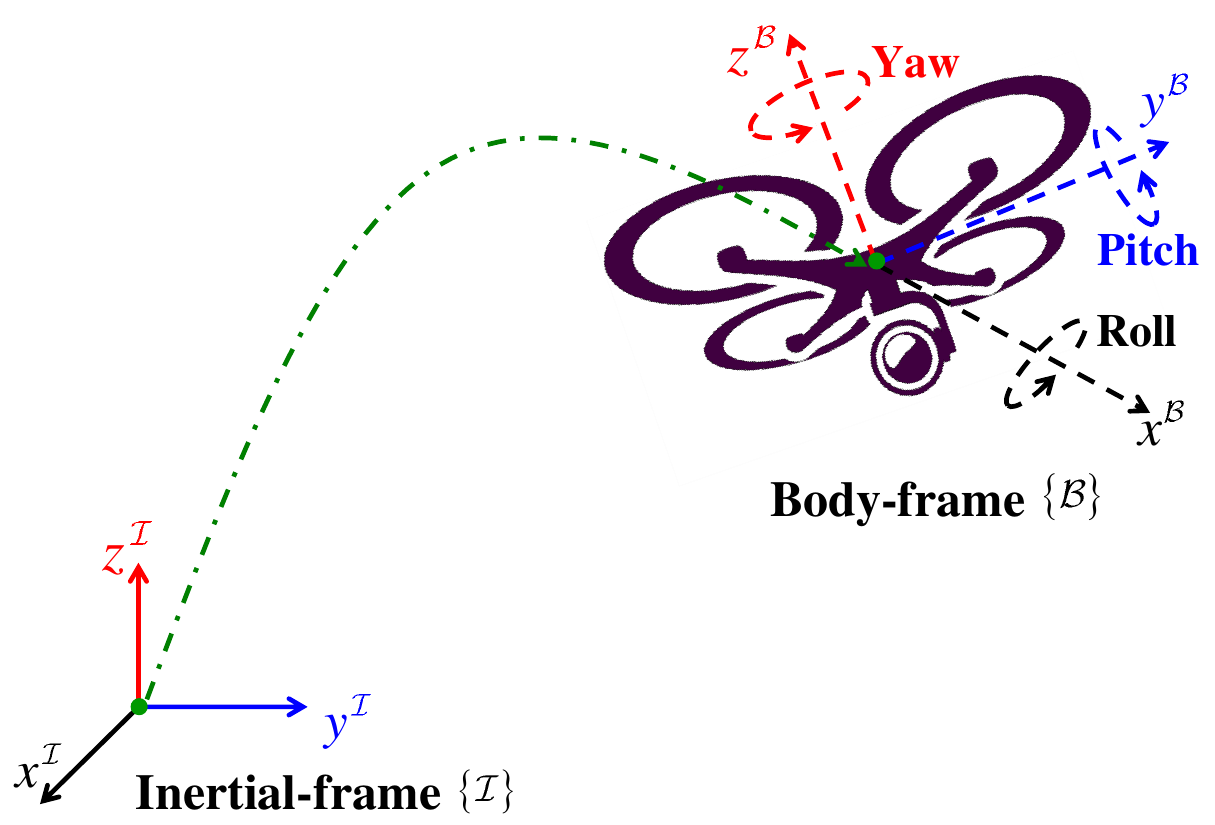}\caption{The orientation of a rigid-body in body-frame relative to inertial-frame
		in the 3D-space \cite{hashim2018SO3Stochastic}. }
	\label{fig:OVERVIEW_SO3_1} 
\end{figure}

Let $R\in\mathbb{SO}\left(3\right)$ denote the attitude (rotational
matrix). The relative orientation of the body-frame $\left\{ \mathcal{B}\right\} $
with respect to the inertial-frame $\left\{ \mathcal{I}\right\} $
is given by the rotational matrix $R$ as illustrated in Figure \ref{fig:OVERVIEW_SO3_1}.

The inertial-frame is fixed and commonly known as the world or the
global coordinate, while the origin of the inertial-frame is denoted
by $o^{\mathcal{I}}$. The body-frame is fixed to the moving rigid-body
and the origin of the body-frame is denoted by $o^{\mathcal{B}}$.
For a body-frame vector ${\rm v}^{\mathcal{B}}\in\mathbb{R}^{3}$
and inertial-frame vector ${\rm v}^{\mathcal{I}}\in\mathbb{R}^{3}$,
the translation from body-frame to inertial-frame is defined by the
attitude matrix $R\in\mathbb{SO}\left(3\right)$

\begin{equation}
	{\rm v}^{\mathcal{B}}=R^{\top}{\rm v}^{\mathcal{I}}\label{eq:OVERVIEW_SO3_dyn1}
\end{equation}
The orientation of the inertial-frame is always fixed which means
that $\dot{{\rm v}}^{\mathcal{I}}=0$. The velocity at the body-frame
in \eqref{eq:OVERVIEW_SO3_dyn1} is expressed by 
\begin{align}
	\dot{{\rm v}}^{\mathcal{B}} & =\dot{R}^{\top}{\rm v}^{\mathcal{I}}\nonumber \\
	& =\dot{R}^{\top}R{\rm v}^{\mathcal{B}}\label{eq:OVERVIEW_SO3_dyn2}
\end{align}
At every time instant $t$, there exists a unique angular velocity
vector $\Omega=\left[\Omega_{x},\Omega_{y},\Omega_{z}\right]^{\top}\in\mathbb{R}^{3}$
such that, for every particle in the body, one has 
\begin{align}
	\dot{{\rm v}}^{\mathcal{B}} & ={\rm v}^{\mathcal{B}}\times\Omega\nonumber \\
	& =\left[{\rm v}^{\mathcal{B}}\right]_{\times}\Omega\label{eq:OVERVIEW_SO3_dyn3}
\end{align}
where $\Omega$ is the body-referenced angular velocity. Hence, combining
\eqref{eq:OVERVIEW_SO3_dyn2} and \eqref{eq:OVERVIEW_SO3_dyn3} yields
\begin{align}
	\dot{R}^{\top}R{\rm v}^{\mathcal{B}} & =-\left[\Omega\right]_{\times}{\rm v}^{\mathcal{B}}\nonumber \\
	\dot{R}^{\top} & =\left[\Omega\right]_{\times}^{\top}R^{\top}\label{eq:OVERVIEW_SO3_dyn4}
\end{align}
where $\left[{\rm v}^{\mathcal{B}}\right]_{\times}\Omega=-\left[\Omega\right]_{\times}{\rm v}^{\mathcal{B}}$.
Thus, the attitude dynamics $\dot{R}\in\mathbb{R}^{3\times3}$ of
a rigid-body are given by

\begin{align}
	\dot{R} & =R\left[\Omega\right]_{\times}\label{eq:OVERVIEW_SO3_dyn}
\end{align}
From \eqref{eq:OVERVIEW_SO3_dyn}, one can find that the anti-symmetric
projection operator of the attitude kinematics $\boldsymbol{\mathcal{P}}_{a}:\mathbb{R}^{3\times3}\rightarrow\mathfrak{so}\left(3\right)$
is 
\begin{align}
	\boldsymbol{\mathcal{P}}_{a}\left(\dot{R}\right)= & \boldsymbol{\mathcal{P}}_{a}\left(R\left[\Omega\right]_{\times}\right)\nonumber \\
	= & \frac{1}{2}\left(R\left[\Omega\right]_{\times}+\left[\Omega\right]_{\times}R^{\top}\right)\nonumber \\
	= & \frac{1}{2}\left({\rm Tr}\left\{ R\right\} \left[\Omega\right]_{\times}-\left[R^{\top}\Omega\right]_{\times}\right)\nonumber \\
	= & \frac{1}{2}\left[\left({\rm Tr}\left\{ R\right\} \mathbf{I}_{3}-R\right)^{\top}\Omega\right]_{\times}\in\mathfrak{so}\left(3\right)\label{eq:-37}
\end{align}
where $R\left[\Omega\right]_{\times}+\left[\Omega\right]_{\times}R^{\top}=\left[\left({\rm Tr}\left\{ R\right\} \mathbf{I}_{3}-R\right)^{\top}\Omega\right]_{\times}$
as defined in identity \eqref{eq:OVERVIEW_PER_Identity_used6}. Thus,
the dynamics in \eqref{eq:OVERVIEW_SO3_dyn} can be defined in vector
form $\mathbf{vex}\left(\boldsymbol{\mathcal{P}}_{a}\left(R\right)\right)\in\mathbb{R}^{3}$
as 
\begin{align}
	\frac{d}{dt}\mathbf{vex}\left(\boldsymbol{\mathcal{P}}_{a}\left(R\right)\right)= & \mathbf{vex}\left(\boldsymbol{\mathcal{P}}_{a}\left(R\left[\Omega\right]_{\times}\right)\right)\nonumber \\
	= & \frac{1}{2}\left({\rm Tr}\left\{ R\right\} \mathbf{I}_{3}-R\right)^{\top}\Omega\label{eq:OVERVIEW_SO3_dyn_VEX}
\end{align}
Similarly, for $\mathcal{A}=\mathcal{A}^{\top}\in\mathbb{R}^{3\times3}$,
one obtains 
\begin{align}
	\frac{d}{dt}\mathbf{vex}\left(\boldsymbol{\mathcal{P}}_{a}\left(\mathcal{A}R\right)\right)= & \mathbf{vex}\left(\boldsymbol{\mathcal{P}}_{a}\left(\mathcal{A}R\left[\Omega\right]_{\times}\right)\right)\nonumber \\
	= & \frac{1}{2}\left({\rm Tr}\left\{ \mathcal{A}R\right\} \mathbf{I}_{3}-\mathcal{A}R\right)^{\top}\Omega\label{eq:OVERVIEW_SO3_dyn_VEX2}
\end{align}

\subsection{Normalized Euclidean distance}

Since, $-1\leq{\rm Tr}\left\{ R\right\} \leq3$ for $R\in\mathbb{SO}\left(3\right)$,
the normalized Euclidean distance of a rotational matrix on $\mathbb{SO}\left(3\right)$
is defined as 
\begin{equation}
	\left\Vert R\right\Vert _{I}=\frac{1}{4}{\rm Tr}\left\{ \mathbf{I}_{3}-R\right\} \in\left[0,1\right]\label{eq:OVERVIEW_SO3_Ecul_Dist}
\end{equation}

\subsection{Normalized Euclidean distance dynamics}

From the normalized Euclidean distance in \eqref{eq:OVERVIEW_SO3_Ecul_Dist},
the attitude dynamics in \eqref{eq:OVERVIEW_SO3_dyn}, and the identity
in \eqref{eq:OVERVIEW_PER_Identity_used7}, one can find that the
dynamics of the normalized Euclidean distance of \eqref{eq:OVERVIEW_SO3_Ecul_Dist}
are equal to 
\begin{align}
	\frac{d}{dt}\left\Vert R\right\Vert _{I} & =\frac{d}{dt}\frac{1}{4}{\rm Tr}\left\{ \mathbf{I}_{3}-R\right\} \nonumber \\
	& =-\frac{1}{4}{\rm Tr}\left\{ \dot{R}\right\} \nonumber \\
	& =-\frac{1}{4}{\rm Tr}\left\{ R\left[\Omega\right]_{\times}\right\} \nonumber \\
	& =-\frac{1}{4}{\rm Tr}\left\{ R\boldsymbol{\mathcal{P}}_{a}\left(\left[\Omega\right]_{\times}\right)\right\} \nonumber \\
	& =\frac{1}{2}\mathbf{vex}\left(\boldsymbol{\mathcal{P}}_{a}\left(R\right)\right)^{\top}\Omega\label{eq:OVERVIEW_SO3_dyn_norm}
\end{align}

\subsubsection{Discrete Attitude dynamics}

The attitude dynamics defined in continuous form in \eqref{eq:OVERVIEW_SO3_dyn}
could be expressed in discrete form using exact integration by 
\begin{equation}
	R\left[k+1\right]=R\left[k\right]\exp\left(\left[\Omega\left[k\right]\right]_{\times}\Delta t\right)\label{eq:OVERVIEW_SO3_dyn_discrete}
\end{equation}
where $k\in\mathbb{N}$ refers to the $k$th sample, $\Delta t$ is
a time step which is normally small, and $R\left[k\right]$ and $\Omega\left[k\right]$
refer to the true attitude and angular velocity at the $k$th sample,
respectively.

Now we introduce important properties which will be helpful in the
process of designing attitude filter and/or attitude control. \begin{lem}
	\label{Lem:OVERVIEW_6} Let $R\in\mathbb{SO}\left(3\right)$, $u\in\mathbb{S}^{2}$
	denote a unit vector such that $\left\Vert u\right\Vert =1$, and
	$\alpha\in\mathbb{R}$ denote the angle of rotation about $u$. Then,
	the following holds: 
	\begin{align}
		\mathbf{vex}\left(\boldsymbol{\mathcal{P}}_{a}\left(R\right)\right) & =2{\rm cos}\left(\frac{\alpha}{2}\right)\sin\left(\frac{\alpha}{2}\right)u\label{eq:OVERVIEW_SO3_lem6_1}\\
		\left\Vert R\right\Vert _{I} & =\frac{1}{2}\left(1-\cos\left(\alpha\right)\right)=\sin^{2}\left(\frac{\alpha}{2}\right)\label{eq:OVERVIEW_SO3_lem6_2}\\
		\left\Vert \mathbf{vex}\left(\boldsymbol{\mathcal{P}}_{a}\left(R\right)\right)\right\Vert ^{2} & =4{\rm cos}^{2}\left(\frac{\alpha}{2}\right)\sin^{2}\left(\frac{\alpha}{2}\right)\label{eq:OVERVIEW_SO3_lem6_3}
	\end{align}
	\textbf{Proof. }See\textbf{ \eqref{eq:OVERVIEW_att_ang_VEX}},\textbf{
		\eqref{eq:OVERVIEW_att_ang_RI}}, and\textbf{ \eqref{eq:OVERVIEW_att_ang_Vex1}.
} \end{lem} \begin{lem} \label{Lem:OVERVIEW_1} Let $R\in\mathbb{SO}\left(3\right)$
	and $\rho\in\mathbb{R}^{3}$ be the Rodriguez parameters vector. Then,
	the following holds: 
	\begin{align}
		\mathbf{vex}\left(\boldsymbol{\mathcal{P}}_{a}\left(R\right)\right) & =2\frac{\rho}{1+\left\Vert \rho\right\Vert ^{2}}\label{eq:OVERVIEW_SO3_lem1_1}\\
		\left\Vert R\right\Vert _{I} & =\frac{\left\Vert \rho\right\Vert ^{2}}{1+\left\Vert \rho\right\Vert ^{2}}\label{eq:OVERVIEW_SO3_lem1_2}\\
		\left\Vert \mathbf{vex}\left(\boldsymbol{\mathcal{P}}_{a}\left(R\right)\right)\right\Vert ^{2} & =4\frac{\left\Vert \rho\right\Vert ^{2}}{\left(1+\left\Vert \rho\right\Vert ^{2}\right)^{2}}\label{eq:OVERVIEW_SO3_lem1_4}
	\end{align}
	\textbf{Proof. }See\textbf{ \eqref{eq:OVERVIEW_ROD_VEX_Pa}},\textbf{
		\eqref{eq:OVERVIEW_ROD_TR2}},\textbf{ }and\textbf{ \eqref{eq:OVERVIEW_ROD_VEX2_1}.
} \end{lem} The following Lemma (Lemma \ref{Lem:OVERVIEW_5}) is
true if either Lemma \ref{Lem:OVERVIEW_6} or Lemma \ref{Lem:OVERVIEW_1}
holds: \begin{lem} \label{Lem:OVERVIEW_5} Let $R\in\mathbb{SO}\left(3\right)$.
	Then, the following holds: 
	\begin{align}
		\left\Vert \mathbf{vex}\left(\boldsymbol{\mathcal{P}}_{a}\left(R\right)\right)\right\Vert ^{2} & =4\left(1-\left\Vert R\right\Vert _{I}\right)\left\Vert R\right\Vert _{I}\label{eq:OVERVIEW_SO3_lem5_1}
	\end{align}
	\textbf{Proof. }See\textbf{ \eqref{eq:OVERVIEW_att_ang_Vex2}}, or\textbf{
		\eqref{eq:OVERVIEW_ROD_VEX2_1}. } \end{lem} The following Lemma
(Lemma \ref{Lem:OVERVIEW_7}) is true if both Lemma \ref{Lem:OVERVIEW_6}
and Lemma \ref{Lem:OVERVIEW_1} hold: \begin{lem} \label{Lem:OVERVIEW_7}
	Let $\rho\in\mathbb{R}^{3}$ be the Rodriguez parameters vector, $u\in\mathbb{S}^{2}$
	denote a unit vector such that $\left\Vert u\right\Vert =1$, and
	$\alpha\in\mathbb{R}$ denote the angle of rotation about $u$. Then,
	the following holds: 
	\begin{align}
		\rho & ={\rm tan}\left(\frac{\alpha}{2}\right)u\label{eq:OVERVIEW_SO3_lem7_1}\\
		\alpha & =2\text{ }{\rm tan}^{-1}\left(\left\Vert \rho\right\Vert \right)\label{eq:OVERVIEW_SO3_lem7_2}\\
		u & ={\rm cot}\left(\frac{\alpha}{2}\right)\rho\label{eq:OVERVIEW_SO3_lem7_3}
	\end{align}
	\textbf{Proof. }See\textbf{ \eqref{eq:OVERVIEW_ROD_from_a_u}},\textbf{
		\eqref{eq:OVERVIEW_ROD_2_a_u1}}, and\textbf{ \eqref{eq:OVERVIEW_ROD_2_a_u2}.
} \end{lem} \begin{lem} \label{Lem:OVERVIEW_2} Let $R\in\mathbb{SO}\left(3\right)$
	and $Q=\left[q_{0},q^{\top}\right]^{\top}\in\mathbb{S}^{3}$ be a
	unit-quaternion vector with $q_{0}\in\mathbb{R}$ and $q\in\mathbb{R}^{3}$.
	Then, the following holds: 
	\begin{align}
		\mathbf{vex}\left(\boldsymbol{\mathcal{P}}_{a}\left(R\right)\right) & =2q_{0}q\label{eq:OVERVIEW_SO3_lem2_1}\\
		\left\Vert R\right\Vert _{I} & =1-q_{0}^{2}\label{eq:OVERVIEW_SO3_lem2_2}\\
		\left\Vert \mathbf{vex}\left(\boldsymbol{\mathcal{P}}_{a}\left(R\right)\right)\right\Vert ^{2} & =4q_{0}^{2}\left\Vert q\right\Vert ^{2}\label{eq:OVERVIEW_SO3_lem2_3}
	\end{align}
	\textbf{Proof. }See\textbf{ \eqref{eq:OVERVIEW_Q_TR2}},\textbf{ \eqref{eq:OVERVIEW_Q_VEX}},
	and\textbf{ \eqref{eq:OVERVIEW_Q_VEX2}. } \end{lem} The following
Lemma (Lemma \ref{Lem:OVERVIEW_8}) is true if both Lemma \ref{Lem:OVERVIEW_6}
and Lemma \ref{Lem:OVERVIEW_2} hold: \begin{lem} \label{Lem:OVERVIEW_8}
	Let $Q=\left[q_{0},q^{\top}\right]^{\top}\in\mathbb{S}^{3}$ be a
	unit-quaternion vector with $q_{0}\in\mathbb{R}$ and $q\in\mathbb{R}^{3}$,
	and let $u\in\mathbb{S}^{2}$ denote a unit vector such that $\left\Vert u\right\Vert =1$,
	and $\alpha\in\mathbb{R}$ denote the angle of rotation about $u$.
	Then, the following holds: 
	\begin{align}
		\alpha & =2\cos^{-1}\left(q_{0}\right)\label{eq:OVERVIEW_SO3_lem8_1}\\
		u & =\frac{1}{\sin\left(\alpha/2\right)}q\label{eq:OVERVIEW_SO3_lem8_2}\\
		q_{0} & =\cos\left(\alpha/2\right)\label{eq:OVERVIEW_SO3_lem8_3}\\
		q & =u\sin\left(\alpha/2\right)\label{eq:OVERVIEW_SO3_lem8_4}
	\end{align}
	\textbf{Proof. }See\textbf{ \eqref{eq:OVERVIEW_Q_Q_2_ang1}}, \textbf{\eqref{eq:OVERVIEW_Q_Q_2_ang2}},\textbf{
	}and\textbf{ \eqref{eq:OVERVIEW_Q_ang_2_Q}. } \end{lem} The following
Lemma (Lemma \ref{Lem:OVERVIEW_4}) is true if Lemma \ref{Lem:OVERVIEW_1}
and Lemma \ref{Lem:OVERVIEW_2} hold: \begin{lem} \label{Lem:OVERVIEW_4}
	Let $\rho\in\mathbb{R}^{3}$ be the Rodriguez parameters vector and
	let $Q=\left[q_{0},q^{\top}\right]^{\top}\in\mathbb{S}^{3}$ be a
	unit-quaternion vector with $q_{0}\in\mathbb{R}$ and $q\in\mathbb{R}^{3}$.
	Then, the following holds: 
	\begin{align}
		q_{0} & =\pm\frac{1}{\sqrt{1+\left\Vert \rho\right\Vert ^{2}}}\label{eq:OVERVIEW_SO3_lem4_1}\\
		q & =\pm\frac{\rho}{\sqrt{1+\left\Vert \rho\right\Vert ^{2}}}\label{eq:OVERVIEW_SO3_lem4_2}\\
		\rho & =\frac{q}{q_{0}}\label{eq:OVERVIEW_SO3_lem4_3}
	\end{align}
	\textbf{Proof. }See\textbf{ }Lemma \ref{Lem:OVERVIEW_1} and Lemma
	\ref{Lem:OVERVIEW_2}.\textbf{ }It should be remarked that Equation
	\eqref{eq:OVERVIEW_SO3_lem4_3} is only valid for $q_{0}\neq0$. \end{lem}
\begin{lem} \label{Lem:OVERVIEW_3} Let $R\in\mathbb{SO}\left(3\right)$
	and $\rho\in\mathbb{R}^{3}$ be the Rodriguez parameters vector. Define
	$M=M^{\top}\in\mathbb{R}^{3\times3}$, where $M$ is with rank 3,
	and ${\rm Tr}\left\{ M\right\} =3$. Define $\bar{\mathbf{M}}={\rm Tr}\left\{ M\right\} \mathbf{I}_{3}-M$
	and let the minimum singular value of $\bar{\mathbf{M}}$ be $\underline{\lambda}:=\underline{\lambda}\left(\bar{\mathbf{M}}\right)$.
	Then, the following holds: 
	\begin{align}
		\left\Vert MR\right\Vert _{I} & =\frac{1}{2}\frac{\rho^{\top}\bar{\mathbf{M}}\rho}{1+\left\Vert \rho\right\Vert ^{2}}\label{eq:OVERVIEW_SO3_lem3_1}\\
		\mathbf{vex}\left(\boldsymbol{\mathcal{P}}_{a}\left(MR\right)\right) & =\frac{\left(\mathbf{I}_{3}+\left[\rho\right]_{\times}\right)^{\top}\bar{\mathbf{M}}}{1+\left\Vert \rho\right\Vert ^{2}}\rho\label{eq:OVERVIEW_SO3_lem3_2}\\
		\left\Vert \mathcal{\mathbf{vex}}\left(\boldsymbol{\mathcal{P}}_{a}\left(MR\right)\right)\right\Vert ^{2} & =\frac{\rho^{\top}\bar{\mathbf{M}}\left(\mathbf{I}_{3}-\left[\rho\right]_{\times}^{2}\right)\bar{\mathbf{M}}\rho}{\left(1+\left\Vert \rho\right\Vert ^{2}\right)^{2}}\label{eq:OVERVIEW_SO3_lem3_3}\\
		\left\Vert MR\right\Vert _{I} & \leq\frac{2}{\underline{\lambda}}\frac{\left\Vert \mathcal{\mathbf{vex}}\left(\boldsymbol{\mathcal{P}}_{a}\left(MR\right)\right)\right\Vert ^{2}}{1+{\rm Tr}\left\{ M^{-1}MR\right\} }\label{eq:OVERVIEW_SO3_lem3_4}
	\end{align}
	\textbf{Proof. }See\textbf{ \eqref{eq:OVERVIEW_ROD_EXPL_append_MBR_I}},\textbf{
		\eqref{eq:OVERVIEW_ROD_EXPL_append_MBR_VEX}},\textbf{ \eqref{eq:OVERVIEW_ROD_EXPL_append_MBR_VEX2}},
	and\textbf{ \eqref{eq:OVERVIEW_ROD_EXPL_Proof_VEX_MI2}. } \end{lem}
Some of the above mentioned results and mapping to/from $\mathbb{SO}\left(3\right)$
from/to other attitude representations could be found in \cite{stuelpnagel1964parametrization},
\cite{wen1991attitude}, \cite{shuster1993survey}, \cite{hashim2018SO3Stochastic},
\cite{hashim2019SE3Det}, and \cite{mohamed2019filters}.

\subsection{Attitude error and attitude error dynamics}

Consider the attitude dynamics in \textbf{\eqref{eq:OVERVIEW_SO3_dyn}}
\[
\dot{R}=R\left[\Omega\right]_{\times}
\]
Let us introduce desired/estimator attitude dynamics 
\begin{equation}
	\dot{R}_{\star}=R_{\star}\left[\Omega_{\star}\right]_{\times}\label{eq:OVERVIEW_SO3_dyn_EST}
\end{equation}
Consider the error in attitude to be given by 
\begin{equation}
	\tilde{R}=R^{\top}R_{\star}\label{eq:OVERVIEW_SO3_Error}
\end{equation}
The attitude filter/control aims to drive $\tilde{R}\rightarrow\mathbf{I}_{3}$.
The dynamics of the attitude error can be found to be 
\begin{align}
	\dot{\tilde{R}} & =\dot{R}^{\top}R_{\star}+R^{\top}\dot{R}_{\star}\nonumber \\
	& =\left[\Omega\right]_{\times}^{\top}R^{\top}R_{\star}+R^{\top}R_{\star}\left[\Omega_{\star}\right]_{\times}\nonumber \\
	& =-\left[\Omega\right]_{\times}\tilde{R}+\tilde{R}\left[\Omega_{\star}\right]_{\times}\label{eq:OVERVIEW_SO3_dot_Error}
\end{align}
where $\left[\Omega\right]_{\times}^{\top}=-\left[\Omega\right]_{\times}$
as given in \textbf{\eqref{eq:OVERVIEW_PER_SO3_3}}.

\noindent %
\noindent\makebox[1\linewidth]{%
	\rule{0.6\textwidth}{1.4pt}%
}

\section{Euler Parameterization \label{sec:Euler-Parameterization}}

\begin{figure}[h]
	\centering{}\includegraphics[scale=0.5]{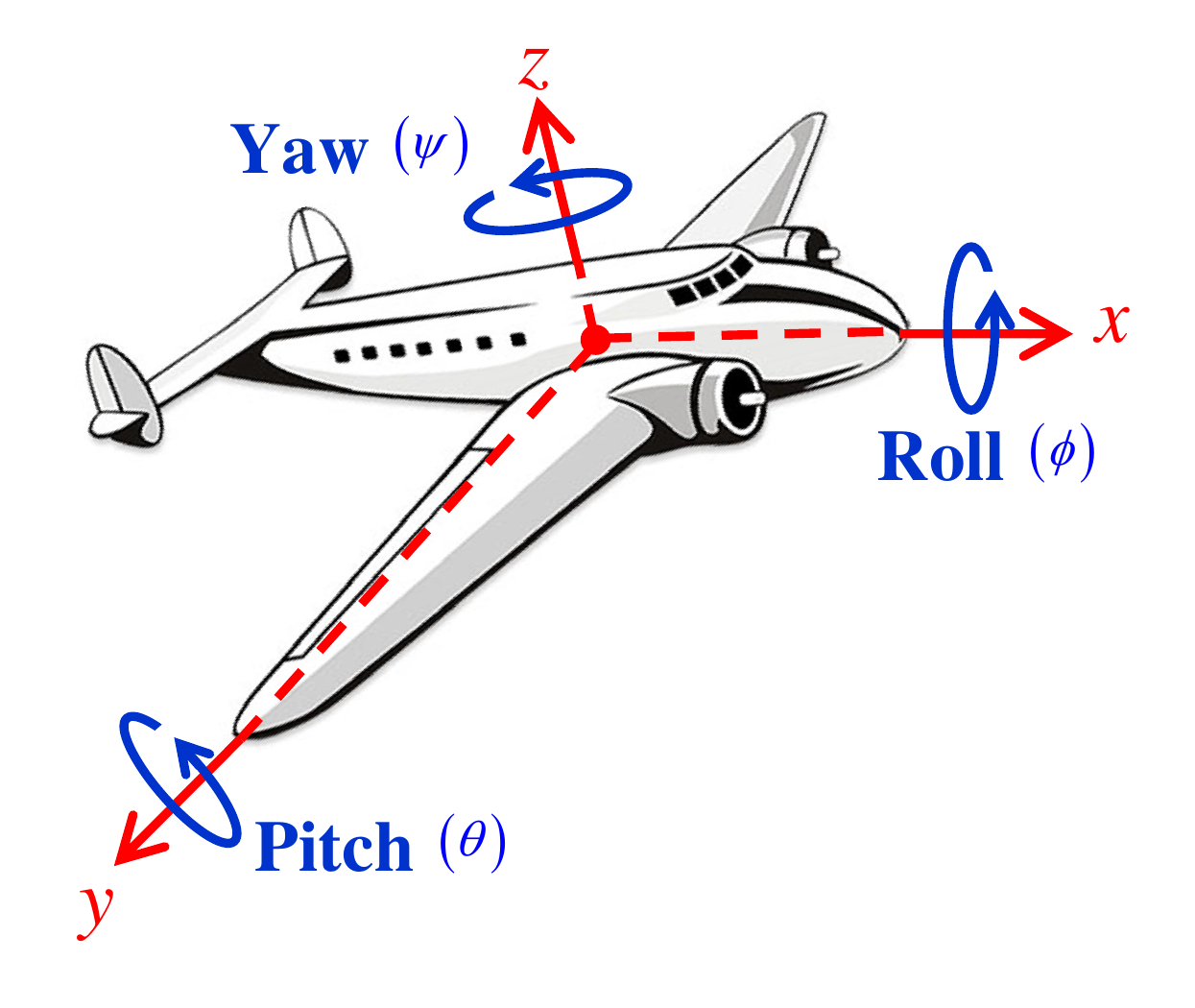}\caption{Graphical representation of Euler angles with respect to the reference-axis
		of the body-frame. }
	\label{fig:OVERVIEW_SO3_2} 
\end{figure}

The set of Euler angles is extensively used and widely known for attitude
representation of the rigid-body. These angles are easy to visualize
and understand allowing many researchers to use Euler angles for attitude
parameterization such as \cite{grood1983joint,mokhtari2004dynamic}.
The naming convention is to label the angle based on the axis of rotation
associated with. Roll angle represents rotation about the $x$-axis
and is denoted by $\phi$, pitch angle refers to rotation about the
$y$-axis and is denoted by $\theta$, and yaw angle is the rotation
about the $z$-axis and is denoted by $\psi$. Figure \ref{fig:OVERVIEW_SO3_2}
illustrates the orientation of the rigid-body relative to the three
axes and Euler angles associated with each axis of rotation.

The names of these angles: roll, pitch and yaw are widely used in
all aircraft applications. The rotational matrix could be described
by a direction cosine matrix \cite{spong2008robot}. Let the world
axis of the inertial-frame be $p^{\mathcal{I}}=\left[x^{\mathcal{I}},y^{\mathcal{I}},z^{\mathcal{I}}\right]^{\top}$
and the body axis of the rigid-body fixed frame be $p^{\mathcal{B}}=\left[x^{\mathcal{B}},y^{\mathcal{B}},z^{\mathcal{B}}\right]^{\top}$.
Then, the direction cosine matrix of an angle $\gamma$ is given by
\begin{equation}
	R=\left[\begin{array}{ccc}
		\cos\left(\gamma_{x^{\mathcal{B}},x^{\mathcal{I}}}\right) & \cos\left(\gamma_{x^{\mathcal{B}},y^{\mathcal{I}}}\right) & \cos\left(\gamma_{x^{\mathcal{B}},z^{\mathcal{I}}}\right)\\
		\cos\left(\gamma_{y^{\mathcal{B}},x^{\mathcal{I}}}\right) & \cos\left(\gamma_{y^{\mathcal{B}},y^{\mathcal{I}}}\right) & \cos\left(\gamma_{y^{\mathcal{B}},z^{\mathcal{I}}}\right)\\
		\cos\left(\gamma_{z^{\mathcal{B}},x^{\mathcal{I}}}\right) & \cos\left(\gamma_{z^{\mathcal{B}},y^{\mathcal{I}}}\right) & \cos\left(\gamma_{z^{\mathcal{B}},z^{\mathcal{I}}}\right)
	\end{array}\right]\label{eq:-35}
\end{equation}
For instance, consider the cosine matrix in \eqref{eq:-35}, the rotation
by an angle $\phi$ about the $x^{\mathcal{I}}$-axis is given by
the following cosine matrix 
\[
R=\left[\begin{array}{ccc}
	\cos\left(0\right) & \cos\left(\pi/2\right) & \cos\left(\pi/2\right)\\
	\cos\left(\pi/2\right) & \cos\left(\phi\right) & \cos\left(\pi/2+\phi\right)\\
	\cos\left(\pi/2\right) & \cos\left(\pi/2-\phi\right) & \cos\left(\phi\right)
\end{array}\right]=\left[\begin{array}{ccc}
	1 & 0 & 0\\
	0 & \cos\left(\phi\right) & -\sin\left(\phi\right)\\
	0 & \sin\left(\phi\right) & \cos\left(\phi\right)
\end{array}\right]
\]
More details on the direction cosine matrix visit \cite{spong2008robot}.

\subsection{Attitude parameterization through Euler angles}

\noindent Take, for example, the inertial-frame and the body-frame
as defined in Figure \ref{fig:OVERVIEW_SO3_2}. For simplicity, consider
the following notation 
\[
s=\sin\text{, }c=\cos\text{, }t=\tan
\]
With respect to the cosine matrix in \eqref{eq:-35}, the rotation
about $x^{\mathcal{I}}$, $y^{\mathcal{I}}$, and $z^{\mathcal{I}}$,
respectively, is given by 
\begin{itemize}
	\item \textit{rolling }about $x^{\mathcal{I}}$ by an angle $\phi$ 
	\begin{equation}
		\left[\begin{array}{c}
			x^{\mathcal{I}}\\
			y^{\mathcal{I}}\\
			z^{\mathcal{I}}
		\end{array}\right]=\left[\begin{array}{ccc}
			1 & 0 & 0\\
			0 & c\phi & -s\phi\\
			0 & s\phi & c\phi
		\end{array}\right]\left[\begin{array}{c}
			x^{\mathcal{B}}\\
			y^{\mathcal{B}}\\
			z^{\mathcal{B}}
		\end{array}\right]\Rightarrow R_{x,\phi}=\left[\begin{array}{ccc}
			1 & 0 & 0\\
			0 & c\phi & -s\phi\\
			0 & s\phi & c\phi
		\end{array}\right]\label{eq:-13}
	\end{equation}
	\item \textit{pitching }about $y^{\mathcal{I}}$ by an angle $\theta$ 
	\begin{equation}
		\left[\begin{array}{c}
			x^{\mathcal{I}}\\
			y^{\mathcal{I}}\\
			z^{\mathcal{I}}
		\end{array}\right]=\left[\begin{array}{ccc}
			c\theta & 0 & s\theta\\
			0 & 1 & 0\\
			-s\theta & 0 & c\theta
		\end{array}\right]\left[\begin{array}{c}
			x^{\mathcal{B}}\\
			y^{\mathcal{B}}\\
			z^{\mathcal{B}}
		\end{array}\right]\Rightarrow R_{y,\theta}=\left[\begin{array}{ccc}
			c\theta & 0 & s\theta\\
			0 & 1 & 0\\
			-s\theta & 0 & c\theta
		\end{array}\right]\label{eq:-12}
	\end{equation}
	\item \textit{yawing }about $z^{\mathcal{I}}$ by an angle $\psi$ 
	\begin{equation}
		\left[\begin{array}{c}
			x^{\mathcal{I}}\\
			y^{\mathcal{I}}\\
			z^{\mathcal{I}}
		\end{array}\right]=\left[\begin{array}{ccc}
			c\psi & s\psi & 0\\
			-s\psi & c\psi & 0\\
			0 & 0 & 1
		\end{array}\right]\left[\begin{array}{c}
			x^{\mathcal{B}}\\
			y^{\mathcal{B}}\\
			z^{\mathcal{B}}
		\end{array}\right]\Rightarrow R_{z,\psi}=\left[\begin{array}{ccc}
			c\psi & -s\psi & 0\\
			s\psi & c\psi & 0\\
			0 & 0 & 1
		\end{array}\right]\label{eq:-11}
	\end{equation}
\end{itemize}
The transformation matrix is obtained by 
\begin{enumerate}
	\item Transformation about the fixed $\left(x,\phi\right)$. 
	\item Next, transformation about the fixed $\left(y,\theta\right)$. 
	\item Next, transformation about the fixed $\left(z,\psi\right)$. 
\end{enumerate}
to yield 
\begin{align*}
	\mathcal{R}_{\xi} & =R_{z,\psi}R_{y,\theta}R_{x,\phi}=R_{\left(\phi,\theta,\psi\right)}\\
	& =\left[\begin{array}{ccc}
		c\psi & -s\psi & 0\\
		s\psi & c\psi & 0\\
		0 & 0 & 1
	\end{array}\right]\left[\begin{array}{ccc}
		c\theta & 0 & s\theta\\
		0 & 1 & 0\\
		-s\theta & 0 & c\theta
	\end{array}\right]\left[\begin{array}{ccc}
		1 & 0 & 0\\
		0 & c\phi & -s\phi\\
		0 & s\phi & c\phi
	\end{array}\right]
\end{align*}
which is equivalent to 
\begin{equation}
	\mathcal{R}_{\xi}=\left[\begin{array}{ccc}
		c\theta c\psi & -c\phi s\psi+s\phi s\theta c\psi & s\phi s\psi+c\phi s\theta c\psi\\
		c\theta s\psi & c\phi c\psi+s\phi s\theta s\psi & -s\phi c\psi+c\phi s\theta s\psi\\
		-s\theta & s\phi c\theta & c\phi c\theta
	\end{array}\right]\in\mathbb{SO}\left(3\right)\label{eq:OVERVIEW_EUL_R}
\end{equation}
with $\mathcal{R}_{\xi}\in\mathbb{SO}\left(3\right)$ being the attitude
representation using Euler angles $\left(\phi,\theta,\psi\right)$.
Thus, the relation between the inertial-frame and body-frame is equivalent
to 
\begin{align*}
	\left[\begin{array}{c}
		x^{\mathcal{I}}\\
		y^{\mathcal{I}}\\
		z^{\mathcal{I}}
	\end{array}\right] & =\left[\begin{array}{ccc}
		c\theta c\psi & -c\phi s\psi+s\phi s\theta c\psi & s\phi s\psi+c\phi s\theta c\psi\\
		c\theta s\psi & c\phi c\psi+s\phi s\theta s\psi & -s\phi c\psi+c\phi s\theta s\psi\\
		-s\theta & s\phi c\theta & c\phi c\theta
	\end{array}\right]\left[\begin{array}{c}
		x^{\mathcal{B}}\\
		y^{\mathcal{B}}\\
		z^{\mathcal{B}}
	\end{array}\right]\\
	{\rm v}^{\mathcal{I}} & =\mathcal{R}_{\xi}{\rm v}^{\mathcal{B}}
\end{align*}

\subsection{Euler angles from $\mathbb{SO}\left(3\right)$}

\noindent In order to obtain Euler angles from a given rotational
matrix, from \eqref{eq:OVERVIEW_EUL_R} one has 
\begin{equation}
	\frac{R_{\left(3,2\right)}}{R_{\left(3,3\right)}}=\frac{\sin(\phi)\cos\left(\theta\right)}{\cos(\phi)\cos\left(\theta\right)}\label{eq:-14}
\end{equation}
Hence, the roll angle can be found to be 
\begin{equation}
	\phi=\arctan\left(\frac{R_{\left(3,2\right)}}{R_{\left(3,3\right)}}\right)\label{eq:OVERVIEW_EUL_phi}
\end{equation}
Next, from \eqref{eq:OVERVIEW_EUL_R} one has 
\begin{align}
	\sin\left(\theta\right) & =-R_{\left(3,1\right)}\\
	{\rm cos}\left(\theta\right) & =\sqrt{R_{\left(3,2\right)}^{2}+R_{\left(3,3\right)}^{2}}
\end{align}

\noindent Thus, the pitch angle is equivalent to 
\begin{equation}
	\theta=\arctan\left(\frac{-R_{\left(3,1\right)}}{\sqrt{R_{\left(3,2\right)}^{2}+R_{\left(3,3\right)}^{2}}}\right)\label{eq:OVERVIEW_EUL_theta}
\end{equation}
Finally, from \eqref{eq:OVERVIEW_EUL_R} one has 
\begin{equation}
	\frac{R_{\left(2,1\right)}}{R_{\left(1,1\right)}}=\frac{\cos\left(\theta\right)\sin(\psi)}{\cos\left(\theta\right)\cos(\psi)}\label{eq:-16}
\end{equation}
Thereby, the yaw angle could be found to be 
\begin{equation}
	\psi=\arctan\left(\frac{R_{\left(2,1\right)}}{R_{\left(1,1\right)}}\right)\label{eq:OVERVIEW_EUL_psi}
\end{equation}
Summarizing the above-mentioned results, Euler angles can be obtained
from a given rotational matrix $R$ (\textbf{\textcolor{blue}{the
		mapping from $\mathbb{SO}\left(3\right)$ to $\xi$}}) in the following
manner: 
\begin{equation}
	\left[\begin{array}{c}
		\phi\\
		\theta\\
		\psi
	\end{array}\right]=\left[\begin{array}{c}
		\arctan\left(\frac{R_{\left(3,2\right)}}{R_{\left(3,3\right)}}\right)\\
		\arctan\left(\frac{-R_{\left(3,1\right)}}{\sqrt{R_{\left(3,2\right)}^{2}+R_{\left(3,3\right)}^{2}}}\right)\\
		\arctan\left(\frac{R_{\left(2,1\right)}}{R_{\left(1,1\right)}}\right)
	\end{array}\right]\label{eq:SO3_EUL}
\end{equation}

\subsection{Angular velocities transformation}

\noindent Let us define the Euler angle vector by $\xi=\left[\phi,\theta,\psi\right]^{\top}$
and the body-fixed angular velocity vector by $\Omega=\left[\Omega_{x},\Omega_{y},\Omega_{z}\right]^{\top}$.
The Euler rate vector $\dot{\xi}=\left[\dot{\phi},\dot{\theta},\dot{\psi}\right]^{\top}$
is related to the body-fixed angular velocity vector ($\Omega$) through
a transformation matrix $\mathcal{J}$ such that: 
\begin{equation}
	\dot{\xi}=\mathcal{J}\Omega\label{eq:OVERVIEW_EUL_dot}
\end{equation}
The body-fixed angular velocity vector ($\Omega$) is related to the
Euler rate vector ($\dot{\xi}$) by 
\begin{align*}
	\Omega & =\mathcal{J}^{-1}\dot{\xi}\\
	& =\left[\begin{array}{c}
		\dot{\phi}\\
		0\\
		0
	\end{array}\right]+R_{x,\phi}\left[\begin{array}{c}
		0\\
		\dot{\theta}\\
		0
	\end{array}\right]+R_{x,\phi}R_{y,\theta}\left[\begin{array}{c}
		0\\
		0\\
		\dot{\psi}
	\end{array}\right]\\
	& =\left[\begin{array}{c}
		\dot{\phi}\\
		0\\
		0
	\end{array}\right]+\left[\begin{array}{ccc}
		1 & 0 & 0\\
		0 & c\phi & s\phi\\
		0 & -s\phi & c\phi
	\end{array}\right]\left[\begin{array}{c}
		0\\
		\dot{\theta}\\
		0
	\end{array}\right]+\left[\begin{array}{ccc}
		1 & 0 & 0\\
		0 & c\phi & s\phi\\
		0 & -s\phi & c\phi
	\end{array}\right]\left[\begin{array}{ccc}
		c\theta & 0 & -s\theta\\
		0 & 1 & 0\\
		s\theta & 0 & c\theta
	\end{array}\right]\left[\begin{array}{c}
		0\\
		0\\
		\dot{\psi}
	\end{array}\right]\\
	& =\left[\begin{array}{c}
		\dot{\phi}\\
		0\\
		0
	\end{array}\right]+\left[\begin{array}{c}
		0\\
		\dot{\theta}c\phi\\
		-\dot{\theta}s\phi
	\end{array}\right]+\left[\begin{array}{c}
		-\dot{\psi}s\theta\\
		\dot{\psi}s\phi c\theta\\
		\dot{\psi}c\phi c\theta
	\end{array}\right]
\end{align*}
which is equivalent to 
\begin{equation}
	\Omega=\mathcal{J}^{-1}\dot{\xi}=\left[\begin{array}{ccc}
		1 & 0 & -s\theta\\
		0 & c\phi & s\phi c\theta\\
		0 & -s\phi & c\phi c\theta
	\end{array}\right]\left[\begin{array}{c}
		\dot{\phi}\\
		\dot{\theta}\\
		\dot{\psi}
	\end{array}\right]\label{eq:OVERVIEW_EUL_Omega}
\end{equation}
From \eqref{eq:OVERVIEW_EUL_Omega}, and for a given initial Euler
angles vector $\xi\left(0\right)$, the Euler rate can be found to
be \cite{spong2008robot} 
\begin{equation}
	\dot{\xi}=\mathcal{J}\Omega=\left[\begin{array}{ccc}
		1 & s\phi t\theta & c\phi t\theta\\
		0 & c\phi & -s\phi\\
		0 & s\phi\sec\theta & c\phi\sec\theta
	\end{array}\right]\left[\begin{array}{c}
		\Omega_{x}\\
		\Omega_{y}\\
		\Omega_{z}
	\end{array}\right]\label{eq:OVERVIEW_EUL_Eul_dot}
\end{equation}
which means that the transformation from the body-fixed angular velocity
vector ($\Omega$) to Euler rate vector ($\dot{\xi}$) is given by
\begin{equation}
	\mathcal{J}=\left[\begin{array}{ccc}
		1 & s\phi t\theta & c\phi t\theta\\
		0 & c\phi & -s\phi\\
		0 & s\phi\sec\theta & c\phi\sec\theta
	\end{array}\right]\label{eq:OVERVIEW_EUL_J}
\end{equation}

\subsection{Detailed derivation of Euler rate}

\noindent Let us detail the transformation matrix $\mathcal{J}$.
Recall the attitude dynamics in \eqref{eq:OVERVIEW_SO3_dyn} 
\begin{align*}
	\dot{R} & =R\left[\Omega\right]_{\times}
\end{align*}
Also, recall the transformation matrix in \eqref{eq:OVERVIEW_EUL_R}
\[
R=\left[\begin{array}{ccc}
	r_{11} & r_{12} & r_{13}\\
	r_{21} & r_{22} & r_{23}\\
	r_{31} & r_{32} & r_{33}
\end{array}\right]=\left[\begin{array}{ccc}
	c\theta c\psi & -c\phi s\psi+s\phi s\theta c\psi & s\phi s\psi+c\phi s\theta c\psi\\
	c\theta s\psi & c\phi c\psi+s\phi s\theta s\psi & -s\phi c\psi+c\phi s\theta s\psi\\
	-s\theta & s\phi c\theta & c\phi c\theta
\end{array}\right]
\]
The derivative of the above equations is as follows: 
\begin{equation}
	\dot{R}=\left[\begin{array}{ccc}
		\dot{r}_{11} & \dot{r}_{12} & \dot{r}_{13}\\
		\dot{r}_{21} & \dot{r}_{22} & \dot{r}_{23}\\
		\dot{r}_{31} & \dot{r}_{32} & \dot{r}_{33}
	\end{array}\right]\label{eq:OVERVIEW_EUL_Rdot}
\end{equation}
such that 
\begin{align*}
	\dot{r}_{11} & =-\dot{\theta}s\theta c\psi-\dot{\psi}c\theta s\psi\\
	\dot{r}_{12} & =\dot{\phi}\left(s\phi s\psi+c\phi s\theta c\psi\right)+\dot{\theta}s\phi c\theta c\psi-\dot{\psi}\left(c\phi c\psi+s\phi s\theta s\psi\right)\\
	\dot{r}_{13} & =\dot{\phi}\left(c\phi s\psi-s\phi s\theta c\psi\right)+\dot{\theta}c\phi c\theta c\psi+\dot{\psi}\left(s\phi c\psi-c\phi s\theta s\psi\right)\\
	\dot{r}_{21} & =-\dot{\theta}s\theta s\psi+\dot{\psi}c\theta c\psi\\
	\dot{r}_{22} & =\dot{\phi}\left(c\phi s\theta s\psi-s\phi c\psi\right)+\dot{\theta}s\phi c\theta s\psi+\dot{\psi}\left(s\phi s\theta c\psi-c\phi s\psi\right)\\
	\dot{r}_{23} & =-\dot{\phi}\left(c\phi c\psi+s\phi s\theta s\psi\right)+\dot{\theta}c\phi c\theta s\psi+\dot{\psi}\left(s\phi s\psi+c\phi s\theta c\psi\right)\\
	\dot{r}_{31} & =-\dot{\theta}c\theta\\
	\dot{r}_{32} & =\dot{\phi}c\phi c\theta-\dot{\theta}s\phi s\theta\\
	\dot{r}_{33} & =-\dot{\phi}s\phi c\theta-\dot{\theta}c\phi s\theta
\end{align*}
From \eqref{eq:OVERVIEW_EUL_R} and \eqref{eq:OVERVIEW_EUL_Rdot}
one can find 
\begin{align}
	R^{\top}\dot{R} & =\left[\begin{array}{ccc}
		0 & -\dot{\psi}c\phi c\theta+\dot{\theta}s\phi & \dot{\theta}c\phi+\dot{\psi}s\phi c\theta\\
		\dot{\psi}c\phi c\theta-\dot{\theta}s\phi & 0 & -\dot{\phi}+\dot{\psi}s\theta\\
		-\dot{\theta}c\phi-\dot{\psi}s\phi c\theta & \dot{\phi}-\dot{\psi}s\theta & 0
	\end{array}\right]\label{eq:OVERVIEW_EUL_RRdot}
\end{align}
From \eqref{eq:OVERVIEW_SO3_dyn} and \eqref{eq:OVERVIEW_EUL_RRdot}
one has 
\begin{align*}
	R^{\top}\dot{R} & =\left[\Omega\right]_{\times}\\
	\left[\begin{array}{ccc}
		0 & -\dot{\psi}c\phi c\theta+\dot{\theta}s\phi & \dot{\theta}c\phi+\dot{\psi}s\phi c\theta\\
		\dot{\psi}c\phi c\theta-\dot{\theta}s\phi & 0 & -\dot{\phi}+\dot{\psi}s\theta\\
		-\dot{\theta}c\phi-\dot{\psi}s\phi c\theta & \dot{\phi}-\dot{\psi}s\theta & 0
	\end{array}\right] & =\left[\begin{array}{ccc}
		0 & -\Omega_{z} & \Omega_{y}\\
		\Omega_{z} & 0 & -\Omega_{x}\\
		-\Omega_{y} & \Omega_{x} & 0
	\end{array}\right]
\end{align*}
Thus, the body-fixed angular velocity vector can be expressed as 
\begin{align}
	\mathbf{vex}\left(\left[\Omega\right]_{\times}\right) & =\mathbf{vex}\left(R^{\top}\dot{R}\right)\nonumber \\
	\left[\begin{array}{c}
		\Omega_{x}\\
		\Omega_{y}\\
		\Omega_{z}
	\end{array}\right] & =\left[\begin{array}{c}
		\dot{\phi}-\dot{\psi}s\theta\\
		\dot{\theta}c\phi+\dot{\psi}s\phi c\theta\\
		\dot{\psi}c\phi c\theta-\dot{\theta}s\phi
	\end{array}\right]\nonumber \\
	\Omega & =\mathcal{J}^{-1}\dot{\xi}=\left[\begin{array}{ccc}
		1 & 0 & -s\theta\\
		0 & c\phi & s\phi c\theta\\
		0 & -s\phi & c\phi c\theta
	\end{array}\right]\left[\begin{array}{c}
		\dot{\phi}\\
		\dot{\theta}\\
		\dot{\psi}
	\end{array}\right]\label{eq:OVERVIEW_EUL_Omega2}
\end{align}
Therefore, the Euler rate becomes 
\begin{equation}
	\dot{\xi}=\mathcal{J}\Omega=\left[\begin{array}{ccc}
		1 & s\phi t\theta & c\phi t\theta\\
		0 & c\phi & -s\phi\\
		0 & s\phi\sec\theta & c\phi\sec\theta
	\end{array}\right]\left[\begin{array}{c}
		\Omega_{x}\\
		\Omega_{y}\\
		\Omega_{z}
	\end{array}\right]\label{eq:OVERVIEW_EUL_Eul_dot-1}
\end{equation}
and the transformation matrix is equivalent to 
\begin{equation}
	\mathcal{J}=\left[\begin{array}{ccc}
		1 & s\phi t\theta & c\phi t\theta\\
		0 & c\phi & -s\phi\\
		0 & s\phi\sec\theta & c\phi\sec\theta
	\end{array}\right]\label{eq:OVERVIEW_EUL_J-1}
\end{equation}

\subsection{Euler parameterization problem}

Euler angles are kinematically singular. In other words, the transformation
of the Euler angles rates ($\dot{\xi}$) to the angular velocity ($\Omega$)
is locally defined, while there is no global definition for the transformation
matrix ($\mathcal{J}$). At certain configurations, the mapping from
$\mathbb{SO}\left(3\right)$ to Euler angles ($\xi$) could result
in infinite solutions. For example, consider the following orthogonal
matrix 
\begin{equation}
	R=\left[\begin{array}{ccc}
		0 & 0 & 1\\
		0 & 1 & 0\\
		-1 & 0 & 0
	\end{array}\right]\in\mathbb{SO}\left(3\right)\left\{ \begin{array}{c}
		{\rm det}\left(R\right)=+1\\
		R^{\top}R=RR^{\top}=\mathbf{I}_{3}
	\end{array}\right.\label{eq:OVERVIEW_EUL_5}
\end{equation}
It is impossible to obtain the true Euler angles of the orthogonal
matrix in \eqref{eq:OVERVIEW_EUL_5}, and there exist infinite solutions,
since

\[
R_{\left(\phi,\theta,\psi\right)}=R_{\left(90,90,90\right)}=R_{\left(0,90,0\right)}=R_{\left(45,90,45\right)}=R_{\left(\star,90,\star\right)}=\cdots=\left[\begin{array}{ccc}
	0 & 0 & 1\\
	0 & 1 & 0\\
	-1 & 0 & 0
\end{array}\right]
\]
where $\star$ denotes any angle. Hence, for $R=\left[\begin{array}{ccc}
	0 & 0 & 1\\
	0 & 1 & 0\\
	-1 & 0 & 0
\end{array}\right]$, there exists infinite number of solutions for the Euler angles,
and as a result the true Euler angles cannot be obtained. Thus, despite
Euler angles being capable of describing every attitude, the representation
produced by the set of Euler angles is not unique. Consequently, these
representations are limited to local attitude maneuvers, due to the
fact that continuous control maneuvers are not guaranteed for certain
configurations. In addition, there are angular velocities ($\Omega$)
that cannot be obtained by the means of the time derivatives of Euler
angles ($\dot{\xi}$).

\noindent %
\noindent\makebox[1\linewidth]{%
	\rule{0.6\textwidth}{1.4pt}%
}

\section{Angle-axis Parameterization\label{sec:OVERVIEW_Ang_Axis}}

\subsection{Mapping: From/To unit vector and angle of rotation to/from other
	representations}

\begin{figure}[h]
	\centering{}\includegraphics[scale=0.45]{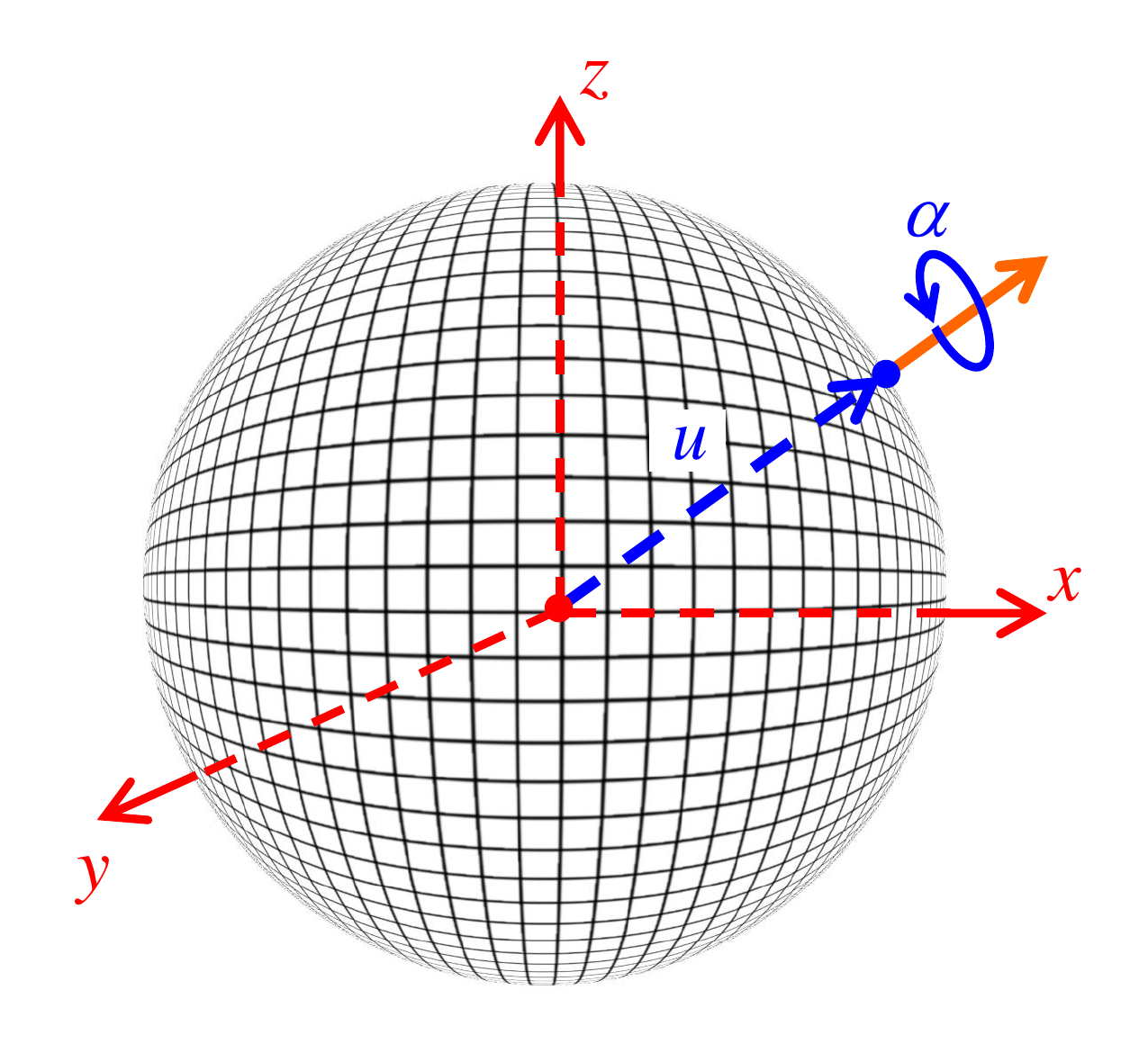}\caption{The orientation of a rigid-body in body-frame relative to inertial-frame. }
	\label{fig:OVERVIEW_SO3_3} 
\end{figure}

The relative orientation of any two frames can always be expressed
in terms of a single rotation about a given normalized vector with
a given rotation angle. Let $u=\left[u_{1},u_{2},u_{3}\right]^{\top}\in\mathbb{S}^{2}$
denote a unit vector such that $\left\Vert u\right\Vert =\sqrt{u_{1}^{2}+u_{2}^{2}+u_{3}^{2}}=1$,
and $\alpha\in\mathbb{R}$ denote the angle of rotation about $u$.
Then for a given rotation angle ($\alpha$) and unit axis ($u$),
the corresponding rotational matrix $\mathcal{R}_{\alpha}:\mathbb{R}\times\mathbb{S}^{2}\rightarrow\mathbb{SO}\left(3\right)$
(\textbf{\textcolor{blue}{the mapping from $\left(\alpha,u\right)$
		to $\mathbb{SO}\left(3\right)$}}) can be expressed by the following
formula: 
\begin{align}
	\mathcal{R}_{\alpha}\left(\alpha,u\right) & =\exp\left(\left[\alpha u\right]_{\times}\right)\in\mathbb{SO}\left(3\right)\label{eq:OVERVIEW_att_ang3-1}
\end{align}
Figure \ref{fig:OVERVIEW_SO3_3} visualizes the angle-axis parameterization
such that any point on a unit sphere can be described by a rotation
angle ($\alpha$) and the unit axis ($u$).

This isomorphism lies at the core of the Euler's theorem, which states
that the axis of rotation $u\in\mathbb{S}^{2}$ with rotational angle
$\alpha\in\mathbb{R}$ yields the rotational matrix $\mathcal{R}_{\alpha}\left(\alpha,u\right)=\exp\left(\left[\alpha u\right]_{\times}\right)\in\mathbb{SO}\left(3\right)$.
The isomorphism in \eqref{eq:OVERVIEW_att_ang3-1} can be expressed
with respect to Rodriguez formula by \cite{murray1994mathematicall}
\begin{align}
	\mathcal{R}_{\alpha}\left(\alpha,u\right) & =\mathbf{I}_{3}+\sin\left(\alpha\right)\left[u\right]_{\times}+\left(1-\cos\left(\alpha\right)\right)\left[u\right]_{\times}^{2}\in\mathbb{SO}\left(3\right)\label{eq:OVERVIEW_att_ang3-1-1}
\end{align}
Hence{\small{}{}{}{} 
	\begin{align}
		\mathcal{R}_{\alpha}\left(\alpha,u\right)= & \mathbf{I}_{3}+\sin\left(\alpha\right)\left[\begin{array}{ccc}
			0 & -u_{3} & u_{2}\\
			u_{3} & 0 & -u_{1}\\
			-u_{2} & u_{1} & 0
		\end{array}\right]+\left(1-\cos\left(\alpha\right)\right)\left[\begin{array}{ccc}
			0 & -u_{3} & u_{2}\\
			u_{3} & 0 & -u_{1}\\
			-u_{2} & u_{1} & 0
		\end{array}\right]\left[\begin{array}{ccc}
			0 & -u_{3} & u_{2}\\
			u_{3} & 0 & -u_{1}\\
			-u_{2} & u_{1} & 0
		\end{array}\right]\nonumber \\
		= & \left[\begin{array}{ccc}
			1 & 0 & 0\\
			0 & 1 & 0\\
			0 & 0 & 1
		\end{array}\right]+\sin\left(\alpha\right)\left[\begin{array}{ccc}
			0 & -u_{3} & u_{2}\\
			u_{3} & 0 & -u_{1}\\
			-u_{2} & u_{1} & 0
		\end{array}\right]\nonumber \\
		& +\left(1-\cos\left(\alpha\right)\right)\left[\begin{array}{ccc}
			-\left(u_{2}^{2}+u_{3}^{2}\right) & u_{1}u_{2} & u_{1}u_{3}\\
			u_{1}u_{2} & -\left(u_{1}^{2}+u_{3}^{2}\right) & u_{2}u_{3}\\
			u_{1}u_{3} & u_{2}u_{3} & -\left(u_{1}^{2}+u_{2}^{2}\right)
		\end{array}\right]\label{eq:OVERVIEW_att_ang1}
	\end{align}
}Consider the following properties 
\begin{align*}
	\sin^{2}\left(\alpha\right) & =\frac{1}{2}-\frac{1}{2}\cos(2\alpha)\\
	2\sin^{2}(\alpha/2) & =1-\cos\left(\alpha\right)
\end{align*}
from \eqref{eq:OVERVIEW_att_ang1} one has{\small{}{}{}{} 
	\begin{align*}
		\mathcal{R}_{\alpha}\left(\alpha,u\right)= & \left[\begin{array}{ccc}
			1 & 0 & 0\\
			0 & 1 & 0\\
			0 & 0 & 1
		\end{array}\right]+\sin\left(\alpha\right)\left[\begin{array}{ccc}
			0 & -u_{3} & u_{2}\\
			u_{3} & 0 & -u_{1}\\
			-u_{2} & u_{1} & 0
		\end{array}\right]\\
		& +2\sin^{2}\left(\alpha/2\right)\left[\begin{array}{ccc}
			-\left(u_{2}^{2}+u_{3}^{2}\right) & u_{1}u_{2} & u_{1}u_{3}\\
			u_{1}u_{2} & -\left(u_{1}^{2}+u_{3}^{2}\right) & u_{2}u_{3}\\
			u_{1}u_{3} & u_{2}u_{3} & -\left(u_{1}^{2}+u_{2}^{2}\right)
		\end{array}\right]
	\end{align*}
}which is equivalent to {\small{}{}{}{} 
	\begin{align}
		& \mathcal{R}_{\alpha}\left(\alpha,u\right)=\nonumber \\
		& \hspace{1em}\left[\begin{array}{ccc}
			1-2\left(u_{2}^{2}+u_{3}^{2}\right)\sin^{2}\left(\alpha/2\right) & -\sin\left(\alpha\right)u_{3}+2u_{1}u_{2}\sin^{2}\left(\alpha/2\right) & \sin\left(\alpha\right)u_{2}+2u_{1}u_{3}\sin^{2}\left(\alpha/2\right)\\
			\sin\left(\alpha\right)u_{3}+2u_{1}u_{2}\sin^{2}\left(\alpha/2\right) & 1-2\left(u_{1}^{2}+u_{3}^{2}\right)\sin^{2}\left(\alpha/2\right) & -\sin\left(\alpha\right)u_{1}+2u_{2}u_{3}\sin^{2}\left(\alpha/2\right)\\
			-\sin\left(\alpha\right)u_{2}+2u_{1}u_{3}\sin^{2}\left(\alpha/2\right) & \sin\left(\alpha\right)u_{1}+2u_{2}u_{3}\sin^{2}\left(\alpha/2\right) & 1-2\left(u_{1}^{2}+u_{2}^{2}\right)\sin^{2}\left(\alpha/2\right)
		\end{array}\right]\label{eq:OVERVIEW_att_ang2}
	\end{align}
}The angle of rotation $\alpha$ and the unit axis $u$ can be easily
obtained knowing the orientation matrix $\mathcal{R}_{\alpha}\left(\alpha,u\right)\in\mathbb{SO}\left(3\right)$
in \eqref{eq:OVERVIEW_att_ang2}. From \eqref{eq:OVERVIEW_att_ang1},
the angle $\alpha$ (\textbf{\textcolor{blue}{the mapping from $\mathbb{SO}\left(3\right)$
		to $\left(\alpha,u\right)$}}) is given by 
\begin{align}
	\alpha & ={\rm cos}^{-1}\left(\frac{1}{2}\left({\rm Tr}\left\{ R\right\} -1\right)\right)\label{eq:OVERVIEW_att_ang_alpha}
\end{align}
Knowing $\alpha$, the unit vector $u$ can be obtained by 
\begin{align}
	u & =\frac{1}{{\rm sin}\left(\alpha\right)}\mathbf{vex}\left(\frac{1}{2}\left(R-R^{\top}\right)\right)\nonumber \\
	& =\frac{1}{{\rm sin}\left(\alpha\right)}\mathbf{vex}\left(\boldsymbol{\mathcal{P}}_{a}\left(R\right)\right)\label{eq:OVERVIEW_att_ang_u}
\end{align}
From \eqref{eq:OVERVIEW_att_ang1}, the normalized Euclidean distance
can be defined in terms of angle-axis components by 
\begin{align}
	\left\Vert R\right\Vert _{I} & =\frac{1}{2}\left(1-\cos\left(\alpha\right)\right)\nonumber \\
	& =\sin^{2}\left(\frac{\alpha}{2}\right)\label{eq:OVERVIEW_att_ang_RI}
\end{align}
This shows \eqref{eq:OVERVIEW_SO3_lem6_2} in Lemma \ref{Lem:OVERVIEW_6}.
Hence, the following holds

\noindent 
\begin{align}
	{\rm cos}^{2}\left(\frac{\alpha}{2}\right) & =1-\sin^{2}\left(\frac{\alpha}{2}\right)\nonumber \\
	& =1-\left\Vert R\right\Vert _{I}\label{eq:OVERVIEW_att_ang3}
\end{align}
From \eqref{eq:OVERVIEW_att_ang3} the anti-symmetric operator of
the rotational matrix can be expressed with regards to angle-axis
components by

\begin{align*}
	\boldsymbol{\mathcal{P}}_{a}\left(R\right) & =\sin\left(\alpha\right)\left[u\right]_{\times}\\
	& =2{\rm cos}\left(\frac{\alpha}{2}\right)\sin\left(\frac{\alpha}{2}\right)\left[u\right]_{\times}
\end{align*}
As such, the vex operator can be found to be 
\begin{equation}
	\mathbf{vex}\left(\boldsymbol{\mathcal{P}}_{a}\left(R\right)\right)=2{\rm cos}\left(\frac{\alpha}{2}\right)\sin\left(\frac{\alpha}{2}\right)u\label{eq:OVERVIEW_att_ang_VEX}
\end{equation}
This proves \eqref{eq:OVERVIEW_SO3_lem6_1} in Lemma \ref{Lem:OVERVIEW_6}.
The norm square of the result in \eqref{eq:OVERVIEW_att_ang_VEX}
is 
\begin{align}
	\left\Vert \mathbf{vex}\left(\boldsymbol{\mathcal{P}}_{a}\left(R\right)\right)\right\Vert ^{2} & =\mathbf{vex}\left(\boldsymbol{\mathcal{P}}_{a}\left(R\right)\right)^{\top}\mathbf{vex}\left(\boldsymbol{\mathcal{P}}_{a}\left(R\right)\right)\nonumber \\
	& =u^{\top}u\sin^{2}\left(\alpha\right),\hspace{1em}\left\Vert u\right\Vert ^{2}=1\nonumber \\
	& =4{\rm cos}^{2}\left(\frac{\alpha}{2}\right)\sin^{2}\left(\frac{\alpha}{2}\right)\label{eq:OVERVIEW_att_ang_Vex1}
\end{align}
Considering the results of \eqref{eq:OVERVIEW_att_ang_RI} and \eqref{eq:OVERVIEW_att_ang3},
one can find 
\begin{align}
	\left\Vert \mathbf{vex}\left(\boldsymbol{\mathcal{P}}_{a}\left(R\right)\right)\right\Vert ^{2} & =4{\rm cos}^{2}\left(\frac{\alpha}{2}\right)\sin^{2}\left(\frac{\alpha}{2}\right)\nonumber \\
	& =4\left(1-\sin^{2}\left(\frac{\alpha}{2}\right)\right)\sin^{2}\left(\frac{\alpha}{2}\right)\nonumber \\
	& =4\left(1-\left\Vert R\right\Vert _{I}\right)\left\Vert R\right\Vert _{I}\label{eq:OVERVIEW_att_ang_Vex2}
\end{align}
This confirms \eqref{eq:OVERVIEW_SO3_lem5_1} in Lemma \ref{Lem:OVERVIEW_5}.

\subsection{Problems of angle-axis parameterization: }

The problem of expressing the attitude using angle-axis components
is that some rotations are not defined in terms of angle-axis parameterization.
For instance, consider the following results 
\begin{align*}
	\alpha & ={\rm cos}^{-1}\left(\frac{1}{2}\left({\rm Tr}\left\{ R\right\} -1\right)\right)\\
	u & =\frac{1}{2{\rm sin}\left(\alpha\right)}\left[\begin{array}{c}
		R_{32}-R_{23}\\
		R_{13}-R_{31}\\
		R_{21}-R_{12}
	\end{array}\right]
\end{align*}
Thus, the following four orientations cannot use angle-axis components
in their definition 
\begin{align}
	\text{at }R & =\left[\begin{array}{ccc}
		1 & 0 & 0\\
		0 & 1 & 0\\
		0 & 0 & 1
	\end{array}\right]\in\mathbb{SO}\left(3\right)\Rightarrow\alpha=0\in\mathbb{R},\hspace{1em}u=\left[\infty,\infty,\infty\right]^{\top}\notin\mathbb{S}^{2}\label{eq:-18}\\
	\text{at }R & =\left[\begin{array}{ccc}
		1 & 0 & 0\\
		0 & -1 & 0\\
		0 & 0 & -1
	\end{array}\right]\in\mathbb{SO}\left(3\right)\Rightarrow\alpha=\pi\in\mathbb{R},\hspace{1em}u=\left[\infty,\infty,\infty\right]^{\top}\notin\mathbb{S}^{2}\label{eq:-19}\\
	\text{at }R & =\left[\begin{array}{ccc}
		-1 & 0 & 0\\
		0 & 1 & 0\\
		0 & 0 & -1
	\end{array}\right]\in\mathbb{SO}\left(3\right)\Rightarrow\alpha=\pi\in\mathbb{R},\hspace{1em}u=\left[\infty,\infty,\infty\right]^{\top}\notin\mathbb{S}^{2}\label{eq:-20}\\
	\text{at }R & =\left[\begin{array}{ccc}
		-1 & 0 & 0\\
		0 & -1 & 0\\
		0 & 0 & 1
	\end{array}\right]\in\mathbb{SO}\left(3\right)\Rightarrow\alpha=\pi\in\mathbb{R},\hspace{1em}u=\left[\infty,\infty,\infty\right]^{\top}\notin\mathbb{S}^{2}\label{eq:-21}
\end{align}
Therefore, $\alpha$ should not be a multiple of $\pi$, while the
representation of angle-axis parameterization is valid for all $\alpha\neq k\pi\forall k=0,1,2,3,\ldots$.

\noindent %
\noindent\makebox[1\linewidth]{%
	\rule{0.6\textwidth}{1.4pt}%
}

\section{Rodriguez Vector Parameterization\label{sec:OVERVIEW_ROD}}

\subsection{Mapping: From/To Rodriguez vector to/from other representations}

The Rodriguez vector $\rho=\left[\rho_{1},\rho_{2},\rho_{3}\right]^{\top}\in\mathbb{R}^{3}$,
proposed by Olinde Rodriguez \cite{rodrigues1840lois}, can be employed
for attitude representation. It should be noted, that naming conventions
vary and in some articles Rodriguez vector is termed Gibbs vector
\cite{gibbs1961scientific,wilson1901vector}. The orthogonal matrix
$\mathcal{R}_{\rho}$ is obtained through mapping the vector on $\mathbb{R}^{3}$
to $\mathbb{SO}\left(3\right)$ such that $\mathcal{R}_{\rho}:\mathbb{R}^{3}\rightarrow\mathbb{SO}\left(3\right)$
\cite{tsiotras1997higher}. The Cayley transform is given by \cite{tsiotras1997higher}
\begin{equation}
	\mathcal{R}_{\rho}=\left(\mathbf{I}_{3}+\left[\rho\right]_{\times}\right)\left(\mathbf{I}_{3}-\left[\rho\right]_{\times}\right)^{-1},\hspace{1em}\mathcal{R}_{\rho}\in\mathbb{SO}\left(3\right)\label{eq:OVERVIEW_ROD_SO3_R}
\end{equation}
It is important to present the following properties \cite{cayley1845xiii}: 
\begin{enumerate}
	\item The matrix multiplication in \eqref{eq:OVERVIEW_ROD_SO3_R} is commutative,
	so $\mathcal{R}_{\rho}$ can be alternatively defined as 
	\begin{align*}
		\mathcal{R}_{\rho} & =\left(\mathbf{I}_{3}+\left[\rho\right]_{\times}\right)\left(\mathbf{I}_{3}-\left[\rho\right]_{\times}\right)^{-1}\in\mathbb{SO}\left(3\right)\\
		& =\left(\mathbf{I}_{3}-\left[\rho\right]_{\times}\right)^{-1}\left(\mathbf{I}_{3}+\left[\rho\right]_{\times}\right)\in\mathbb{SO}\left(3\right)
	\end{align*}
	and 
	\begin{align*}
		\mathcal{R}_{\rho}^{\top} & =\left(\mathbf{I}_{3}-\left[\rho\right]_{\times}\right)\left(\mathbf{I}_{3}+\left[\rho\right]_{\times}\right)^{-1}\in\mathbb{SO}\left(3\right)\\
		& =\left(\mathbf{I}_{3}+\left[\rho\right]_{\times}\right)^{-1}\left(\mathbf{I}_{3}-\left[\rho\right]_{\times}\right)\in\mathbb{SO}\left(3\right)
	\end{align*}
	\item The determinant of $\mathcal{R}_{\rho}$ must always be $+1$. 
	\item For $\det\left(\mathcal{R}_{\rho}\right)=-1$, Rodriguez vector cannot
	be defined. 
	\item $\mathcal{R}_{\rho}^{\top}\mathcal{R}_{\rho}=\mathcal{R}_{\rho}\mathcal{R}_{\rho}^{\top}=\mathcal{R}_{\rho}^{-1}\mathcal{R}_{\rho}=\mathcal{R}_{\rho}\mathcal{R}_{\rho}^{-1}=\mathbf{I}_{3}$. 
\end{enumerate}
The related map from Rodriguez parameters vector form to $\mathbb{SO}\left(3\right)$
($\mathcal{R}_{\rho}:\mathbb{R}^{3}\rightarrow\mathbb{SO}\left(3\right)$)
is \cite{wilson1901vector} 
\begin{align}
	\mathcal{R}_{\rho}= & \frac{1}{1+\left\Vert \rho\right\Vert ^{2}}\left(\left(1-\left\Vert \rho\right\Vert ^{2}\right)\mathbf{I}_{3}+2\rho\rho^{\top}+2\left[\rho\right]_{\times}\right)\label{eq:OVERVIEW_ROD_SO3}
\end{align}
To proof the result in \eqref{eq:OVERVIEW_ROD_SO3}, we solve for
\eqref{eq:OVERVIEW_ROD_SO3_R}. The following derivation shows mapping
$\mathcal{R}_{\rho}:\mathbb{R}^{3}\rightarrow\mathbb{SO}\left(3\right)$

\noindent 
\begin{align*}
	\left(\mathbf{I}_{3}-\left[\rho\right]_{\times}\right)^{-1} & =\left[\begin{array}{ccc}
		1 & \rho_{3} & -\rho_{2}\\
		-\rho_{3} & 1 & \rho_{1}\\
		\rho_{2} & -\rho_{1} & 1
	\end{array}\right]^{-1}\\
	& =\frac{1}{1+\rho_{1}^{2}+\rho_{2}^{2}+\rho_{3}^{2}}\left[\begin{array}{ccc}
		1+\rho_{1}^{2} & \rho_{1}\rho_{2}-\rho_{3} & \rho_{1}\rho_{3}+\rho_{2}\\
		\rho_{1}\rho_{2}+\rho_{3} & 1+\rho_{2}^{2} & \rho_{2}\rho_{3}-\rho_{1}\\
		\rho_{1}\rho_{3}-\rho_{2} & \rho_{2}\rho_{3}+\rho_{1} & 1+\rho_{3}^{2}
	\end{array}\right]\\
	& =\frac{1}{1+\left\Vert \rho\right\Vert ^{2}}\left(\mathbf{I}_{3}+\rho\rho^{\top}+\left[\rho\right]_{\times}\right)
\end{align*}
Considering the identities in \eqref{eq:OVERVIEW_PER_SO3_5} and \eqref{eq:OVERVIEW_PER_SO3_VEX_7},
we have $\left[\rho\right]_{\times}\rho=0$ and $\left[\rho\right]_{\times}^{2}=\rho\rho^{\top}-\left\Vert \rho\right\Vert ^{2}\mathbf{I}_{3}$
such that 
\begin{align*}
	\mathcal{R}_{\rho}= & \left(\mathbf{I}_{3}+\left[\rho\right]_{\times}\right)\left(\mathbf{I}_{3}-\left[\rho\right]_{\times}\right)^{-1}\\
	= & \frac{1}{1+\left\Vert \rho\right\Vert ^{2}}\left(\mathbf{I}_{3}+\left[\rho\right]_{\times}\right)\left(\mathbf{I}_{3}+\rho\rho^{\top}+\left[\rho\right]_{\times}\right)\\
	= & \frac{1}{1+\left\Vert \rho\right\Vert ^{2}}\left(\mathbf{I}_{3}+2\left[\rho\right]_{\times}+\left[\rho\right]_{\times}\rho\rho^{\top}+\rho\rho^{\top}+\left[\rho\right]_{\times}^{2}\right)\\
	= & \frac{1}{1+\left\Vert \rho\right\Vert ^{2}}\left(\left(1-\left\Vert \rho\right\Vert ^{2}\right)\mathbf{I}_{3}+2\left[\rho\right]_{\times}+2\rho\rho^{\top}\right)
\end{align*}
The end result is $\mathcal{R}_{\rho}:\mathbb{R}^{3}\rightarrow\mathbb{SO}\left(3\right)$
(\textbf{\textcolor{blue}{the mapping from $\rho$ to $\mathbb{SO}\left(3\right)$}})
\textcolor{red}{{} } 
\begin{align}
	\mathcal{R}_{\rho}= & \frac{1}{1+\left\Vert \rho\right\Vert ^{2}}\left(\left(1-\left\Vert \rho\right\Vert ^{2}\right)\mathbf{I}_{3}+2\left[\rho\right]_{\times}+2\rho\rho^{\top}\right)\label{eq:OVERVIEW_ROD_SO3_1}
\end{align}
which can be written in detailed form as follows 
\begin{equation}
	\mathcal{R}_{\rho}=\frac{1}{1+\left\Vert \rho\right\Vert ^{2}}\left[\begin{array}{ccc}
		1+\rho_{1}^{2}-\rho_{2}^{2}-\rho_{3}^{2} & 2\left(\rho_{1}\rho_{2}-\rho_{3}\right) & 2\left(\rho_{1}\rho_{3}+\rho_{2}\right)\\
		2\left(\rho_{1}\rho_{2}+\rho_{3}\right) & 1+\rho_{2}^{2}-\rho_{3}^{2}-\rho_{1}^{2} & 2\left(\rho_{2}\rho_{3}-\rho_{1}\right)\\
		2\left(\rho_{1}\rho_{3}-\rho_{2}\right) & 2\left(\rho_{2}\rho_{3}+\rho_{1}\right) & 1+\rho_{3}^{2}-\rho_{1}^{2}-\rho_{2}^{2}
	\end{array}\right]\label{eq:-29}
\end{equation}

\noindent Now, let us present the proof of orthogonality

\begin{align}
	\mathcal{R}_{\rho}^{\top}\mathcal{R}_{\rho} & =\left(\left(\mathbf{I}_{3}-\left[\rho\right]_{\times}\right)^{-1}\left(\mathbf{I}_{3}+\left[\rho\right]_{\times}\right)\right)^{\top}\left(\left(\mathbf{I}_{3}+\left[\rho\right]_{\times}\right)\left(\mathbf{I}_{3}-\left[\rho\right]_{\times}\right)^{-1}\right)\nonumber \\
	& =\left(\mathbf{I}_{3}-\left[\rho\right]_{\times}\right)\left(\mathbf{I}_{3}+\left[\rho\right]_{\times}\right)^{-1}\left(\mathbf{I}_{3}+\left[\rho\right]_{\times}\right)\left(\mathbf{I}_{3}-\left[\rho\right]_{\times}\right)^{-1}\nonumber \\
	& =\left(\mathbf{I}_{3}-\left[\rho\right]_{\times}\right)\left(\mathbf{I}_{3}-\left[\rho\right]_{\times}\right)^{-1}\nonumber \\
	& =\mathbf{I}_{3}\label{eq:-9}
\end{align}
Hence{} 
\begin{equation}
	\mathcal{R}_{\rho}^{\top}\mathcal{R}_{\rho}=\mathbf{I}_{3}\label{eq:OVERVIEW_ROD_ORTH}
\end{equation}
Let us prove the inverse mapping: $\rho:\mathbb{SO}\left(3\right)\rightarrow\mathbb{R}^{3}$

\noindent 
\begin{align*}
	\mathcal{R}_{\rho} & =\left(\mathbf{I}_{3}+\left[\rho\right]_{\times}\right)\left(\mathbf{I}_{3}-\left[\rho\right]_{\times}\right)^{-1}\\
	\mathcal{R}_{\rho}\left(\mathbf{I}_{3}-\left[\rho\right]_{\times}\right) & =\mathbf{I}_{3}+\left[\rho\right]_{\times}\\
	\mathcal{R}_{\rho}-\mathbf{I}_{3} & =\mathcal{R}_{\rho}\left[\rho\right]_{\times}+\left[\rho\right]_{\times}\\
	\mathcal{R}_{\rho}-\mathbf{I}_{3} & =\left(\mathcal{R}_{\rho}+\mathbf{I}_{3}\right)\left[\rho\right]_{\times}\\
	\left[\rho\right]_{\times} & =\left(\mathcal{R}_{\rho}+\mathbf{I}_{3}\right)^{-1}\left(\mathcal{R}_{\rho}-\mathbf{I}_{3}\right)
\end{align*}
Hence, one can find $\rho:\mathbb{SO}\left(3\right)\rightarrow\mathbb{R}^{3}$
(\textbf{\textcolor{blue}{the mapping from $\mathbb{SO}\left(3\right)$
		to $\rho$}}) as{} \cite{rodrigues1840lois}{} 
\begin{equation}
	\rho=\mathbf{vex}\left(\left(\mathcal{R}_{\rho}+\mathbf{I}_{3}\right)^{-1}\left(\mathcal{R}_{\rho}-\mathbf{I}_{3}\right)\right)\label{eq:OVERVIEW_ROD}
\end{equation}
which is also equivalent to 
\begin{equation}
	\rho=\frac{1}{1+{\rm Tr}\left\{ R\right\} }\left[\begin{array}{c}
		R_{32}-R_{23}\\
		R_{13}-R_{31}\\
		R_{21}-R_{12}
	\end{array}\right]\label{eq:OVERVIEW_SO3_ROD}
\end{equation}

\noindent From, \eqref{eq:OVERVIEW_EUL_R}, \eqref{eq:SO3_EUL}, and
\eqref{eq:OVERVIEW_ROD_SO3}, the related map from Rodriguez vector
to Euler angles $\xi:\mathbb{R}^{3}\rightarrow\mathbb{R}^{3}$ (\textbf{\textcolor{blue}{the
		mapping from $\rho$ to $\xi$}}) can be defined by

\noindent 
\begin{equation}
	\left[\begin{array}{c}
		\phi\\
		\theta\\
		\psi
	\end{array}\right]=\left[\begin{array}{c}
		\arctan\left(\frac{2\rho_{2}\rho_{3}+2\rho_{1}}{1+\rho_{3}^{2}-\rho_{1}^{2}-\rho_{2}^{2}}\right)\\
		\arctan\left(\frac{2\rho_{2}-2\rho_{1}\rho_{3}}{\sqrt{4\left(\rho_{2}\rho_{3}+\rho_{1}\right)^{2}+\left(1+\rho_{3}^{2}-\rho_{1}^{2}-\rho_{2}^{2}\right)^{2}}}\right)\\
		\arctan\left(\frac{2\rho_{1}\rho_{2}+2\rho_{3}}{1+\rho_{1}^{2}-\rho_{2}^{2}-\rho_{3}^{2}}\right)
	\end{array}\right]\label{eq:-30}
\end{equation}
Let us define the relationship between $\boldsymbol{\mathcal{P}}_{a}\left(R\right)$
and the normalized Euclidean distance $\left\Vert R\right\Vert _{I}$.
With direct substitution of \eqref{eq:OVERVIEW_ROD_SO3} in \eqref{eq:OVERVIEW_SO3_Ecul_Dist}
one obtains\textcolor{red}{{} } 
\begin{align}
	\left\Vert R\right\Vert _{I} & =\frac{1}{4}{\rm Tr}\left\{ \mathbf{I}_{3}-\mathcal{R}_{\rho}\right\} \nonumber \\
	& =\frac{1}{4}{\rm Tr}\left\{ \mathbf{I}_{3}-\frac{1}{1+\left\Vert \rho\right\Vert ^{2}}\left(\left(1-\left\Vert \rho\right\Vert ^{2}\right)\mathbf{I}_{3}+2\left[\rho\right]_{\times}+2\rho\rho^{\top}\right)\right\} \nonumber \\
	& =\frac{1}{4}\frac{1}{1+\left\Vert \rho\right\Vert ^{2}}{\rm Tr}\left\{ 2\left\Vert \rho\right\Vert ^{2}\mathbf{I}_{3}-2\rho\rho^{\top}-2\left[\rho\right]_{\times}\right\} \nonumber \\
	& =\frac{1}{4}\frac{1}{1+\left\Vert \rho\right\Vert ^{2}}\left(6\left\Vert \rho\right\Vert ^{2}-2\left\Vert \rho\right\Vert ^{2}\right)\label{eq:OVERVIEW_ROD_TR211}
\end{align}
Where ${\rm Tr}\left\{ \left[\rho\right]_{\times}\right\} =0$ and
${\rm Tr}\left\{ \rho\rho^{\top}\right\} =\left\Vert \rho\right\Vert ^{2}$.
Thus, the normalized Euclidean distance $\left\Vert R\right\Vert _{I}$
in terms of Rodriguez vector is equivalent to 
\begin{equation}
	\left\Vert R\right\Vert _{I}=\frac{\left\Vert \rho\right\Vert ^{2}}{1+\left\Vert \rho\right\Vert ^{2}}\label{eq:OVERVIEW_ROD_TR2}
\end{equation}
This proves \eqref{eq:OVERVIEW_SO3_lem1_2} in Lemma \ref{Lem:OVERVIEW_1}.
Likewise, the anti-symmetric projection operator of attitude $R$
in \eqref{eq:OVERVIEW_ROD_SO3} can be defined as 
\begin{align*}
	\boldsymbol{\mathcal{P}}_{a}\left(R\right)= & \frac{1}{2}\frac{1}{1+\left\Vert \rho\right\Vert ^{2}}\left(\mathcal{R}_{\rho}-\mathcal{R}_{\rho}\right)\\
	= & \frac{1}{2}\frac{1}{1+\left\Vert \rho\right\Vert ^{2}}\left(2\left[\rho\right]_{\times}+2\left[\rho\right]_{\times}\right)
\end{align*}
which is equivalent to 
\begin{align}
	\boldsymbol{\mathcal{P}}_{a}\left(R\right)= & 2\frac{1}{1+\left\Vert \rho\right\Vert ^{2}}\left[\rho\right]_{\times}\label{eq:OVERVIEW_ROD_Pa}
\end{align}
Thereby, the vex operator of \eqref{eq:OVERVIEW_ROD_Pa} becomes 
\begin{equation}
	\mathbf{vex}\left(\boldsymbol{\mathcal{P}}_{a}\left(R\right)\right)=2\frac{\rho}{1+\left\Vert \rho\right\Vert ^{2}}\label{eq:OVERVIEW_ROD_VEX_Pa}
\end{equation}
This proves \eqref{eq:OVERVIEW_SO3_lem1_1} in Lemma \ref{Lem:OVERVIEW_1}.
The square norm of \eqref{eq:OVERVIEW_ROD_VEX_Pa} is 
\begin{equation}
	\left\Vert \mathbf{vex}\left(\boldsymbol{\mathcal{P}}_{a}\left(R\right)\right)\right\Vert ^{2}=4\frac{\left\Vert \rho\right\Vert ^{2}}{\left(1+\left\Vert \rho\right\Vert ^{2}\right)^{2}}\label{eq:OVERVIEW_ROD_VEX2_1}
\end{equation}
This shows \eqref{eq:OVERVIEW_SO3_lem1_4} in Lemma \ref{Lem:OVERVIEW_1}.
Also, from \eqref{eq:OVERVIEW_ROD_TR2}, and \eqref{eq:OVERVIEW_ROD_VEX_Pa},
one can find that 
\begin{equation}
	\rho=\frac{\mathbf{vex}\left(\boldsymbol{\mathcal{P}}_{a}\left(R\right)\right)}{2\left(1+\left\Vert \rho\right\Vert ^{2}\right)}\label{eq:OVERVIEW_ROD_rho_1}
\end{equation}
and from \eqref{eq:OVERVIEW_ROD_TR2}, $\left\Vert \rho\right\Vert ^{2}$
can be defined in terms of normalized Euclidean distance $\left\Vert R\right\Vert _{I}$
by\textcolor{red}{{} } 
\begin{equation}
	\left\Vert \rho\right\Vert ^{2}=\frac{\left\Vert R\right\Vert _{I}}{1-\left\Vert R\right\Vert _{I}}\label{eq:OVERVIEW_ROD_rho_2}
\end{equation}
Thereby, from \eqref{eq:OVERVIEW_ROD_rho_1} and \eqref{eq:OVERVIEW_ROD_rho_2},
one has $\rho:\mathbb{SO}\left(3\right)\rightarrow\mathbb{R}^{3}$
(\textbf{\textcolor{blue}{the mapping from $\mathbb{SO}\left(3\right)$
		to $\rho$}}) 
\begin{equation}
	\rho=\frac{\mathbf{vex}\left(\boldsymbol{\mathcal{P}}_{a}\left(R\right)\right)}{2\left(1-\left\Vert R\right\Vert _{I}\right)}\label{eq:OVERVIEW_ROD__SO3}
\end{equation}
And from \eqref{eq:OVERVIEW_ROD_TR2} and \eqref{eq:OVERVIEW_ROD_rho_2},
we have 
\begin{align}
	\left\Vert \mathbf{vex}\left(\boldsymbol{\mathcal{P}}_{a}\left(R\right)\right)\right\Vert ^{2} & =4\times\frac{\left\Vert \rho\right\Vert ^{2}}{1+\left\Vert \rho\right\Vert ^{2}}\times\frac{1}{1+\left\Vert \rho\right\Vert ^{2}}\nonumber \\
	& =4\left\Vert R\right\Vert _{I}\left(1-\left\Vert R\right\Vert _{I}\right)\label{eq:OVERVIEW_ROD_VEX2_2}
\end{align}
This confirms \eqref{eq:OVERVIEW_SO3_lem5_1} in Lemma \ref{Lem:OVERVIEW_5}
and \eqref{eq:OVERVIEW_SO3_lem1_4} in Lemma \ref{Lem:OVERVIEW_1}.
According to the results in \eqref{eq:OVERVIEW_ROD_rho_1}, \eqref{eq:OVERVIEW_att_ang3},
and \eqref{eq:OVERVIEW_att_ang_VEX}, one can find

\begin{align}
	\rho & =\frac{\mathbf{vex}\left(\boldsymbol{\mathcal{P}}_{a}\left(R\right)\right)}{2\left(1-\left\Vert R\right\Vert _{I}\right)}\nonumber \\
	& =\frac{2{\rm cos}\left(\frac{\alpha}{2}\right)\sin\left(\frac{\alpha}{2}\right)u}{2{\rm cos}^{2}\left(\frac{\alpha}{2}\right)}\nonumber \\
	& ={\rm tan}\left(\frac{\alpha}{2}\right)u\label{eq:-7_1}
\end{align}
As such, the related map from angle-axis parameterization to Rodriguez
vector $\rho:\mathbb{R}\times\mathbb{S}^{2}\rightarrow\mathbb{R}^{3}$
(\textbf{\textcolor{blue}{the mapping from $\left(\alpha,u\right)$
		to $\rho$}}) can be defined by \cite{rodrigues1840lois,wilson1901vector,gibbs1961scientific}
\begin{align}
	\rho & ={\rm tan}\left(\frac{\alpha}{2}\right)u\label{eq:OVERVIEW_ROD_from_a_u}
\end{align}
which proves \eqref{eq:OVERVIEW_SO3_lem7_1}. In order to find the
inverse mapping of \eqref{eq:OVERVIEW_ROD_from_a_u}, recall $\left\Vert R\right\Vert _{I}=\sin^{2}\left(\frac{\alpha}{2}\right)$
from \eqref{eq:OVERVIEW_att_ang_RI} and $\left\Vert \rho\right\Vert ^{2}=\frac{\left\Vert R\right\Vert _{I}}{1-\left\Vert R\right\Vert _{I}}$
from \eqref{eq:OVERVIEW_ROD_rho_2}. Hence, the rotation angle $\alpha$
can be obtained by 
\[
\left\Vert \rho\right\Vert ^{2}=\frac{\sin^{2}\left(\frac{\alpha}{2}\right)}{\cos^{2}\left(\frac{\alpha}{2}\right)}=\frac{\sin^{2}\left(\frac{\alpha}{2}\right)}{\cos^{2}\left(\frac{\alpha}{2}\right)}=\tan^{2}\left(\frac{\alpha}{2}\right)
\]
Therefore, $\alpha$ becomes 
\[
\alpha=2\text{ }{\rm tan}^{-1}\left(\left\Vert \rho\right\Vert \right)=2\text{ }{\rm sin}^{-1}\left(\frac{\left\Vert \rho\right\Vert }{\sqrt{1+\left\Vert \rho\right\Vert ^{2}}}\right)
\]
and from \eqref{eq:OVERVIEW_ROD_from_a_u} the unit axis $u\in\mathbb{S}^{2}$
is 
\[
u={\rm cot}\left(\frac{\alpha}{2}\right)\rho
\]
Thus, the related map from Rodriguez vector to angle-axis parameterization
$\rho:\mathbb{R}^{3}\rightarrow\mathbb{R}\times\mathbb{S}^{2}$ (\textbf{\textcolor{blue}{the
		mapping from $\rho$ to $\left(\alpha,u\right)$}}) can be expressed
as \cite{rodrigues1840lois,wilson1901vector,gibbs1961scientific}
\begin{align}
	\alpha & =2\text{ }{\rm tan}^{-1}\left(\left\Vert \rho\right\Vert \right)=2\text{ }{\rm sin}^{-1}\left(\frac{\left\Vert \rho\right\Vert }{\sqrt{1+\left\Vert \rho\right\Vert ^{2}}}\right)\label{eq:OVERVIEW_ROD_2_a_u1}\\
	u & ={\rm cot}\left(\frac{\alpha}{2}\right)\rho\label{eq:OVERVIEW_ROD_2_a_u2}
\end{align}
the above mentioned results show \eqref{eq:OVERVIEW_SO3_lem7_2} and
\eqref{eq:OVERVIEW_SO3_lem7_3}.

\subsection{Attitude measurements}

\noindent Consideringthe body-frame vector ${\rm v}_{i}^{\mathcal{B}}\in\mathbb{R}^{3}$
and the inertial-frame vector ${\rm v}_{i}^{\mathcal{I}}\in\mathbb{R}^{3}$
of the $i$th measurement, the relationship between these two frames
of the $i$th measurement is given by

\begin{equation}
	{\rm v}_{i}^{\mathcal{B}}=R^{\top}{\rm v}_{i}^{\mathcal{I}}\label{eq:OVERVIEW_SO3_dyn1-1}
\end{equation}

\noindent where $R\in\mathbb{SO}\left(3\right)$ is the attitude matrix.
Let $n$ denote the number of body-frame and inertial-frame vectors
available for measurement for $i=1,\ldots,n$. In general, three body-frame
and inertial-frame non-collinear vector measurements are required
for attitude estimation or reconstruction ($n\geq3$). However, if
two non-collinear vectors are available for measurements ($n=2$),
the third vector can be defined by ${\rm v}_{3}^{\mathcal{B}}={\rm v}_{1}^{\mathcal{B}}\times{\rm v}_{2}^{\mathcal{B}}$
and ${\rm v}_{3}^{\mathcal{I}}={\rm v}_{1}^{\mathcal{I}}\times{\rm v}_{2}^{\mathcal{I}}$,
such that the three vectors are non-collinear. Define 
\begin{align}
	M^{\mathcal{B}} & =\sum_{i=1}^{n}\left(\frac{{\rm v}_{i}^{\mathcal{B}}\left({\rm v}_{i}^{\mathcal{B}}\right)^{\top}}{\left\Vert {\rm v}_{i}^{\mathcal{B}}\right\Vert ^{2}}\right)\in\mathbb{R}^{3\times3}\label{eq:OVERVIEW_ROD_MB}\\
	M^{\mathcal{I}} & =\sum_{i=1}^{n}\left(\frac{{\rm v}_{i}^{\mathcal{I}}\left({\rm v}_{i}^{\mathcal{I}}\right)^{\top}}{\left\Vert {\rm v}_{i}^{\mathcal{I}}\right\Vert ^{2}}\right)\in\mathbb{R}^{3\times3}\label{eq:OVERVIEW_ROD_MI}
\end{align}
with $n\geq2$, one of the above-mentioned equations in \eqref{eq:OVERVIEW_ROD_MB}
or \eqref{eq:OVERVIEW_ROD_MI} is normally employed for attitude estimation.
Let $M:=M^{\star}$, such that $\star$ refers to $\mathcal{I}$ or
$\mathcal{B}$ and $\bar{\mathbf{M}}={\rm Tr}\left\{ M\right\} \mathbf{I}_{3}-M$.
Consider three measurements ($n=3$) such that ${\rm Tr}\left\{ M\right\} =3$
and the normalized Euclidean distance of $MR$ $\left\Vert MR\right\Vert _{I}=\frac{1}{4}{\rm Tr}\left\{ M\left(\mathbf{I}_{3}-R\right)\right\} $.
According to angle-axis parameterization in \eqref{eq:SO3PPF_PER_AX_ANG},
one obtains

\begin{align}
	\left\Vert MR\right\Vert _{I} & =\frac{1}{4}{\rm Tr}\left\{ -M\left(\sin(\theta)\left[u\right]_{\times}+\left(1-\cos(\theta)\right)\left[u\right]_{\times}^{2}\right)\right\} \nonumber \\
	& =-\frac{1}{4}{\rm Tr}\left\{ \left(1-\cos(\theta)\right)M\left[u\right]_{\times}^{2}\right\} \label{eq:OVERVIEW_ROD_EXPL_append3}
\end{align}
where ${\rm Tr}\left\{ M\left[u\right]_{\times}\right\} =0$ as given
in identity \eqref{eq:OVERVIEW_PER_Identity_used5}. One has \cite{murray1994mathematicall}
\begin{equation}
	\left\Vert R\right\Vert _{I}=\frac{1}{4}{\rm Tr}\left\{ \mathbf{I}_{3}-R\right\} =\frac{1}{2}\left(1-{\rm cos}\left(\theta\right)\right)={\rm sin}^{2}\left(\frac{\alpha}{2}\right)\label{eq:OVERVIEW_ROD_EXPL_append4}
\end{equation}
and the Rodriguez parameters vector with respect to angle-axis parameterization
is \cite{shuster1993survey} 
\[
u={\rm cot}\left(\frac{\alpha}{2}\right)\rho
\]
From identity \eqref{eq:OVERVIEW_PER_Identity_used4} $\left[u\right]_{\times}^{2}=-u^{\top}u\mathbf{I}_{3}+uu^{\top}$,
the expression in \eqref{eq:OVERVIEW_ROD_EXPL_append3} becomes 
\begin{align*}
	\left\Vert MR\right\Vert _{I} & =\frac{1}{2}\left\Vert R\right\Vert _{I}u^{\top}\bar{\mathbf{M}}u=\frac{1}{2}\left\Vert R\right\Vert _{I}{\rm cot}^{2}\left(\frac{\alpha}{2}\right)\rho^{\top}\bar{\mathbf{M}}\rho
\end{align*}
From \eqref{eq:OVERVIEW_ROD_EXPL_append4}, one can find ${\rm cos}^{2}\left(\frac{\alpha}{2}\right)=1-\left\Vert R\right\Vert _{I}$
which means 
\[
{\rm tan}^{2}\left(\frac{\alpha}{2}\right)=\frac{\left\Vert R\right\Vert _{I}}{1-\left\Vert R\right\Vert _{I}}
\]
Therefore, the normalized Euclidean distance can be expressed with
respect to Rodriguez parameters vector as 
\begin{align}
	\left\Vert MR\right\Vert _{I} & =\frac{1}{2}\left(1-\left\Vert R\right\Vert _{I}\right)\rho^{\top}\bar{\mathbf{M}}\rho=\frac{1}{2}\frac{\rho^{\top}\bar{\mathbf{M}}\rho}{1+\left\Vert \rho\right\Vert ^{2}}\label{eq:OVERVIEW_ROD_EXPL_append_MBR_I}
\end{align}
This proves \eqref{eq:OVERVIEW_SO3_lem3_1} in Lemma \ref{Lem:OVERVIEW_3}.
The anti-symmetric projection operator may be expressed in terms of
Rodriquez parameters vector using the identity in \eqref{eq:OVERVIEW_PER_SO3_PA_6}
and \eqref{eq:OVERVIEW_PER_ROD} by 
\begin{align*}
	\boldsymbol{\mathcal{P}}_{a}\left(MR\right)= & \frac{M\rho\rho^{\top}-\rho\rho^{\top}M+M\left[\rho\right]_{\times}+\left[\rho\right]_{\times}M}{1+\left\Vert \rho\right\Vert ^{2}}\\
	= & \frac{\left[\left({\rm Tr}\left\{ M\right\} \mathbf{I}_{3}-M+\left[\rho\right]_{\times}M\right)\rho\right]_{\times}}{1+\left\Vert \rho\right\Vert ^{2}}
\end{align*}
It follows that the vex operator of the above expression is 
\begin{align}
	\mathcal{\mathbf{vex}}\left(\boldsymbol{\mathcal{P}}_{a}\left(MR\right)\right) & =\frac{\left(\mathbf{I}_{3}+\left[\rho\right]_{\times}\right)^{\top}}{1+\left\Vert \rho\right\Vert ^{2}}\bar{\mathbf{M}}\rho\label{eq:OVERVIEW_ROD_EXPL_append_MBR_VEX}
\end{align}
This shows \eqref{eq:OVERVIEW_SO3_lem3_2} in Lemma \ref{Lem:OVERVIEW_3}.
Also, one can find 
\[
\mathcal{\mathbf{vex}}\left(\boldsymbol{\mathcal{P}}_{a}\left(MR\right)\right)\mathcal{\mathbf{vex}}\left(\boldsymbol{\mathcal{P}}_{a}\left(MR\right)\right)^{\top}=\frac{\left(\mathbf{I}_{3}+\left[\rho\right]_{\times}\right)^{\top}\bar{\mathbf{M}}\rho\rho^{\top}\bar{\mathbf{M}}\left(\mathbf{I}_{3}+\left[\rho\right]_{\times}\right)}{\left(1+\left\Vert \rho\right\Vert ^{2}\right)^{2}}
\]
The 2-norm of \eqref{eq:OVERVIEW_ROD_EXPL_append_MBR_VEX} can be
obtained by 
\begin{align}
	\left\Vert \mathcal{\mathbf{vex}}\left(\boldsymbol{\mathcal{P}}_{a}\left(MR\right)\right)\right\Vert ^{2} & =\frac{\rho^{\top}\bar{\mathbf{M}}\left(\mathbf{I}_{3}-\left[\rho\right]_{\times}^{2}\right)\bar{\mathbf{M}}\rho}{\left(1+\left\Vert \rho\right\Vert ^{2}\right)^{2}}\label{eq:OVERVIEW_ROD_EXPL_append_MBR_VEX2}
\end{align}
This proves \eqref{eq:OVERVIEW_SO3_lem3_3} in Lemma \ref{Lem:OVERVIEW_3}.
Using identity \eqref{eq:OVERVIEW_PER_Identity_used4} $\left[\rho\right]_{\times}^{2}=-\rho^{\top}\rho\mathbf{I}_{3}+\rho\rho^{\top}$,
one obtains 
\begin{align}
	\left\Vert \mathcal{\mathbf{vex}}\left(\boldsymbol{\mathcal{P}}_{a}\left(MR\right)\right)\right\Vert ^{2} & =\frac{\rho^{\top}\bar{\mathbf{M}}\left(\mathbf{I}_{3}-\left[\rho\right]_{\times}^{2}\right)\bar{\mathbf{M}}\rho}{\left(1+\left\Vert \rho\right\Vert ^{2}\right)^{2}}\nonumber \\
	& =\frac{\rho^{\top}\left(\bar{\mathbf{M}}\right)^{2}\rho}{1+\left\Vert \rho\right\Vert ^{2}}-\frac{\left(\rho^{\top}\bar{\mathbf{M}}\rho\right)^{2}}{\left(1+\left\Vert \rho\right\Vert ^{2}\right)^{2}}\nonumber \\
	& \geq\underline{\lambda}\left(\bar{\mathbf{M}}\right)\left(1-\frac{\left\Vert \rho\right\Vert ^{2}}{1+\left\Vert \rho\right\Vert ^{2}}\right)\frac{\rho^{\top}\bar{\mathbf{M}}\rho}{1+||\rho||^{2}}\nonumber \\
	& \geq2\underline{\lambda}\left(1-\left\Vert R\right\Vert _{I}\right)\left\Vert MR\right\Vert _{I}\label{eq:OVERVIEW_ROD_EXPL_append_VEX_MI2}
\end{align}
where $\underline{\lambda}=\underline{\lambda}\left(\bar{\mathbf{M}}\right)$
is the minimum singular value of $\bar{\mathbf{M}}$ and $\left\Vert R\right\Vert _{I}=\left\Vert \rho\right\Vert ^{2}/\left(1+\left\Vert \rho\right\Vert ^{2}\right)$
as defined in \eqref{eq:OVERVIEW_ROD_TR2}. It can be found that 
\begin{align}
	1-\left\Vert R\right\Vert _{I} & =\frac{1}{4}\left(1+{\rm Tr}\left\{ M^{-1}MR\right\} \right)\label{eq:SO3STCH_EXPL_append_rho2}
\end{align}
Therefore, from \eqref{eq:OVERVIEW_ROD_EXPL_append_VEX_MI2} and \eqref{eq:SO3STCH_EXPL_append_rho2}
the following inequality holds \textcolor{red}{{} } 
\begin{align}
	\left\Vert \mathcal{\mathbf{vex}}\left(\boldsymbol{\mathcal{P}}_{a}\left(MR\right)\right)\right\Vert ^{2} & \geq\frac{\underline{\lambda}}{2}\left(1+{\rm Tr}\left\{ M^{-1}MR\right\} \right)\left\Vert MR\right\Vert _{I}\label{eq:OVERVIEW_ROD_EXPL_Proof_VEX_MI2}
\end{align}
which confirms \eqref{eq:OVERVIEW_SO3_lem3_4} in Lemma \ref{Lem:OVERVIEW_1}.
It should be remarked that both $M^{-1}$ and $MR$ can be obtained
by a set of vectorial measurements helping the designer to avoid the
process of attitude reconstruction. For more details visit \cite{hashim2018Conf1,hashim2019SO3Det,mohamed2019filters,hashim2019SE3Det}.

\subsection{Rodriguez vector dynamics}

The kinematic relationship between Rodriguez vector and angular velocity
can be expressed as follows \cite{shuster1993survey} 
\begin{align*}
	\frac{d}{dt}\rho & =\frac{1}{2}\left(\Omega-\Omega\times\rho+\left(\Omega\cdot\rho\right)\rho\right)\\
	& =\frac{1}{2}\left(\Omega+\left[\rho\right]_{\times}\Omega+\rho\rho^{\top}\Omega\right)\\
	& =\frac{1}{2}\left(\mathbf{I}_{3}+\left[\rho\right]_{\times}+\rho\rho^{\top}\right)\Omega
\end{align*}
Hence, the dynamics are governed by{} 
\begin{align}
	\frac{d}{dt}\rho & =\frac{1}{2}\left(\mathbf{I}_{3}+\left[\rho\right]_{\times}+\rho\rho^{\top}\right)\Omega\label{eq:OVERVIEW_ROD_dot}
\end{align}

\subsection{Attitude error and attitude error dynamics in the sense of Rodriguez
	vector}

Consider the attitude dynamics in \textbf{\eqref{eq:OVERVIEW_SO3_dyn}}
\[
\dot{R}=R\left[\Omega\right]_{\times}
\]
consider 
\[
\Omega_{\star}=\Omega+\beta
\]
where $\Omega_{\star}$ could represent the uncertain measurements
of $\Omega$ and $\beta$ could be considered as an unknown variable.
Let us introduce desired attitude dynamics (estimator attitude dynamics
) by 
\begin{equation}
	\dot{R}_{\star}=R_{\star}\left[\Omega_{\star}-\hat{\beta}\right]_{\times}\label{eq:OVERVIEW_SO3_dyn_EST-1}
\end{equation}
where $\hat{\beta}$ is the estimate of $\beta$. Let the error in
attitude be given by 
\begin{equation}
	\tilde{R}=RR_{\star}^{\top}\label{eq:OVERVIEW_SO3_Error-1}
\end{equation}
and define the error in $\beta$ by $\tilde{\beta}=\beta-\hat{\beta}$.
Hence, the dynamics of the attitude error can be found to be 
\begin{align}
	\dot{\tilde{R}} & =\dot{R}R_{\star}^{\top}+R\dot{R}_{\star}^{\top}\nonumber \\
	& =R\left[\Omega\right]_{\times}R_{\star}^{\top}-R\left[\Omega+\beta-\hat{\beta}\right]_{\times}R_{\star}^{\top}\nonumber \\
	& =-R\left[\tilde{\beta}\right]_{\times}R_{\star}^{\top}\label{eq:OVERVIEW_SO3_dot_Error-1}
\end{align}
where $\left[\Omega+\beta-\hat{\beta}\right]_{\times}^{\top}=-\left[\Omega+\beta-\hat{\beta}\right]_{\times}$
as defined in \eqref{eq:OVERVIEW_PER_SO3_3}. One can find 
\begin{align}
	\dot{\tilde{R}} & =-R\left[\tilde{\beta}\right]_{\times}R_{\star}^{\top}\nonumber \\
	& =-R\mathbf{I}_{3}\left[\tilde{\beta}\right]_{\times}\mathbf{I}_{3}R_{\star}^{\top}\nonumber \\
	& =-RR_{\star}^{\top}R_{\star}\left[\tilde{\beta}\right]_{\times}R_{\star}^{\top}R_{\star}R_{\star}^{\top}\nonumber \\
	& =-RR_{\star}^{\top}\left[R_{\star}\tilde{\beta}\right]_{\times}\nonumber \\
	& =-\tilde{R}\left[R_{\star}\tilde{\beta}\right]_{\times}\label{eq:OVERVIEW_SO3_dot_Error-1-1}
\end{align}
with $\left[R_{\star}\tilde{\beta}\right]_{\times}=R_{\star}\left[\tilde{\beta}\right]_{\times}R_{\star}^{\top}$
being given in identity \eqref{eq:OVERVIEW_PER_Identity_used1}. Thus,
in view of \textbf{\eqref{eq:OVERVIEW_SO3_dyn} }and \textbf{\eqref{eq:OVERVIEW_ROD_dot}},
one can write \textbf{\eqref{eq:OVERVIEW_SO3_dot_Error-1-1}} in terms
of Rodriguez vector dynamics as 
\begin{align}
	\frac{d}{dt}\tilde{\rho} & =-\frac{1}{2}\left(\mathbf{I}_{3}+\left[\tilde{\rho}\right]_{\times}+\tilde{\rho}\tilde{\rho}^{\top}\right)R_{\star}\beta\label{eq:OVERVIEW_ROD_dot_ERR}
\end{align}
where $\tilde{\rho}$ is the error in Rodriguez vector associated
with $\tilde{R}$. For more detailed derivations visit \cite{hashim2018SO3Stochastic,mohamed2019filters}.
The objective of attitude filter/control is to drive the error in
Rodriguez vector to zero ($\tilde{\rho}\rightarrow0$). Driving $\tilde{\rho}\rightarrow0$
implies that $\tilde{R}\rightarrow\mathbf{I}_{3}$, since we have
\begin{align}
	\mathcal{R}_{\tilde{\rho}}\left(\tilde{\rho}\right) & =\frac{1}{1+\left\Vert \tilde{\rho}\right\Vert ^{2}}\left(\left(1-\left\Vert \tilde{\rho}\right\Vert ^{2}\right)\mathbf{I}_{3}+2\tilde{\rho}\tilde{\rho}^{\top}+2\left[\tilde{\rho}\right]_{\times}\right)\nonumber \\
	\tilde{\rho}=0 & \Leftrightarrow\mathcal{R}_{\tilde{\rho}}\left(\tilde{\rho}\right)=\mathbf{I}_{3}\label{eq:-33}
\end{align}

\subsection{Problems of Rodriguez vector parameterization: }

Although Rodriguez vector provides a unique representation of the
attitude, the Rodriguez vector and the modified Rodriguez vector are
not defined for 180° of rotation. To be more specific, the parameters
of the Rodriguez vector are not defined for any of the following three
rotational matrices 
\begin{align}
	\text{at }R & =\left[\begin{array}{ccc}
		1 & 0 & 0\\
		0 & -1 & 0\\
		0 & 0 & -1
	\end{array}\right]\in\mathbb{SO}\left(3\right)\Rightarrow\hspace{1em}\rho=\left[\infty,\infty,\infty\right]^{\top}\notin\mathbb{R}^{3}\label{eq:OVERVIEW_ROD_Prob1}\\
	\text{at }R & =\left[\begin{array}{ccc}
		-1 & 0 & 0\\
		0 & 1 & 0\\
		0 & 0 & -1
	\end{array}\right]\in\mathbb{SO}\left(3\right)\Rightarrow\hspace{1em}\rho=\left[\infty,\infty,\infty\right]^{\top}\notin\mathbb{R}^{3}\label{eq:OVERVIEW_ROD_Prob2}\\
	\text{at }R & =\left[\begin{array}{ccc}
		-1 & 0 & 0\\
		0 & -1 & 0\\
		0 & 0 & 1
	\end{array}\right]\in\mathbb{SO}\left(3\right)\Rightarrow\hspace{1em}\rho=\left[\infty,\infty,\infty\right]^{\top}\notin\mathbb{R}^{3}\label{eq:OVERVIEW_ROD_Prob3}
\end{align}
which means that as $\alpha\rightarrow\pi$, the Rodriguez vector
$\rho\rightarrow\infty$. Therefore, the mapping from $\mathbb{SO}\left(3\right)$
to Rodriguez vector cannot be achieved for any rotational matrix in
\eqref{eq:OVERVIEW_ROD_Prob1}, \eqref{eq:OVERVIEW_ROD_Prob2}, or
\eqref{eq:OVERVIEW_ROD_Prob3}. In this regard, for $R$ equivalent
to \eqref{eq:OVERVIEW_ROD_Prob1}, \eqref{eq:OVERVIEW_ROD_Prob2},
or \eqref{eq:OVERVIEW_ROD_Prob3}, continuous control laws cannot
be globally defined using the three-parameter vector.

\noindent %
\noindent\makebox[1\linewidth]{%
	\rule{0.6\textwidth}{1.4pt}%
}

\section{Unit-quaternion\label{sec:Unit-quaternion}}

\noindent Unit-quaternion has proven to be an effective tool for the
tracking control of UAVs, such as \cite{lizarralde1996attitude,mayhew2011quaternion,joshi1995robust}.
Unit-quaternion also showed impressive results in attitude estimators.
The unit-quaternion vector has been employed in deterministic attitude
filters, for instance \cite{mahony2008nonlinear,euston2008complementary}.
However, unit-quaternion has been used more extensively in Gaussian
attitude filters, namely Kalman filter \cite{choukroun2006novel},
extended Kalman filter \cite{lefferts1982kalman,marins2001extended},
multiplicative extended Kalman filter \cite{markley2003attitude},
unscented Kalman filter \cite{cheon2007unscented}, and invariant
extended Kalman filter \cite{barrau2015intrinsic}. The main advantage
of the unit-quaternion vector is that it gives nonsingular representation
of the attitude. The vector $Q$ is said to be a unit-quaternion such
that $Q\in\mathbb{S}^{3}$ if the following two conditions are met
\cite{cayley1845xiii,shuster1993survey} 
\begin{enumerate}
	\item $Q\in\mathbb{R}^{4}$. 
	\item $\left\Vert Q\right\Vert =1$. 
\end{enumerate}
\noindent The unit-quaternion is a four-element representation of
the attitude, and these four elements do not have intuitive physical
meanings. The unit-quaternion is denoted by 
\begin{align}
	Q & =\left[\begin{array}{cccc}
		q_{0} & q_{1} & q_{2} & q_{3}\end{array}\right]^{\top}=\left[\begin{array}{c}
		q_{0}\\
		q
	\end{array}\right]\in\mathbb{S}^{3}\label{eq:OVERVIEW_Q}
\end{align}
where $q=\left[q_{1},q_{2},q_{3}\right]^{\top}\in\mathbb{R}^{3}$
and $q_{0}\in\mathbb{R}$. The unit-quaternion vector $Q$ is defined
by 
\[
\mathbb{S}^{3}=\left\{ \left.Q\in\mathbb{R}^{4}\right|\left\Vert Q\right\Vert =\sqrt{q_{0}^{2}+q_{1}^{2}+q_{2}^{2}+q_{3}^{2}}=1\right\} 
\]
where $Q$ is non-Euclidean, lies in the three-sphere ($\mathbb{S}^{3}$)
and is given by 
\[
Q=q_{0}+iq_{1}+jq_{2}+kq_{3}
\]
with $i$, $j$, and $k$ being hyper-imaginary numbers which satisfy
the following rules

\noindent 
\begin{align*}
	i^{2} & =j^{2}=k^{2}=ijk=-1\\
	ij & =i\times j=-j\times i=-ji=k,\hspace{1em}ji=-k\\
	jk & =j\times k=-k\times j=-kj=i,\hspace{1em}kj=-i\\
	ki & =k\times i=-i\times k=-ik=j,\hspace{1em}ik=-j
\end{align*}
where multiplication follows the \textquotedbl{} natural order\textquotedbl{}
convention \cite{shuster1993survey,breckenridge1999quaternions}.
Quaternion multiplication is in general associative, that is, 
\begin{equation}
	Q_{1}\odot Q_{2}\odot Q_{3}=\left(Q_{1}\odot Q_{2}\right)\odot Q_{3}=Q_{1}\odot\left(Q_{2}\odot Q_{3}\right)\label{eq:}
\end{equation}
where $Q_{1},Q_{2},Q_{3}\in\mathbb{S}^{3}$ and $\odot$ denotes multiplication
operator between two quaternion vectors. Quaternion multiplication
is not commutative such that 
\begin{equation}
	Q_{1}\odot Q_{2}\ne Q_{2}\odot Q_{1},\hspace{1em}\forall Q_{1},Q_{2}\in\mathbb{S}^{3}\label{eq:-1}
\end{equation}
Hence 
\begin{align*}
	P & =p_{0}+p\\
	& =p_{0}+ip_{1}+jp_{2}+kp_{3}\\
	M & =m_{0}+m\\
	& =m_{0}+im_{1}+jm_{2}+km_{3}\\
	PM & =\left(p_{0}+ip_{1}+jp_{2}+kp_{3}\right)\left(m_{0}+im_{1}+jm_{2}+km_{3}\right)\\
	& =p_{0}m_{0}-p\cdot m+p_{0}m+m_{0}p+p\times m
\end{align*}
The inverse of the rotation is defined by the complex conjugate or
inverse of a unit-quaternion, which is given by 
\begin{equation}
	Q^{*}=Q^{-1}=\left[\begin{array}{c}
		q_{0}\\
		-q
	\end{array}\right]\in\mathbb{S}^{3}\label{eq:OVERVIEW_Q_inv}
\end{equation}

\subsection{Quaternion multiplication}

\noindent Analogously to linear matrix multiplication of rotation
matrices, the composition of successive rotations represented by unit-quaternion
is obtained by the distributive and associative, but not commutative,
quaternion multiplication. To define this operation, consider two
unit-quaternion vectors 
\[
Q_{1}=\left[\begin{array}{c}
	q_{01}\\
	q_{1}
\end{array}\right],\hspace{1em}Q_{2}=\left[\begin{array}{c}
	q_{02}\\
	q_{2}
\end{array}\right],\hspace{1em}\forall Q_{1},Q_{2}\in\mathbb{S}^{3}
\]
The quaternion product between $Q_{1}$ and $Q_{2}$, denoted by $Q_{3}\in\mathbb{S}^{3}$,
is given by 
\begin{align}
	Q_{3} & =Q_{1}\odot Q_{2}=\left[\begin{array}{c}
		q_{01}\\
		q_{1}
	\end{array}\right]\odot\left[\begin{array}{c}
		q_{02}\\
		q_{2}
	\end{array}\right]\nonumber \\
	& =\left[\begin{array}{c}
		q_{01}q_{02}-q_{1}^{\top}q_{2}\\
		q_{01}q_{2}+q_{02}q_{1}+\left[q_{1}\right]_{\times}q_{2}
	\end{array}\right]\label{eq:OVERVIEW_Q_Mul}
\end{align}
The neutral element of the unit-quaternion is denoted by $Q_{{\rm I}}\in\mathbb{S}^{3}$,
which is defined by 
\begin{align}
	Q_{{\rm I}}=Q\odot Q^{*} & =\left[\begin{array}{c}
		q_{0}\\
		q
	\end{array}\right]\odot\left[\begin{array}{c}
		q_{0}\\
		-q
	\end{array}\right]\nonumber \\
	& =\left[\begin{array}{c}
		q_{0}q_{0}+q^{\top}q\\
		-q_{0}q+q_{0}q-\left[q\right]_{\times}q
	\end{array}\right]\nonumber \\
	& =\left[\begin{array}{c}
		q_{0}^{2}+\left\Vert q\right\Vert ^{2}\\
		-\left[q\right]_{\times}q
	\end{array}\right]\nonumber \\
	& =\left[\begin{array}{c}
		1\\
		0\\
		0\\
		0
	\end{array}\right]\label{eq:OVERVIEW_Q_Ident}
\end{align}
where $\left[q\right]_{\times}q=\left[0,0,0\right]^{\top}$ as given
in \eqref{eq:OVERVIEW_PER_Identity_used2}. The inverse of quaternion
multiplication is equivalent to 
\[
\left(Q_{1}\odot Q_{2}\right)^{-1}=Q_{2}^{-1}\odot Q_{1}^{-1},\hspace{1em}\forall Q_{1},Q_{2}\in\mathbb{S}^{3}
\]

\subsection{Mapping: From/To unit-quaternion to/from other representations}

\noindent Using the quaternion product, the unit-quaternion may also
be utilized to give the coordinates of a vector in multiple frames
of reference. Actually, having the property ${\rm v}^{\mathcal{B}}=\mathcal{R}_{Q}^{\top}\left(Q\right){\rm v}^{\mathcal{I}}$,
the vector ${\rm v}^{\mathcal{B}}$ may be obtained through the quaternion
product by the following operation \cite{stuelpnagel1964parametrization,shuster1993survey}:

\noindent 
\begin{align}
	\left[\begin{array}{c}
		0\\
		{\rm v}^{\mathcal{B}}
	\end{array}\right] & =\left[\begin{array}{c}
		q_{0}\\
		-q
	\end{array}\right]\odot\left[\begin{array}{c}
		0\\
		{\rm v}^{\mathcal{I}}
	\end{array}\right]\odot\left[\begin{array}{c}
		q_{0}\\
		q
	\end{array}\right]\nonumber \\
	& =\left[\begin{array}{c}
		q_{0}q^{\top}{\rm v}^{\mathcal{I}}-q_{0}x_{1}^{\top}q-\left({\rm v}^{\mathcal{I}}\right)^{\top}\left(q\times q\right)\\
		qq^{\top}{\rm v}^{\mathcal{I}}+q_{0}^{2}{\rm v}^{\mathcal{I}}-q_{0}q\times{\rm v}^{\mathcal{I}}+\left(q_{0}{\rm v}^{\mathcal{I}}+\left[-q\right]_{\times}{\rm v}^{\mathcal{I}}\right)\times q
	\end{array}\right]\nonumber \\
	& =\left[\begin{array}{c}
		q_{0}q^{\top}{\rm v}^{\mathcal{I}}-q_{0}\left({\rm v}^{\mathcal{I}}\right)^{\top}q-\left({\rm v}^{\mathcal{I}}\right)^{\top}\left(q\times q\right)\\
		qq^{\top}{\rm v}^{\mathcal{I}}+q_{0}^{2}{\rm v}^{\mathcal{I}}-q_{0}q\times{\rm v}^{\mathcal{I}}+q_{0}{\rm v}^{\mathcal{I}}\times q+\left[-q\right]_{\times}{\rm v}^{\mathcal{I}}\times q
	\end{array}\right]\nonumber \\
	& =\left[\begin{array}{c}
		0\\
		qq^{\top}{\rm v}^{\mathcal{I}}+q_{0}^{2}{\rm v}^{\mathcal{I}}-2q_{0}q\times{\rm v}^{\mathcal{I}}-q\times{\rm v}^{\mathcal{I}}\times q
	\end{array}\right]\nonumber \\
	& =\left[\begin{array}{c}
		0\\
		qq^{\top}{\rm v}^{\mathcal{I}}+q_{0}^{2}{\rm v}^{\mathcal{I}}-2q_{0}q\times{\rm v}^{\mathcal{I}}-\left(\left(q\cdot q\right){\rm v}^{\mathcal{I}}-\left(q\cdot{\rm v}^{\mathcal{I}}\right)q\right)
	\end{array}\right]\nonumber \\
	& =\left[\begin{array}{c}
		0\\
		qq^{\top}{\rm v}^{\mathcal{I}}+q_{0}^{2}{\rm v}^{\mathcal{I}}-2q_{0}q\times{\rm v}^{\mathcal{I}}-q^{\top}q{\rm v}^{\mathcal{I}}+q^{\top}{\rm v}^{\mathcal{I}}q
	\end{array}\right]\nonumber \\
	& =\left[\begin{array}{c}
		0\\
		\left(q_{0}^{2}-\left\Vert q\right\Vert ^{2}\right){\rm v}^{\mathcal{I}}-2q_{0}\left[q\right]_{\times}{\rm v}^{\mathcal{I}}+2qq^{\top}{\rm v}^{\mathcal{I}}
	\end{array}\right]\nonumber \\
	& =\left[\begin{array}{c}
		0_{1\times3}\\
		\left(q_{0}^{2}-\left\Vert q\right\Vert ^{2}\right)\mathbf{I}_{3}+2qq^{\top}-2q_{0}\left[q\right]_{\times}
	\end{array}\right]{\rm v}^{\mathcal{I}}\label{eq:OVERVIEW_Q_4}
\end{align}
The translation from body-frame (${\rm v}^{\mathcal{B}}$) to inertial-frame
(${\rm v}^{\mathcal{I}}$) is given in \eqref{eq:OVERVIEW_SO3_dyn1}
by ${\rm v}^{\mathcal{B}}=R^{\top}{\rm v}^{\mathcal{I}}$. Thus the
result in \eqref{eq:OVERVIEW_Q_4} indicates that 
\begin{equation}
	\left[\begin{array}{c}
		0\\
		{\rm v}^{\mathcal{B}}
	\end{array}\right]=\left[\begin{array}{c}
		0_{1\times3}\\
		\mathcal{R}_{Q}^{\top}\left(Q\right)
	\end{array}\right]{\rm v}^{\mathcal{I}}\label{eq:OVERVIEW_Q_3}
\end{equation}
with 
\[
\mathcal{R}_{Q}^{\top}\left(Q\right)=\left(q_{0}^{2}-\left\Vert q\right\Vert ^{2}\right)\mathbf{I}_{3}+2qq^{\top}-2q_{0}\left[q\right]_{\times}
\]
Therefore, the coordinate of a moving frame with respect to a reference-frame
(\textbf{\textcolor{blue}{the mapping from $Q$ to $\mathbb{SO}\left(3\right)$}})
$\mathcal{R}_{Q}:\mathbb{S}^{3}\rightarrow\mathbb{SO}\left(3\right)$
is governed by 
\begin{align}
	\mathcal{R}_{Q}\left(Q\right) & =\left(q_{0}^{2}-\left\Vert q\right\Vert ^{2}\right)\mathbf{I}_{3}+2qq^{\top}+2q_{0}\left[q\right]_{\times}\nonumber \\
	& =\mathbf{I}_{3}+2q_{0}\left[q\right]_{\times}+2\left[q\right]_{\times}^{2}\label{eq:OVERVIEW_Q_R}
\end{align}
where $\mathcal{R}_{Q}\in\mathbb{SO}\left(3\right)$ and the attitude
can be represented in terms of unit-quaternion components as 
\begin{align}
	\mathcal{R}_{Q}\left(Q\right) & =\left[\begin{array}{ccc}
		1-2\left(q_{2}^{2}+q_{3}^{2}\right) & 2\left(q_{1}q_{2}-q_{0}q_{3}\right) & 2\left(q_{1}q_{3}+q_{0}q_{2}\right)\\
		2\left(q_{2}q_{1}+q_{0}q_{3}\right) & 1-2\left(q_{1}^{2}+q_{3}^{2}\right) & 2\left(q_{2}q_{3}-q_{0}q_{1}\right)\\
		2\left(q_{3}q_{1}-q_{0}q_{2}\right) & 2\left(q_{3}q_{2}+q_{0}q_{1}\right) & 1-2\left(q_{1}^{2}+q_{2}^{2}\right)
	\end{array}\right]\in\mathbb{SO}\left(3\right)\label{eq:OVERVIEW_Q_R1}
\end{align}
In order to find the normalized Euclidean distance in terms of unit-quaternion
components, it is necessary to find the trace of the rotational matrix
in \eqref{eq:OVERVIEW_Q_R1}. One can find that the trace of $\mathcal{R}_{Q}\left(Q\right)$
is 
\begin{align}
	{\rm Tr}\left\{ \mathcal{R}_{Q}\left(Q\right)\right\}  & =1-2\left(q_{2}^{2}+q_{3}^{2}\right)+1-2\left(q_{1}^{2}+q_{3}^{2}\right)+1-2\left(q_{1}^{2}+q_{2}^{2}\right)\nonumber \\
	& =3-4\left(1-q_{0}^{2}\right)\nonumber \\
	& =4q_{0}^{2}-1\label{eq:OVERVIEW_Q_Tr}
\end{align}
Hence, from \eqref{eq:OVERVIEW_SO3_Ecul_Dist} and \eqref{eq:OVERVIEW_Q_Tr},
it can be found that the normalized Euclidean distance in terms of
unit-quaternion is 
\begin{align*}
	\left\Vert R\right\Vert _{I} & =\frac{1}{4}{\rm Tr}\left\{ \mathbf{I}_{3}-\mathcal{R}_{Q}\left(Q\right)\right\} 
\end{align*}
which is equivalent to

\noindent 
\begin{align}
	\left\Vert R\right\Vert _{I} & =1-q_{0}^{2}\label{eq:OVERVIEW_Q_TR2}
\end{align}
The result in \eqref{eq:OVERVIEW_Q_TR2} proves \eqref{eq:OVERVIEW_SO3_lem2_2}.
Considering the identity in \eqref{eq:OVERVIEW_PER_SO3_PA_6} and
the mapping from quaternion $Q$ to $\mathbb{SO}\left(3\right)$ in
\eqref{eq:OVERVIEW_Q_R1}, the anti-symmetric projection operator
$\boldsymbol{\mathcal{P}}_{a}\left(R\right)$ can be defined in terms
of unit-quaternion as follows 
\begin{equation}
	\boldsymbol{\mathcal{P}}_{a}\left(\mathcal{R}_{Q}\left(Q\right)\right)=2q_{0}\left[\begin{array}{ccc}
		0 & -q_{3} & q_{2}\\
		q_{3} & 0 & -q_{1}\\
		-q_{2} & q_{1} & 0
	\end{array}\right]\label{eq:OVERVIEW_Q_Pa}
\end{equation}
Consequently, one can find the vex operator with respect to unit-quaternion:
\begin{equation}
	\mathbf{vex}\left(\boldsymbol{\mathcal{P}}_{a}\left(\mathcal{R}_{Q}\left(Q\right)\right)\right)=2q_{0}q\in\mathbb{R}^{3}\label{eq:OVERVIEW_Q_VEX}
\end{equation}
The result in \eqref{eq:OVERVIEW_Q_VEX} justifies \eqref{eq:OVERVIEW_SO3_lem2_1}.
Thus, the norm of the vex operator with regards to unit-quaternion
is equivalent to 
\begin{equation}
	\left\Vert \mathbf{vex}\left(\boldsymbol{\mathcal{P}}_{a}\left(\mathcal{R}_{Q}\left(Q\right)\right)\right)\right\Vert ^{2}=4q_{0}^{2}\left\Vert q\right\Vert ^{2}\in\mathbb{R}\label{eq:OVERVIEW_Q_VEX2}
\end{equation}
This proves \eqref{eq:OVERVIEW_SO3_lem2_3}. Recalling $\mathcal{R}_{Q}\left(Q\right)$
from \eqref{eq:OVERVIEW_Q_R} and \eqref{eq:OVERVIEW_Q_R1}, one can
find $\mathcal{R}_{Q}:\mathbb{SO}\left(3\right)\rightarrow\mathbb{S}^{3}$
(\textbf{\textcolor{blue}{the mapping from $\mathbb{SO}\left(3\right)$
		to $Q$}}) \cite{shuster1993survey}

\noindent 
\begin{align}
	q_{0} & =\frac{1}{2}\sqrt{1+R_{\left(1,1\right)}+R_{\left(2,2\right)}+R_{\left(3,3\right)}}\nonumber \\
	q_{1} & =\frac{1}{4q_{0}}\left(R_{\left(3,2\right)}-R_{\left(2,3\right)}\right)\nonumber \\
	q_{2} & =\frac{1}{4q_{0}}\left(R_{\left(1,3\right)}-R_{\left(3,1\right)}\right)\nonumber \\
	q_{3} & =\frac{1}{4q_{0}}\left(R_{\left(2,1\right)}-R_{\left(1,2\right)}\right)\label{eq:OVERVIEW_Q_comp_1}
\end{align}

\noindent Or 
\begin{align}
	q_{1} & =\pm\frac{1}{2}\sqrt{1+R_{\left(1,1\right)}-R_{\left(2,2\right)}-R_{\left(3,3\right)}}\nonumber \\
	q_{2} & =\frac{1}{4q_{1}}\left(R_{\left(1,2\right)}+R_{\left(2,1\right)}\right)\nonumber \\
	q_{3} & =\frac{1}{4q_{1}}\left(R_{\left(1,3\right)}+R_{\left(3,1\right)}\right)\nonumber \\
	q_{0} & =\frac{1}{4q_{1}}\left(R_{\left(3,2\right)}-R_{\left(2,3\right)}\right)\label{eq:OVERVIEW_Q_comp_2}
\end{align}
Or 
\begin{align}
	q_{2} & =\pm\frac{1}{2}\sqrt{1-R_{\left(1,1\right)}+R_{\left(2,2\right)}-R_{\left(3,3\right)}}\nonumber \\
	q_{1} & =\frac{1}{4q_{2}}\left(R_{\left(1,2\right)}+R_{\left(2,1\right)}\right)\nonumber \\
	q_{3} & =\frac{1}{4q_{2}}\left(R_{\left(2,3\right)}+R_{\left(3,2\right)}\right)\nonumber \\
	q_{0} & =\frac{1}{4q_{2}}\left(R_{\left(1,3\right)}-R_{\left(3,1\right)}\right)\label{eq:OVERVIEW_Q_comp_3}
\end{align}
Or 
\begin{align}
	q_{3} & =\pm\frac{1}{2}\sqrt{1-R_{\left(1,1\right)}-R_{\left(2,2\right)}+R_{\left(3,3\right)}}\nonumber \\
	q_{1} & =\frac{1}{4q_{3}}\left(R_{\left(1,3\right)}+R_{\left(3,1\right)}\right)\nonumber \\
	q_{2} & =\frac{1}{4q_{3}}\left(R_{\left(2,3\right)}+R_{\left(3,2\right)}\right)\nonumber \\
	q_{0} & =\frac{1}{4q_{3}}\left(R_{\left(2,1\right)}-R_{\left(1,2\right)}\right)\label{eq:OVERVIEW_Q_comp_4}
\end{align}
For a given $R$, at least one of $q_{0}$, $q_{1}$, $q_{2}$, and
$q_{3}$ is non-zero at any time instant, while the singularity can
always be avoided through the proper choice of one of the above-mentioned
formulas \eqref{eq:OVERVIEW_Q_comp_1}, \eqref{eq:OVERVIEW_Q_comp_2},
\eqref{eq:OVERVIEW_Q_comp_3}, or \eqref{eq:OVERVIEW_Q_comp_4}. The
unit-quaternion is often considered to be an angle-axis representation.
Indeed, from \eqref{eq:OVERVIEW_att_ang1}, \eqref{eq:OVERVIEW_att_ang2},
and \eqref{eq:OVERVIEW_att_ang_u}, the rotation by an angle $\alpha\in\mathbb{R}$
about an arbitrary unit-length vector $u\in\mathbb{S}^{2}$ can be
described by the unit-quaternion $Q:\mathbb{R}\times\mathbb{S}^{2}\rightarrow\mathbb{S}^{3}$
(\textbf{\textcolor{blue}{the mapping from $\left(\alpha,u\right)$
		to $Q$}}) 
\begin{align}
	Q & =\left[\begin{array}{c}
		\cos\left(\alpha/2\right)\\
		u\sin\left(\alpha/2\right)
	\end{array}\right]\nonumber \\
	& =\left[\begin{array}{cccc}
		\cos\left(\alpha/2\right) & u_{1}\sin\left(\alpha/2\right) & u_{2}\sin\left(\alpha/2\right) & u_{3}\sin\left(\alpha/2\right)\end{array}\right]^{\top}\label{eq:OVERVIEW_Q_ang_2_Q}
\end{align}
And from \eqref{eq:OVERVIEW_Q_ang_2_Q} the mapping from angle-axis
representation to unit-quaternion $Q:\mathbb{S}^{3}\rightarrow\mathbb{R}\times\mathbb{S}^{2}$
(\textbf{\textcolor{blue}{the mapping from $Q$ to $\left(\alpha,u\right)$}})
can be accomplished by 
\begin{align}
	\alpha & =2{\rm cos}^{-1}\left(q_{0}\right)\label{eq:OVERVIEW_Q_Q_2_ang1}\\
	u & =\frac{1}{\sin\left(\alpha/2\right)}q\label{eq:OVERVIEW_Q_Q_2_ang2}
\end{align}
Also, unit-quaternion can be mapped to Euclidean vector in the sense
of Euler angles. From, \eqref{eq:OVERVIEW_EUL_R}, \eqref{eq:SO3_EUL},
and \eqref{eq:OVERVIEW_Q_R1}, Euler angles can be defined in terms
of unit-quaternion components $\xi:\mathbb{S}^{3}\rightarrow\mathbb{R}^{3}$
(\textbf{\textcolor{blue}{the mapping from $Q$ to $\xi$}}) by 
\begin{equation}
	\left[\begin{array}{c}
		\phi\\
		\theta\\
		\psi
	\end{array}\right]=\left[\begin{array}{c}
		\arctan\left(\frac{2\left(q_{3}q_{2}+q_{0}q_{1}\right)}{1-2\left(q_{1}^{2}+q_{2}^{2}\right)}\right)\\
		\arcsin\left(2\left(q_{0}q_{2}-q_{3}q_{1}\right)\right)\\
		\arctan\left(\frac{2\left(q_{2}q_{1}+q_{0}q_{3}\right)}{1-2\left(q_{2}^{2}+q_{3}^{2}\right)}\right)
	\end{array}\right]\label{eq:-28}
\end{equation}
The mapping from unit-quaternion to Rodriguez vector and vice-versa
is proven in Lemma \ref{Lem:OVERVIEW_4}. Let us get back to quaternion
multiplication, let $Q_{1},Q_{2}\in\mathbb{S}^{3}$ 
\[
Q_{1}\odot Q_{2}=\left[\begin{array}{c}
	q_{01}\\
	q_{1}
\end{array}\right]\odot\left[\begin{array}{c}
	q_{02}\\
	q_{2}
\end{array}\right]=\left[\begin{array}{c}
	q_{01}q_{02}-q_{1}^{\top}q_{2}\\
	q_{01}q_{2}+q_{02}q_{1}+\left[q_{1}\right]_{\times}q_{2}
\end{array}\right]
\]
hence, one can easily find \textcolor{red}{{} } 
\begin{equation}
	\mathcal{R}_{Q}\left(Q_{1}\right)\mathcal{R}_{Q}\left(Q_{2}\right)=\mathcal{R}_{Q}\left(Q_{1}\odot Q_{2}\right)\label{eq:OVERVIEW_Q_RR}
\end{equation}
One can find that for $X\in\mathbb{R}^{4}$ that

\begin{equation}
	Q\odot X^{*}\odot X\odot Q^{*}=X^{*}\odot X,\hspace{1em}Q\in\mathbb{S}^{3},X\in\mathbb{R}^{4}\label{eq:OVERVIEW_Q_Identity1}
\end{equation}
Let $Q_{1}=\left[q_{01},q_{1}^{\top}\right]^{\top}=2Q\odot X^{*}$
and $Q_{2}=\left[q_{02},q_{2}^{\top}\right]^{\top}=2X\odot Q^{*}$,
for all $Q\in\mathbb{S}^{3}$ and $X\in\mathbb{R}^{4}$, then 
\begin{equation}
	q_{1}=-q_{2},\hspace{1em}q_{1},q_{2}\in\mathbb{R}^{3}\label{eq:OVERVIEW_Q_Identity2}
\end{equation}

\subsection{Attitude measurements}

\noindent Consider the body-frame vector ${\rm v}_{i}^{\mathcal{B}}\in\mathbb{R}^{3}$
and the inertial-frame vector ${\rm v}_{i}^{\mathcal{I}}\in\mathbb{R}^{3}$.
The relationship between a body-frame and an inertial-frame of the
$i$th measurement is given in \eqref{eq:OVERVIEW_SO3_dyn1}

\begin{equation}
	{\rm v}_{i}^{\mathcal{B}}=R^{\top}{\rm v}_{i}^{\mathcal{I}}\label{eq:OVERVIEW_SO3_dyn1-1-1}
\end{equation}

\noindent where $R\in\mathbb{SO}\left(3\right)$ is the attitude matrix
for $i=1,\ldots,n$. The body-frame vector measurement in \eqref{eq:OVERVIEW_SO3_dyn1-1-1}
in terms of unit-quaternion is equivalent to 
\begin{align}
	\left[\begin{array}{c}
		0\\
		{\rm v}_{i}^{\mathcal{B}}
	\end{array}\right] & =Q^{-1}\odot\left[\begin{array}{c}
		0\\
		{\rm v}_{i}^{\mathcal{I}}
	\end{array}\right]\odot Q\label{eq:OVERVIEW_Q_vB}\\
	\overline{{\rm v}}_{i}^{\mathcal{B}} & =Q^{-1}\odot\overline{{\rm v}}_{i}^{\mathcal{I}}\odot Q\label{eq:OVERVIEW_Q_vB_bar}
\end{align}
where 
\[
\overline{{\rm v}}_{i}^{\mathcal{I}}=\left[\begin{array}{c}
	0\\
	{\rm v}_{i}^{\mathcal{I}}
\end{array}\right]\in\mathbb{R}^{4},\hspace{1em}\overline{{\rm v}}_{i}^{\mathcal{B}}=\left[\begin{array}{c}
	0\\
	{\rm v}_{i}^{\mathcal{B}}
\end{array}\right]\in\mathbb{R}^{4}
\]

\subsection{Unit-quaternion attitude dynamics and measurements}

Let $\Omega=\left[\Omega_{x},\Omega_{y},\Omega_{z}\right]^{\top}\in\mathbb{R}^{3}$
be the angular velocity defined relative to the body-frame $\Omega\in\left\{ \mathcal{B}\right\} $,
and let $Q\in\mathbb{S}^{3}$ be a unit-quaternion vector. Consider
the following representations 
\[
\bar{\Omega}=\left[\begin{array}{c}
	0\\
	\Omega
\end{array}\right]
\]
\begin{equation}
	\Gamma\left(\Omega\right)=\left[\begin{array}{cc}
		0 & -\Omega^{\top}\\
		\Omega & -\left[\Omega\right]_{\times}
	\end{array}\right]=\left[\begin{array}{cccc}
		0 & -\Omega_{x} & -\Omega_{y} & -\Omega_{z}\\
		\Omega_{x} & 0 & \Omega_{z} & -\Omega_{y}\\
		\Omega_{y} & -\Omega_{z} & 0 & \Omega_{x}\\
		\Omega_{z} & \Omega_{y} & -\Omega_{x} & 0
	\end{array}\right]\label{eq:-3}
\end{equation}
\begin{equation}
	\Xi\left(Q\right)=\left[\begin{array}{c}
		-q^{\top}\\
		q_{0}\mathbf{I}_{3}+\left[q\right]_{\times}
	\end{array}\right]=\left[\begin{array}{ccc}
		-q_{1} & -q_{2} & -q_{3}\\
		q_{0} & -q_{3} & q_{2}\\
		q_{3} & q_{0} & -q_{1}\\
		-q_{2} & q_{1} & q_{0}
	\end{array}\right]\label{eq:-4}
\end{equation}
\begin{equation}
	\Psi\left(Q\right)=\left[\begin{array}{cc}
		0 & -q^{\top}\\
		q & q_{0}\mathbf{I}_{3}+\left[q\right]_{\times}
	\end{array}\right]=\left[\begin{array}{cccc}
		0 & -q_{1} & -q_{2} & -q_{3}\\
		q_{1} & q_{0} & -q_{3} & q_{2}\\
		q_{2} & q_{3} & q_{0} & -q_{1}\\
		q_{3} & -q_{2} & q_{1} & q_{0}
	\end{array}\right]\label{eq:-5}
\end{equation}
\begin{align}
	\bar{\Psi}\left(Q\right) & =\left[\begin{array}{cc}
		0 & q^{\top}\\
		q & q_{0}\mathbf{I}_{3}+\left[q\right]_{\times}
	\end{array}\right]=\left[\begin{array}{cccc}
		0 & q_{1} & q_{2} & q_{3}\\
		q_{1} & q_{0} & -q_{3} & q_{2}\\
		q_{2} & q_{3} & q_{0} & -q_{1}\\
		q_{3} & -q_{2} & q_{1} & q_{0}
	\end{array}\right]\label{eq:-6}
\end{align}
where the relationship between \eqref{eq:-4} and \eqref{eq:-5} is
given in \cite{lefferts1982kalman}.

\subsubsection{Continuous Unit-quaternion Attitude dynamics}

Thus, the attitude dynamics can be defined by\textcolor{red}{{} } 
\begin{align}
	\dot{Q} & =\frac{1}{2}Q\odot\bar{\Omega}=\frac{1}{2}\Gamma\left(\Omega\right)Q=\frac{1}{2}\Xi\left(Q\right)\Omega=\frac{1}{2}\Psi\left(Q\right)\bar{\Omega}\label{eq:OVERVIEW_Q_dot}
\end{align}
According to the expression in \eqref{eq:OVERVIEW_Q_dot}, the attitude
dynamics are linear and time-variant by the following representation
\begin{equation}
	\dot{Q}=\frac{1}{2}\Gamma\left(\Omega\right)Q\label{eq:OVERVIEW_Q_dot_SS1}
\end{equation}

\subsubsection{Discrete Unit-quaternion Attitude dynamics}

The continuous form of the attitude dynamics given in \eqref{eq:OVERVIEW_Q_dot_SS1}
could be defined in discrete form through exact integration by 
\begin{equation}
	Q\left[k+1\right]=\frac{1}{2}\exp\left(\Gamma\left(\Omega\left[k\right]\right)\Delta t\right)Q\left[k\right]\label{eq:OVERVIEW_Q_dot_discrete}
\end{equation}
where $k$ denotes the $k$th sample, $\Delta t$ denotes a time step
which is normally small, and $Q\left[k\right]$ and $\Omega\left[k\right]$
refer to the true unit-quaternion and angular velocity at the $k$th
sample, respectively.

\subsubsection{Sensor Measurements}

Equation \eqref{eq:OVERVIEW_Q_vB} is considered to be the measured
output obtained from sensors attached to the moving body 
\begin{equation}
	\left[\begin{array}{c}
		0\\
		{\rm v}_{i}^{\mathcal{B}}
	\end{array}\right]=Q^{-1}\odot\left[\begin{array}{c}
		0\\
		{\rm v}_{i}^{\mathcal{I}}
	\end{array}\right]\odot Q\label{eq:OVERVIEW_Q_dot_SS2}
\end{equation}
The expression in \eqref{eq:OVERVIEW_Q_dot_SS2} is nonlinear. Extended
Kalman filter \cite{lefferts1982kalman} is the earliest described
filter that took into account both linear dynamics in \eqref{eq:OVERVIEW_Q_dot_SS1}
and the nonlinear dynamics in \eqref{eq:OVERVIEW_Q_dot_SS2}. However,
the nonlinear representation in \eqref{eq:OVERVIEW_Q_dot_SS2} can
be reformulated as follows \cite{choukroun2006novel} 
\begin{align}
	\left[\begin{array}{c}
		0\\
		{\rm v}_{i}^{\mathcal{B}}
	\end{array}\right] & =Q^{-1}\odot\left[\begin{array}{c}
		0\\
		{\rm v}_{i}^{\mathcal{I}}
	\end{array}\right]\odot Q\nonumber \\
	Q\odot\left[\begin{array}{c}
		0\\
		{\rm v}_{i}^{\mathcal{B}}
	\end{array}\right] & =\left[\begin{array}{c}
		0\\
		{\rm v}_{i}^{\mathcal{I}}
	\end{array}\right]\odot Q\nonumber \\
	\left[\begin{array}{cc}
		0 & -\left({\rm v}_{i}^{\mathcal{B}}\right)^{\top}\\
		{\rm v}_{i}^{\mathcal{B}} & -\left[{\rm v}_{i}^{\mathcal{B}}\right]_{\times}
	\end{array}\right]\left[\begin{array}{c}
		q_{0}\\
		q
	\end{array}\right] & =\left[\begin{array}{cc}
		0 & -\left({\rm v}_{i}^{\mathcal{I}}\right)^{\top}\\
		{\rm v}_{i}^{\mathcal{I}} & \left[{\rm v}_{i}^{\mathcal{I}}\right]_{\times}
	\end{array}\right]\left[\begin{array}{c}
		q_{0}\\
		q
	\end{array}\right]\label{eq:OVERVIEW_Q_dot_SS3}
\end{align}
Thus, the expression in \eqref{eq:OVERVIEW_Q_dot_SS3} is equivalent
to 
\begin{equation}
	0=\left[\begin{array}{cc}
		0 & -\left({\rm v}_{i}^{\mathcal{B}}-{\rm v}_{i}^{\mathcal{I}}\right)^{\top}\\
		{\rm v}_{i}^{\mathcal{B}}-{\rm v}_{i}^{\mathcal{I}} & -\left[{\rm v}_{i}^{\mathcal{B}}+{\rm v}_{i}^{\mathcal{I}}\right]_{\times}
	\end{array}\right]\left[\begin{array}{c}
		q_{0}\\
		q
	\end{array}\right]\label{eq:OVERVIEW_Q_dot_SS4}
\end{equation}
Therefore, the linear and time-variant state-space representation
of the attitude problem in terms of quaternion is equivalent to 
\begin{align}
	\dot{Q} & =\frac{1}{2}\left[\begin{array}{cc}
		0 & -\Omega^{\top}\\
		\Omega & -\left[\Omega\right]_{\times}
	\end{array}\right]\left[\begin{array}{c}
		q_{0}\\
		q
	\end{array}\right]\label{eq:OVERVIEW_Q_dot_SSX}\\
	& =\frac{1}{2}\Gamma\left(\Omega\right)Q\nonumber 
\end{align}
\begin{align}
	Y & =0=\left[\begin{array}{cc}
		0 & -\left({\rm v}_{i}^{\mathcal{B}}-{\rm v}_{i}^{\mathcal{I}}\right)^{\top}\\
		{\rm v}_{i}^{\mathcal{B}}-{\rm v}_{i}^{\mathcal{I}} & -\left[{\rm v}_{i}^{\mathcal{B}}+{\rm v}_{i}^{\mathcal{I}}\right]_{\times}
	\end{array}\right]\left[\begin{array}{c}
		q_{0}\\
		q
	\end{array}\right]\label{eq:OVERVIEW_Q_dot_SSY}\\
	& =0=H\left({\rm v}_{i}^{\mathcal{I}},{\rm v}_{i}^{\mathcal{B}}\right)Q\nonumber 
\end{align}
where $Q\in\mathbb{S}^{3}$ is the state vector, and $Y\in\mathbb{R}^{4}$
is the output vector. Both $\Gamma\left(\Omega\right)$ and $H\left({\rm v}_{i}^{\mathcal{I}},{\rm v}_{i}^{\mathcal{B}}\right)$
are time-variant known matrices obtained from the measurements of
the sensors attached to the moving body. The representation in \eqref{eq:OVERVIEW_Q_dot_SSY}
is valid only for ${\rm v}_{i}^{\mathcal{B}}$ free of noise and bias
components. This is the case because equation \eqref{eq:OVERVIEW_Q_dot_SSY}
disregards noise and bias attached to ${\rm v}_{i}^{\mathcal{B}}$.
Nonetheless, this approach produced impressive results for an uncertain
${\rm v}_{i}^{\mathcal{B}}$ as described in \cite{choukroun2006novel}.

\subsection{Rotational acceleration}

The relationship between rotational acceleration $\dot{\bar{\Omega}}$
and the quaternion derivative can be defined as follows: 
\begin{align*}
	Q^{*}\odot Q & =Q_{{\rm I}}\\
	\dot{Q}^{*}\odot Q+Q^{*}\odot\dot{Q} & =0\\
	2\dot{Q}^{*}\odot Q+2Q^{*}\odot\dot{Q} & =0\\
	2\dot{Q}^{*}\odot Q+\bar{\Omega} & =0\\
	\bar{\Omega} & =-2\dot{Q}^{*}\odot Q
\end{align*}
Now one can find 
\begin{align*}
	\ddot{Q} & =\frac{1}{2}\left(\dot{Q}\odot\bar{\Omega}+Q\odot\dot{\bar{\Omega}}\right)\\
	& =\frac{1}{2}\left(-2\dot{Q}\odot\dot{Q}^{*}\odot Q+Q\odot\dot{\bar{\Omega}}\right)\\
	Q^{*}\odot\ddot{Q} & =\frac{1}{2}\left(-2Q^{*}\odot\dot{Q}\odot\dot{Q}^{*}\odot Q+\dot{\bar{\Omega}}\odot Q_{{\rm I}}\right)\\
	& =\frac{1}{2}\left(-2Q^{*}\odot\dot{Q}\odot\dot{Q}^{*}\odot Q+\dot{\bar{\Omega}}\right)
\end{align*}
\begin{align*}
	\dot{\bar{\Omega}} & =2\left(Q^{*}\odot\ddot{Q}+Q^{*}\odot\dot{Q}\odot\dot{Q}^{*}\odot Q\right)
\end{align*}
as $Q^{*}\odot\dot{Q}\odot\dot{Q}^{*}\odot Q=\dot{Q}\odot\dot{Q}^{*}$
\begin{align}
	\dot{\bar{\Omega}} & =2\left(Q^{*}\odot\ddot{Q}+\dot{Q}\odot\dot{Q}^{*}\right)\label{eq:OVERVIEW_Q_accel}
\end{align}
Note that $\left\Vert \dot{Q}\right\Vert $ is not necessarily equal
to $1$.

\subsection{Unit-quaternion update}

The rotational matrix can be constructed knowing unit-quaternion if
and only if $\left\Vert Q\right\Vert =1$. During the control/estimation
process, the unit-quaternion may lose precision and thereby, $\left\Vert Q\right\Vert \neq1$.
In order to achieve valid mapping from $Q$ to $\mathbb{SO}\left(3\right)$,
it is necessary to maintain $\left\Vert Q\right\Vert =1$ at each
time instant which can be accomplished through the substitution of
$Q$ by 
\begin{equation}
	Q=\frac{Q}{\left\Vert Q\right\Vert }\label{eq:-17}
\end{equation}

\subsection{Unit-quaternion error and error derivative}

Let $Q\in\mathbb{S}^{3}$ be the true unit-quaternion vector and let
the desired/estimator unit-quaternion vector be given by $Q_{\star}=\left[q_{\star0},q_{\star}^{\top}\right]^{\top}\in\mathbb{S}^{3}$
for all $q_{\star0}\in\mathbb{R}$ and $q_{\star}\in\mathbb{R}^{3}$.
The true and desired/estimator unit-quaternion dynamics are given
by 
\begin{align}
	\dot{Q} & =\frac{1}{2}\Psi\left(Q\right)\bar{\Omega}=\frac{1}{2}\left[\begin{array}{cc}
		0 & -q^{\top}\\
		q & q_{0}\mathbf{I}_{3}+\left[q\right]_{\times}
	\end{array}\right]\bar{\Omega}=\frac{1}{2}\left[\begin{array}{c}
		-q^{\top}\Omega\\
		\left(q_{0}\mathbf{I}_{3}+\left[q\right]_{\times}\right)\Omega
	\end{array}\right]\label{eq:-24}\\
	\nonumber \\
	\dot{Q}_{\star} & =\frac{1}{2}\Psi\left(Q_{\star}\right)\bar{\Omega}_{\star}=\frac{1}{2}\left[\begin{array}{cc}
		0 & -q_{\star}^{\top}\\
		q_{\star} & q_{\star0}\mathbf{I}_{3}+\left[q_{\star}\right]_{\times}
	\end{array}\right]\bar{\Omega}_{\star}=\frac{1}{2}\left[\begin{array}{c}
		-q_{\star}^{\top}\Omega_{\star}\\
		\left(q_{\star0}\mathbf{I}_{3}+\left[q_{\star}\right]_{\times}\right)\Omega_{\star}
	\end{array}\right]\label{eq:-25}
\end{align}
where $\bar{\Omega}_{\star}=\left[0,\Omega_{\star}^{\top}\right]^{\top}\in\mathbb{R}^{4}$
and $\Omega_{\star}\in\mathbb{R}^{3}$. From \eqref{eq:-25}, one
could find 
\begin{equation}
	\dot{Q}_{\star}^{-1}=\frac{1}{2}\bar{\Psi}\left(Q_{\star}\right)\bar{\Omega}_{\star}=-\frac{1}{2}\left[\begin{array}{cc}
		0 & q_{\star}^{\top}\\
		q_{\star} & q_{\star0}\mathbf{I}_{3}+\left[q_{\star}\right]_{\times}
	\end{array}\right]\bar{\Omega}_{\star}=-\frac{1}{2}\left[\begin{array}{c}
		q_{\star}^{\top}\Omega_{\star}\\
		\left(q_{\star0}\mathbf{I}_{3}+\left[q_{\star}\right]_{\times}\right)\Omega_{\star}
	\end{array}\right]\label{eq:-26}
\end{equation}
Let the error between the desired and the true unit-quaternion be
defined by 
\begin{align}
	\tilde{Q} & =\left[\begin{array}{c}
		\tilde{q}_{0}\\
		\tilde{q}
	\end{array}\right]=Q_{\star}^{-1}\odot Q\nonumber \\
	& =\left[\begin{array}{c}
		q_{\star0}\\
		-q_{\star}
	\end{array}\right]\odot\left[\begin{array}{c}
		q_{0}\\
		q
	\end{array}\right]\nonumber \\
	& =\left[\begin{array}{c}
		q_{\star0}q_{0}+q_{\star}^{\top}q\\
		q_{\star0}q-q_{0}q_{\star}-\left[q_{\star}\right]_{\times}q
	\end{array}\right]\label{eq:OVERVIEW_Q_error}
\end{align}
where $\tilde{Q}=\left[\tilde{q}_{0},\tilde{q}^{\top}\right]^{\top}\in\mathbb{S}^{3}$
for all $\tilde{q}_{0}\in\mathbb{R}$ and $\tilde{q}\in\mathbb{R}^{3}$.
The objective of attitude filter/control in terms of unit-quaternion
is to drive the error in unit-quaternion $\tilde{Q}\rightarrow Q_{{\rm I}}=\left[1,0,0,0\right]^{\top}$.
Driving $\tilde{Q}\rightarrow Q_{{\rm I}}$ implies $\tilde{R}\rightarrow\mathbf{I}_{3}$,
since we have 
\begin{align}
	\mathcal{R}_{Q}\left(\tilde{Q}\right) & =\left(\tilde{q}_{0}^{2}-\tilde{q}^{\top}\tilde{q}\right)\mathbf{I}_{3}+2\tilde{q}\tilde{q}^{\top}+2\tilde{q}_{0}\left[\tilde{q}\right]_{\times}\nonumber \\
	\tilde{Q}=Q_{{\rm I}} & \Leftrightarrow\mathcal{R}_{Q}\left(\tilde{Q}\right)=\mathbf{I}_{3}\label{eq:-34}
\end{align}
Accordingly, the unit-quaternion error dynamics are given by 
\begin{equation}
	\dot{\tilde{Q}}=\dot{Q}_{\star}^{-1}\odot Q+Q_{\star}^{-1}\odot\dot{Q}\label{eq:OVERVIEW_Q_error_dot}
\end{equation}
Consequently, from \eqref{eq:-24} and \eqref{eq:-26} one has 
\begin{align}
	\dot{Q}_{\star}^{-1}\odot Q & =\frac{1}{2}\left[\begin{array}{c}
		-q_{\star}^{\top}\Omega_{\star}\\
		-q_{\star0}\Omega_{\star}-\left[q_{\star}\right]_{\times}\Omega_{\star}
	\end{array}\right]\odot\left[\begin{array}{c}
		q_{0}\\
		q
	\end{array}\right]\nonumber \\
	& =\frac{1}{2}\left[\begin{array}{c}
		-q_{0}q_{\star}^{\top}\Omega_{\star}+\left(q_{\star0}\Omega_{\star}^{\top}+\Omega_{\star}^{\top}\left[q_{\star}\right]_{\times}^{\top}\right)q\\
		q_{0}\left(-q_{\star0}\Omega_{\star}-\left[q_{\star}\right]_{\times}\Omega_{\star}\right)+\left(-q_{\star}^{\top}\Omega_{\star}\right)q+\left[-q_{\star0}\Omega_{\star}-\left[q_{\star}\right]_{\times}\Omega_{\star}\right]_{\times}q
	\end{array}\right]\nonumber \\
	& =\frac{1}{2}\left[\begin{array}{c}
		\left(q_{\star0}q-q_{0}q_{\star}-\left[q_{\star}\right]_{\times}q\right)^{\top}\Omega_{\star}\\
		-\left(q_{0}q_{\star0}\mathbf{I}_{3}+qq_{\star}^{\top}\right)\Omega_{\star}+\left(q_{\star0}\left[q\right]_{\times}-q_{0}\left[q_{\star}\right]_{\times}-\left[q_{\star}\right]_{\times}\left[q\right]_{\times}\right)\Omega_{\star}
	\end{array}\right]\nonumber \\
	& =\frac{1}{2}\left[\begin{array}{c}
		\tilde{q}^{\top}\Omega_{\star}\\
		-\tilde{q}_{0}\Omega_{\star}+\left[\tilde{q}\right]_{\times}\Omega_{\star}
	\end{array}\right]\label{eq:OVERVIEW_Q_error_dot1}
\end{align}
\begin{align}
	Q_{\star}^{-1}\odot\dot{Q} & =\frac{1}{2}\left[\begin{array}{c}
		q_{\star0}\left(-q^{\top}\Omega\right)-\left(-q_{\star}\right)^{\top}\left(q_{0}\Omega+\left[q\right]_{\times}\Omega\right)\\
		q_{\star0}\left(q_{0}\Omega+\left[q\right]_{\times}\Omega\right)+\left(-q^{\top}\Omega\right)\left(-q_{\star}\right)+\left[-q_{\star}\right]_{\times}\left(q_{0}\Omega+\left[q\right]_{\times}\Omega\right)
	\end{array}\right]\nonumber \\
	& =\frac{1}{2}\left[\begin{array}{c}
		\left(-q_{\star0}q^{\top}+q_{0}q_{\star}^{\top}+q^{\top}\left[q_{\star}\right]_{\times}^{\top}\right)\Omega\\
		\left(q_{\star0}q_{0}\mathbf{I}_{3}+q_{\star}q^{\top}\right)\Omega+\left[q_{\star0}q-q_{0}q_{\star}-\left[q_{\star}\right]_{\times}q\right]_{\times}\Omega
	\end{array}\right]\nonumber \\
	& =\frac{1}{2}\left[\begin{array}{c}
		-\tilde{q}^{\top}\Omega\\
		\tilde{q}_{0}\Omega+\left[\tilde{q}\right]_{\times}\Omega
	\end{array}\right]\label{eq:OVERVIEW_Q_error_dot2}
\end{align}
From \eqref{eq:OVERVIEW_Q_error_dot1} and \eqref{eq:OVERVIEW_Q_error_dot2}
the unit-quaternion error dynamics in \eqref{eq:OVERVIEW_Q_error_dot}
are equivalent to 
\begin{align}
	\dot{\tilde{Q}} & =\dot{Q}_{\star}^{-1}\odot Q+Q_{\star}^{-1}\odot\dot{Q}\nonumber \\
	& =\frac{1}{2}\left[\begin{array}{c}
		\tilde{q}^{\top}\\
		-\tilde{q}_{0}\mathbf{I}_{3}+\left[\tilde{q}\right]_{\times}
	\end{array}\right]\Omega_{\star}+\frac{1}{2}\left[\begin{array}{c}
		-\tilde{q}^{\top}\\
		\tilde{q}_{0}\mathbf{I}_{3}+\left[\tilde{q}\right]_{\times}
	\end{array}\right]\Omega\nonumber \\
	& =\frac{1}{2}\left[\begin{array}{c}
		\tilde{q}^{\top}\left(\Omega_{\star}-\Omega\right)\\
		\tilde{q}_{0}\left(\Omega-\Omega_{\star}\right)+\left[\tilde{q}\right]_{\times}\left(\Omega_{\star}+\Omega\right)
	\end{array}\right]\label{eq:OVERVIEW_Q_error_dot3}
\end{align}
Therefore, the problem can be summarized as follows 
\begin{align*}
	\tilde{Q} & =\left[\begin{array}{c}
		\tilde{q}_{0}\\
		\tilde{q}
	\end{array}\right]\\
	& =\left[\begin{array}{c}
		q_{\star0}q_{0}+q_{\star}^{\top}q\\
		q_{\star0}q-q_{0}q_{\star}-\left[q_{\star}\right]_{\times}q
	\end{array}\right]\\
	\dot{\tilde{Q}} & =\frac{1}{2}\left[\begin{array}{c}
		\tilde{q}^{\top}\left(\Omega_{\star}-\Omega\right)\\
		\tilde{q}_{0}\left(\Omega-\Omega_{\star}\right)+\left[\tilde{q}\right]_{\times}\left(\Omega_{\star}+\Omega\right)
	\end{array}\right]
\end{align*}
The above-mentioned expression can be simplified if the error in angular
velocity is selected using $\tilde{\Omega}=\Omega-\mathcal{R}_{Q}\left(\tilde{Q}\right)\Omega_{\star}$.
Therefore, one can find 
\begin{align*}
	\tilde{q}^{\top}\left(\Omega_{\star}-\Omega\right) & =\tilde{q}^{\top}\Omega_{\star}-\tilde{q}^{\top}\mathcal{R}_{Q}\left(\tilde{Q}\right)\Omega_{\star}+\tilde{q}^{\top}\mathcal{R}_{Q}\left(\tilde{Q}\right)\Omega_{\star}-\tilde{q}^{\top}\Omega\\
	& =\tilde{q}^{\top}\left(\mathcal{R}_{Q}\left(\tilde{Q}\right)\Omega_{\star}-\Omega\right)+\tilde{q}^{\top}\left(\mathbf{I}_{3}-\mathcal{R}_{Q}\left(\tilde{Q}\right)\right)\Omega_{\star}\\
	& =\tilde{q}^{\top}\left(\mathcal{R}_{Q}\left(\tilde{Q}\right)\Omega_{\star}-\Omega\right)+\tilde{q}^{\top}\underbrace{\left(\mathbf{I}_{3}-\mathbf{I}_{3}-2\tilde{q}_{0}\left[\tilde{q}\right]_{\times}-2\left[\tilde{q}\right]_{\times}^{2}\right)}_{=0}\Omega_{\star}\\
	& =-\tilde{q}^{\top}\left(\Omega-\mathcal{R}_{Q}\left(\tilde{Q}\right)\Omega_{\star}\right)
\end{align*}

\noindent 
\begin{align*}
	\left[\tilde{q}\right]_{\times}\left(\Omega_{\star}+\Omega\right)+\tilde{q}_{0}\left(\Omega-\Omega_{\star}\right) & =\left[\tilde{q}\right]_{\times}\Omega_{\star}+\left[\tilde{q}\right]_{\times}\Omega+\left[\tilde{q}\right]_{\times}\mathcal{R}_{Q}\left(\tilde{Q}\right)\Omega_{\star}-\left[\tilde{q}\right]_{\times}\mathcal{R}_{Q}\left(\tilde{Q}\right)\Omega_{\star}\\
	& \hspace{1em}+\tilde{q}_{0}\Omega-\tilde{q}_{0}\Omega_{\star}+\tilde{q}_{0}\mathcal{R}_{Q}\left(\tilde{Q}\right)\Omega_{\star}-\tilde{q}_{0}\mathcal{R}_{Q}\left(\tilde{Q}\right)\Omega_{\star}\\
	& =\left[\tilde{q}\right]_{\times}\left(\Omega-\mathcal{R}_{Q}\left(\tilde{Q}\right)\Omega_{\star}\right)+\tilde{q}_{0}\left(\Omega-\mathcal{R}_{Q}\left(\tilde{Q}\right)\Omega_{\star}\right)
\end{align*}
hence 
\begin{align}
	\dot{\tilde{Q}} & =\left[\begin{array}{c}
		\dot{\tilde{q}}_{0}\\
		\dot{\tilde{q}}
	\end{array}\right]=\frac{1}{2}\left[\begin{array}{c}
		-\tilde{q}^{\top}\left(\Omega-\mathcal{R}_{Q}\left(\tilde{Q}\right)\Omega_{\star}\right)\\
		\left(\left[\tilde{q}\right]_{\times}+\tilde{q}_{0}\mathbf{I}_{3}\right)\left(\Omega-\mathcal{R}_{Q}\left(\tilde{Q}\right)\Omega_{\star}\right)
	\end{array}\right]\nonumber \\
	& =\frac{1}{2}\left[\begin{array}{c}
		-\tilde{q}^{\top}\\
		\left[\tilde{q}\right]_{\times}+\tilde{q}_{0}\mathbf{I}_{3}
	\end{array}\right]\tilde{\Omega}\label{eq:-2}
\end{align}
which is equivalent to 
\begin{align}
	\dot{\tilde{Q}} & =\frac{1}{2}\Psi\left(\tilde{Q}\right)\bar{\tilde{\Omega}}\label{eq:OVERVIEW_Q_error_dot_Final}
\end{align}
where $\bar{\tilde{\Omega}}=\left[0,\tilde{\Omega}^{\top}\right]^{\top}$.

\subsection{Problem of unit-quaternion }

Despite providing a global representation of the attitude and being
free of non-singularity in the attitude parameterization, unit-quaternion
vector is non-unique. Two different unit-quaternion vectors can result
in the same rotational matrix such that 
\begin{align}
	\mathcal{R}_{Q}\left(Q\right) & =\mathcal{R}_{Q}\left(-Q\right)\hspace{1em}\forall Q\in\mathbb{S}^{3}\label{eq:-22}
\end{align}
such as 
\begin{align}
	\mathcal{R}_{Q}\left(S\right) & =\mathcal{R}_{Q}\left(P\right),\hspace{1em}\forall S=\left[q_{0},q_{1},q_{2},q_{3}\right]^{\top},P=\left[-q_{0},-q_{1},-q_{2},-q_{3}\right]^{\top}\label{eq:-23}
\end{align}
For example, consider the following two unit-quaternion vectors 
\begin{align*}
	S & =\left[\begin{array}{c}
		0.7794\\
		-0.1440\\
		0.4623\\
		-0.3976
	\end{array}\right]\in\mathbb{S}^{3},\hspace{1em}P=-S=\left[\begin{array}{c}
		-0.7794\\
		0.1440\\
		-0.4623\\
		0.3976
	\end{array}\right]\in\mathbb{S}^{3}
\end{align*}
One can easily find $\mathcal{R}_{Q}:\mathbb{S}^{3}\rightarrow\mathbb{SO}\left(3\right)$
\[
\mathcal{R}_{Q}\left(S\right)=\mathcal{R}_{Q}\left(P\right)=\left[\begin{array}{ccc}
	0.2563 & 0.4867 & 0.8351\\
	-0.7529 & 0.6423 & -0.1433\\
	-0.6061 & -0.5921 & 0.5311
\end{array}\right]\in\mathbb{SO}\left(3\right)
\]
Hence, the controller should be carefully designed when using unit-quaternion
for attitude parameterization \cite{bhat2000topological}. In brief,
physical attitude $R\in\mathbb{SO}\left(3\right)$, which has a unique
orientation, is represented by a pair of antipodal quaternions $\pm Q\in\mathbb{S}^{3}$,
such that $R=\mathcal{R}_{Q}\left(\pm Q\right)\in\mathbb{SO}\left(3\right)$.

\noindent %
\noindent\makebox[1\linewidth]{%
	\rule{0.6\textwidth}{1.4pt}%
}

\section{Simulation \label{sec:Simulation}}

In this section, two examples are presented to illustrate that not
any angular velocity ($\Omega$) can be obtained by the means of the
time derivatives of Euler angles ($\dot{\xi}$). It also shows that
Rodriguez vector has a unique representations, in spite of the fact,
that it cannot achieve some configurations. It should be noted that
the true attitude dynamics follow \cite{shuster1993survey} 
\[
\dot{R}=R\left[\Omega\right]_{\times}
\]
The attitude $R$ obtained from the above mentioned dynamics is the
true attitude. Thus, any method of attitude parameterization, such
as Euler angles, Rodriguez vector, or unit-quaternion obtained from
the dynamics in \eqref{eq:OVERVIEW_ROD_dot}, \eqref{eq:OVERVIEW_ROD_dot}
and \eqref{eq:OVERVIEW_Q_dot}, respectively, should be suitable for
acquiring the true representation of the true attitude $R$.

\subsection*{Example 1:}

Consider this angular velocity 
\begin{equation}
	\Omega=\left[\begin{array}{c}
		0.1{\rm sin}\left(0.3376t\right)\\
		0.07{\rm sin}\left(0.6079t+\pi\right)\\
		0.05{\rm sin}\left(0.7413t+\frac{\pi}{3}\right)
	\end{array}\right]\left({\rm rad/sec}\right)\label{eq:sim_ang}
\end{equation}
Let the initial attitude be given by 
\begin{equation}
	R\left(0\right)=\left[\begin{array}{ccc}
		0.9479 & -0.2040 & 0.2448\\
		0.2177 & 0.9756 & -0.0297\\
		-0.2328 & 0.0814 & 0.9691
	\end{array}\right]\in\mathbb{SO}\left(3\right)\label{eq:R0}
\end{equation}
The initial condition of $R\left(0\right)$ in terms of Euler angles,
Rodriguez vector and unit-quaternion are listed in Table \ref{tab:Initial-conditions1}.

\begin{table}[H]
	\centering{}\caption{\label{tab:Initial-conditions1}Initial conditions - Example 1}
	\begin{tabular}{ccc}
		\toprule 
		\addlinespace
		Representation  & Mapping  & Numerical values\tabularnewline\addlinespace
		\midrule
		\midrule 
		\addlinespace
		Euler angles  & $\begin{array}{c}
			\xi:\mathbb{SO}\left(3\right)\rightarrow\mathbb{R}^{3}\\
			R\left(0\right)\rightarrow\xi\left(0\right)
		\end{array}$  & $\xi\left(0\right)=\left[\begin{array}{c}
			\phi\\
			\theta\\
			\psi
		\end{array}\right]=\left[\begin{array}{c}
			4.8035\\
			13.4601\\
			12.9329
		\end{array}\right]\left({\rm deg}\right)$\tabularnewline\addlinespace
		\midrule 
		\addlinespace
		Rodriguez vector  & $\begin{array}{c}
			\rho:\mathbb{SO}\left(3\right)\rightarrow\mathbb{R}^{3}\\
			R\left(0\right)\rightarrow\rho\left(0\right)
		\end{array}$  & $\rho\left(0\right)=\left[\begin{array}{c}
			\rho_{1}\\
			\rho_{2}\\
			\rho_{3}
		\end{array}\right]=\left[\begin{array}{c}
			0.0286\\
			0.1227\\
			0.1083
		\end{array}\right]$\tabularnewline\addlinespace
		\midrule 
		\addlinespace
		Unit-quaternion  & $\begin{array}{c}
			Q:\mathbb{SO}\left(3\right)\rightarrow\mathbb{S}^{3}\\
			R\left(0\right)\rightarrow Q\left(0\right)
		\end{array}$  & $Q\left(0\right)=\left[\begin{array}{c}
			q_{0}\\
			q
		\end{array}\right]=\left[\begin{array}{c}
			0.9865\\
			0.0282\\
			0.1210\\
			0.1069
		\end{array}\right]$\tabularnewline\addlinespace
		\bottomrule
	\end{tabular}
\end{table}

\begin{figure}[h!]
	\centering{}\includegraphics[scale=0.43]{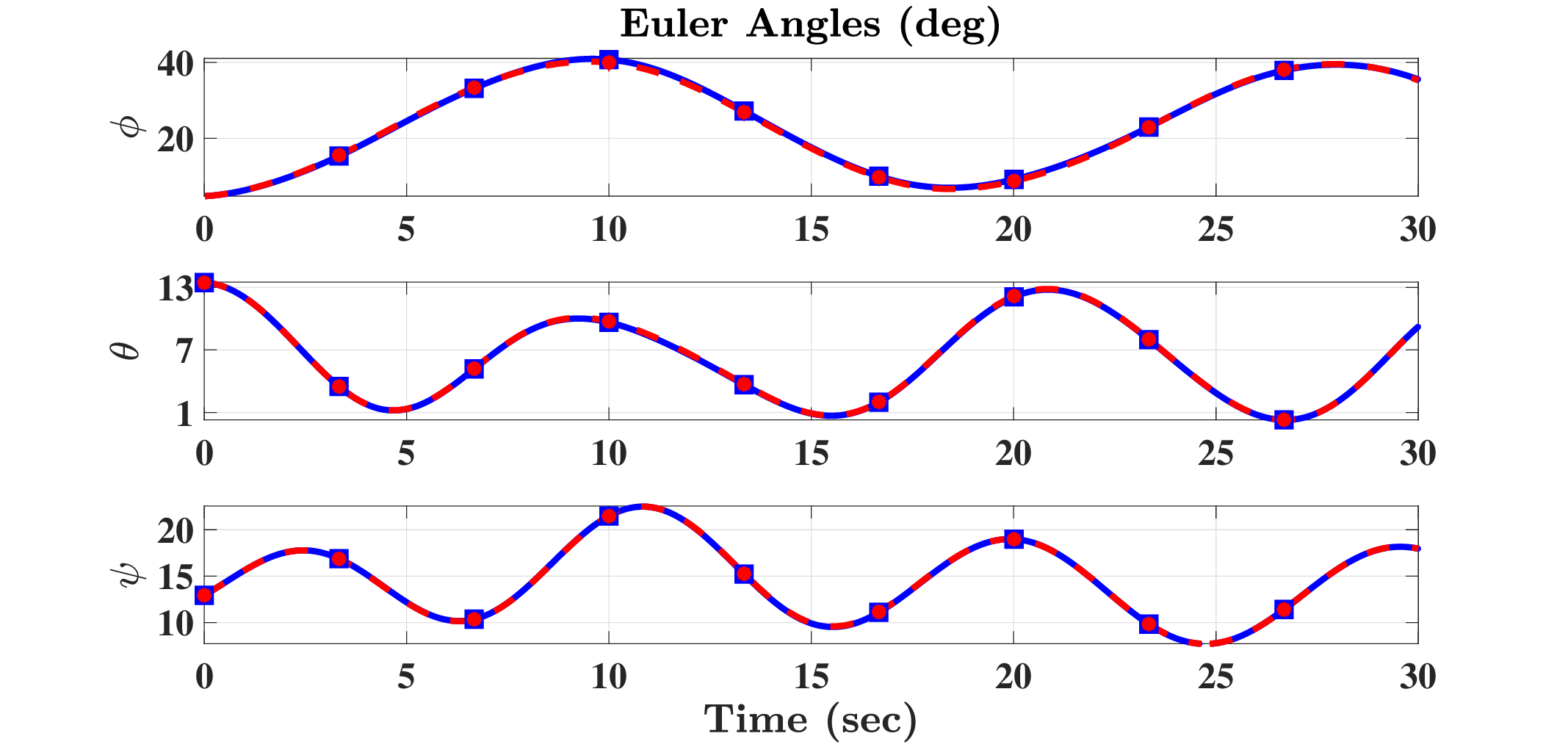}\caption{Euler Angles (\textcolor{blue}{$\dot{R}=R\left[\Omega\right]_{\times},\hspace{1em}R\rightarrow\xi$})
		vs (\textcolor{red}{$\dot{\xi}=\mathcal{J}\Omega$}) - Example 1}
	\label{fig:OVERVIEW_SO3_SIMU2} 
\end{figure}

\begin{figure}[h!]
	\centering{}\includegraphics[scale=0.43]{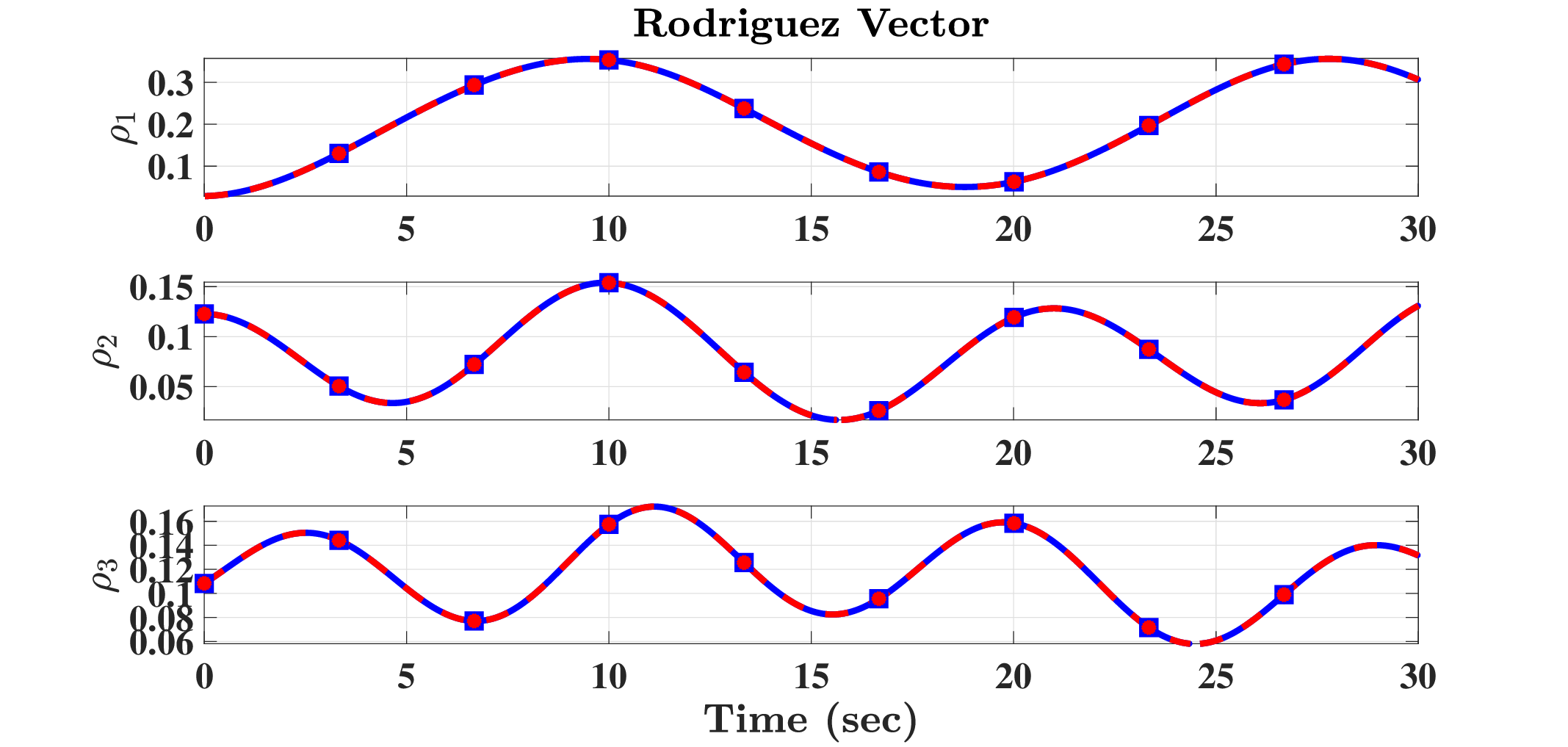}\caption{Rodriguez vector (\textcolor{blue}{$\dot{R}=R\left[\Omega\right]_{\times},\hspace{1em}R\rightarrow\rho$})
		vs (\textcolor{red}{$\dot{\rho}=\frac{1}{2}\left(\mathbf{I}_{3}+\left[\rho\right]_{\times}+\rho\rho^{\top}\right)\Omega$})
		- Example 1}
	\label{fig:OVERVIEW_SO3_SIMU3} 
\end{figure}

\begin{figure}[h!]
	\centering{}\includegraphics[scale=0.43]{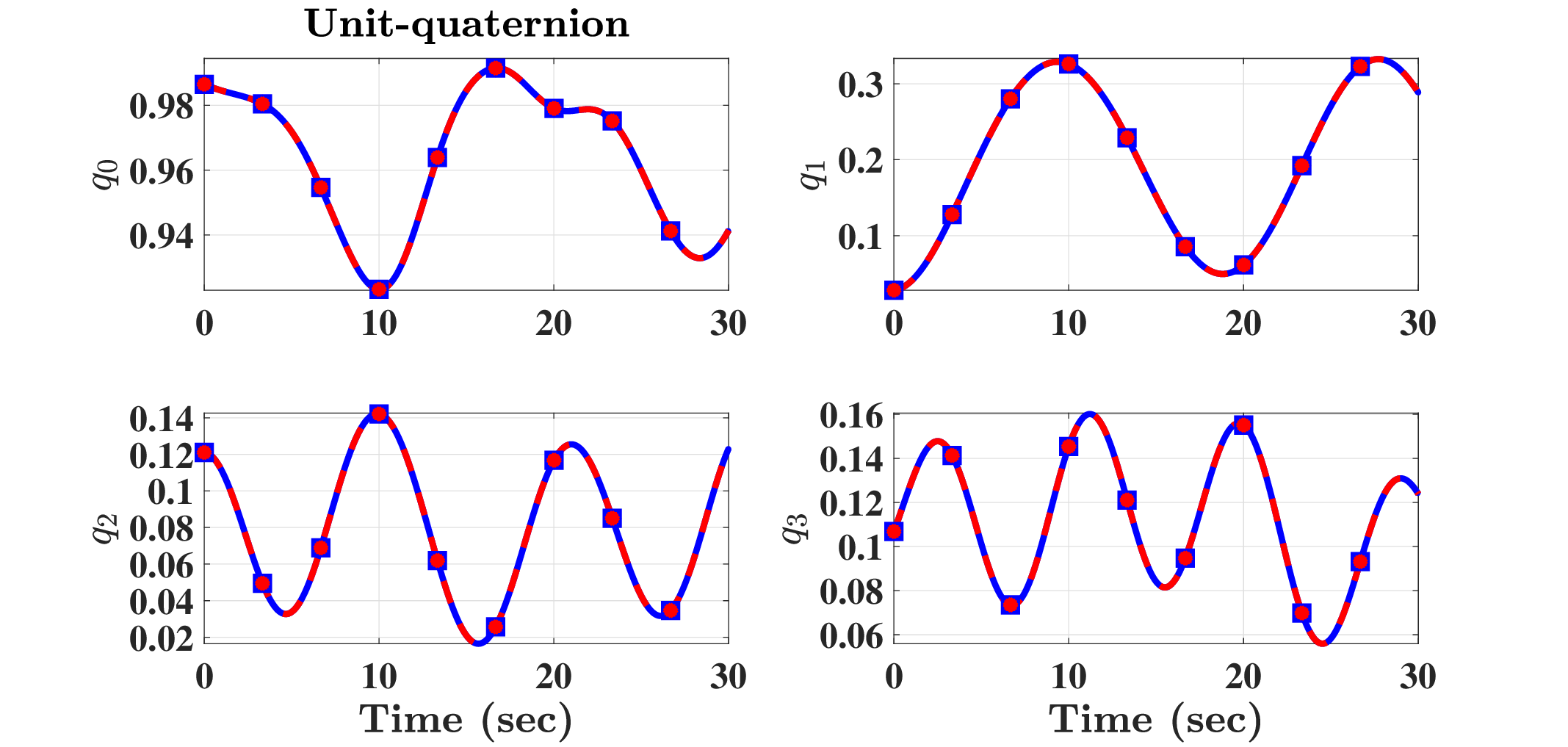}\caption{Unit-Quaternion (\textcolor{blue}{$\dot{R}=R\left[\Omega\right]_{\times},\hspace{1em}R\rightarrow Q$})
		vs (\textcolor{red}{$\dot{Q}=\frac{1}{2}\Gamma\left(\Omega\right)Q$})
		- Example 1}
	\label{fig:OVERVIEW_SO3_SIMU4} 
\end{figure}

The attitude $R$ has been obtained from the dynamics in \eqref{eq:OVERVIEW_SO3_dyn}
($\dot{R}=R\left[\Omega\right]_{\times}$) given $R\left(0\right)$
in \eqref{eq:R0} and the angular velocity in \eqref{eq:sim_ang}.
The mapping of the attitude $R$ from $\mathbb{SO}\left(3\right)$
to Euler angles, Rodriguez vector, and unit-quaternion is obtained
from \eqref{eq:SO3_EUL}, \eqref{eq:OVERVIEW_SO3_ROD}, and \eqref{eq:OVERVIEW_Q_comp_1},
respectively, and depicted in Figure \ref{fig:OVERVIEW_SO3_SIMU2},
\ref{fig:OVERVIEW_SO3_SIMU3}, and \ref{fig:OVERVIEW_SO3_SIMU4},
in blue colors. Euler angles obtained from Euler dynamics in \eqref{eq:OVERVIEW_EUL_dot}
are plotted in red color in Figure \ref{fig:OVERVIEW_SO3_SIMU2} against
blue-colored Euler angles obtained through the mapping in \eqref{eq:SO3_EUL}
from $\mathbb{SO}\left(3\right)$. Rodriguez vector obtained from
Rodriguez vector dynamics in \eqref{eq:OVERVIEW_ROD_dot} is plotted
in red color in Figure \ref{fig:OVERVIEW_SO3_SIMU3} against Rodriguez
vector obtained through the mapping in \eqref{eq:OVERVIEW_SO3_ROD}
from $\mathbb{SO}\left(3\right)$ in blue colors. Unit-quaternion
obtained from unit-quaternion dynamics in \eqref{eq:OVERVIEW_Q_dot}
is plotted in red color in Figure \ref{fig:OVERVIEW_SO3_SIMU3} against
unit-quaternion obtained through the mapping in \eqref{eq:OVERVIEW_Q_comp_1},
from $\mathbb{SO}\left(3\right)$ in blue colors. Table \ref{tab:colr_dynamics1}
summarizes the comparison between the three representations with colors
corresponding to those used in Figure \ref{fig:OVERVIEW_SO3_SIMU2},
\ref{fig:OVERVIEW_SO3_SIMU3}, and \ref{fig:OVERVIEW_SO3_SIMU4}.
In Example 1, low gain and rate of angular velocity is considered.
Figure \ref{fig:OVERVIEW_SO3_SIMU2}, \ref{fig:OVERVIEW_SO3_SIMU3},
and \ref{fig:OVERVIEW_SO3_SIMU4} show \textbf{accurate tracking}
between the mapping from $\mathbb{SO}\left(3\right)$ and the dynamics
obtained from \eqref{eq:SO3_EUL}, \eqref{eq:OVERVIEW_SO3_ROD}, and
\eqref{eq:OVERVIEW_Q_comp_1}, respectively. 
\begin{table}[H]
	\centering{}\caption{\label{tab:colr_dynamics1}Representation, color notation and related
		mapping - Example 1}
	\begin{tabular}{ccccc}
		\toprule 
		\addlinespace
		Representation  & Dynamics  & Color notation  & Mapping  & Figure \tabularnewline\addlinespace
		\midrule
		\midrule 
		\addlinespace
		\textcolor{blue}{$\mathbb{SO}\left(3\right)$} & \textcolor{blue}{$\dot{R}=R\left[\Omega\right]_{\times}$} & \textcolor{blue}{blue} & \textcolor{blue}{$R\left(t\right)\rightarrow\xi\left(t\right)$}  & \textcolor{blue}{5}\tabularnewline\addlinespace
		\midrule 
		\addlinespace
		&  &  & \textcolor{blue}{$R\left(t\right)\rightarrow\rho\left(t\right)$}  & \textcolor{blue}{6}\tabularnewline\addlinespace
		\midrule 
		\addlinespace
		&  &  & \textcolor{blue}{$R\left(t\right)\rightarrow Q\left(t\right)$}  & \textcolor{blue}{7}\tabularnewline\addlinespace
		\midrule 
		\addlinespace
		\textcolor{red}{Euler angles}  & \textcolor{red}{$\dot{\xi}=\mathcal{J}\Omega$}  & \textcolor{red}{Red}  & \textcolor{red}{$\int\dot{\xi}=\xi\left(t\right)$}  & \textcolor{red}{\ref{fig:OVERVIEW_SO3_SIMU2}}\tabularnewline\addlinespace
		\midrule 
		\addlinespace
		\textcolor{red}{Rodriguez vector}  & \textcolor{red}{$\dot{\rho}=\frac{1}{2}\left(\mathbf{I}_{3}+\left[\rho\right]_{\times}+\rho\rho^{\top}\right)\Omega$}  & \textcolor{red}{Red}  & \textcolor{red}{$\int\dot{\rho}=\rho\left(t\right)$}  & \textcolor{red}{\ref{fig:OVERVIEW_SO3_SIMU3}}\tabularnewline\addlinespace
		\midrule 
		\addlinespace
		\textcolor{red}{Unit-quaternion}  & \textcolor{red}{$\dot{Q}=\frac{1}{2}\Gamma\left(\Omega\right)Q$}  & \textcolor{red}{Red}  & \textcolor{red}{$\int\dot{Q}=Q\left(t\right)$}  & \textcolor{red}{\ref{fig:OVERVIEW_SO3_SIMU4} }\tabularnewline\addlinespace
		\bottomrule
	\end{tabular}
\end{table}

\newpage{}

\begin{figure}[h!]
	\centering{}\includegraphics[scale=0.43]{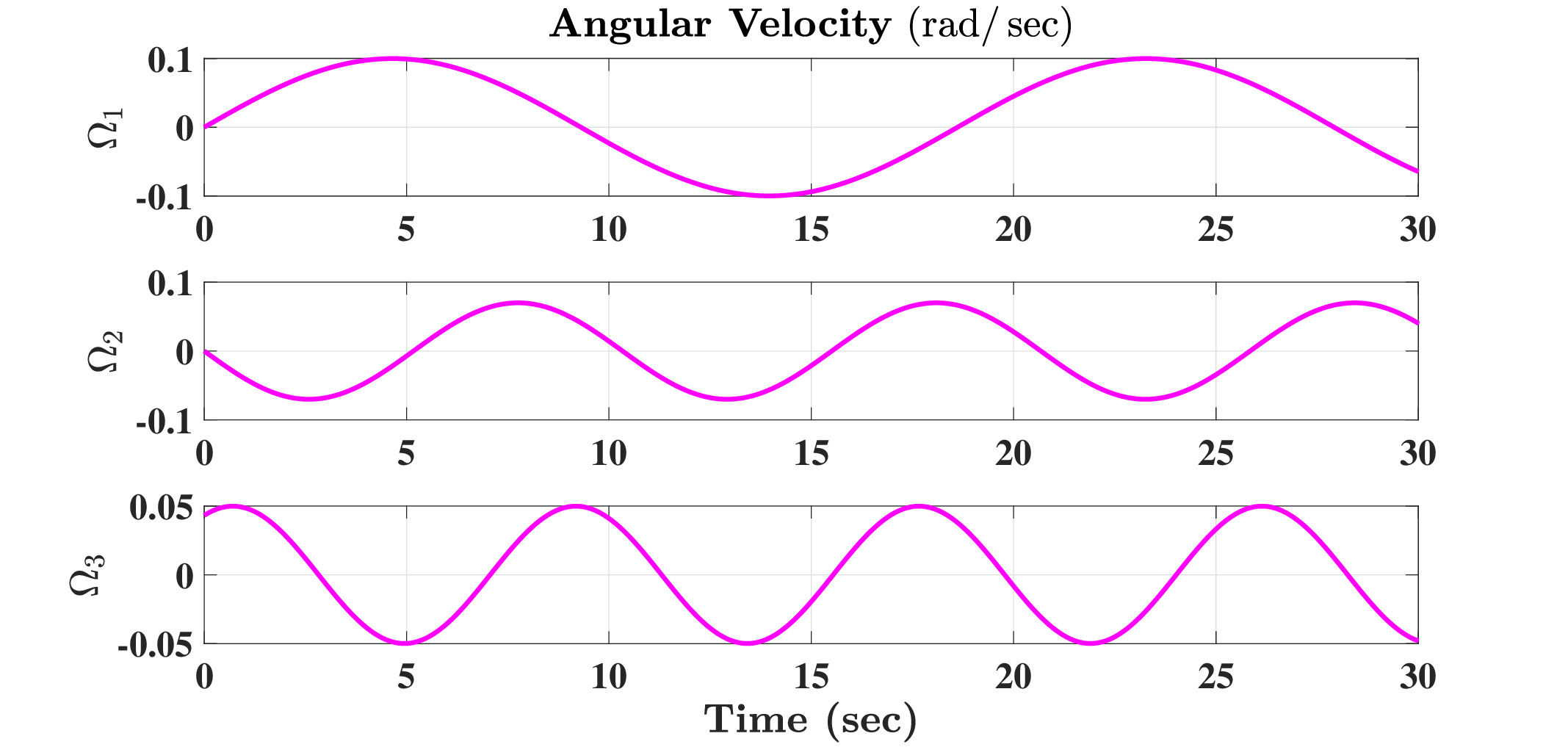}\caption{Angular velocity - Example 1}
	\label{fig:OVERVIEW_SO3_SIMU1} 
\end{figure}

\newpage{} 

\newpage{}

\subsection*{Example 2:}

In order to illustrate the problem of attitude parameterization through
Euler angles ($\xi$) obtained from Euler rates ($\dot{\xi}$), consider
the following angular velocity 
\begin{equation}
	\Omega=\left[\begin{array}{c}
		0.3{\rm sin}\left(0.8422t\right)\\
		0.21{\rm sin}\left(0.3682t+\pi\right)\\
		0.15{\rm sin}\left(1.4516t+\frac{\pi}{3}\right)
	\end{array}\right]\left({\rm rad/sec}\right)\label{eq:sim_ang-1}
\end{equation}
The angular velocity in \eqref{eq:sim_ang-1} is faster in rate and
higher in gains than \eqref{eq:sim_ang}. Let the initial attitude
be given by 
\begin{equation}
	R\left(0\right)=\left[\begin{array}{ccc}
		0.6679 & -0.1808 & 0.7219\\
		0.6552 & 0.6030 & -0.4551\\
		-0.3530 & 0.7770 & 0.5213
	\end{array}\right]\in\mathbb{SO}\left(3\right)\label{eq:R0-1}
\end{equation}
The initial conditions of $R\left(0\right)$ in terms of Euler angles,
Rodriguez vector and unit-quaternion are listed in Table \ref{tab:Initial-conditions2}.

\begin{table}[H]
	\centering{}\caption{\label{tab:Initial-conditions2}Initial conditions - Example 2}
	\begin{tabular}{ccc}
		\toprule 
		Representation  & Mapping  & Numerical values\tabularnewline
		\midrule
		\midrule 
		Euler angles  & $\begin{array}{c}
			\xi:\mathbb{SO}\left(3\right)\rightarrow\mathbb{R}^{3}\\
			R\left(0\right)\rightarrow\xi\left(0\right)
		\end{array}$  & $\xi\left(0\right)=\left[\begin{array}{c}
			\phi\\
			\theta\\
			\psi
		\end{array}\right]=\left[\begin{array}{c}
			56.1428\\
			20.6724\\
			44.4471
		\end{array}\right]\left({\rm deg}\right)$\tabularnewline
		\midrule
		Rodriguez vector  & $\begin{array}{c}
			\rho:\mathbb{SO}\left(3\right)\rightarrow\mathbb{R}^{3}\\
			R\left(0\right)\rightarrow\rho\left(0\right)
		\end{array}$  & $\rho\left(0\right)=\left[\begin{array}{c}
			\rho_{1}\\
			\rho_{2}\\
			\rho_{3}
		\end{array}\right]=\left[\begin{array}{c}
			0.4413\\
			0.3850\\
			0.2994
		\end{array}\right]$\tabularnewline
		\midrule
		Unit-quaternion  & $\begin{array}{c}
			Q:\mathbb{SO}\left(3\right)\rightarrow\mathbb{S}^{3}\\
			R\left(0\right)\rightarrow Q\left(0\right)
		\end{array}$  & $Q\left(0\right)=\left[\begin{array}{c}
			q_{0}\\
			q
		\end{array}\right]=\left[\begin{array}{c}
			0.8355\\
			0.3687\\
			0.3216\\
			0.2502
		\end{array}\right]$\tabularnewline
		\bottomrule
	\end{tabular}
\end{table}

\begin{figure}[h!]
	\centering{}\includegraphics[scale=0.43]{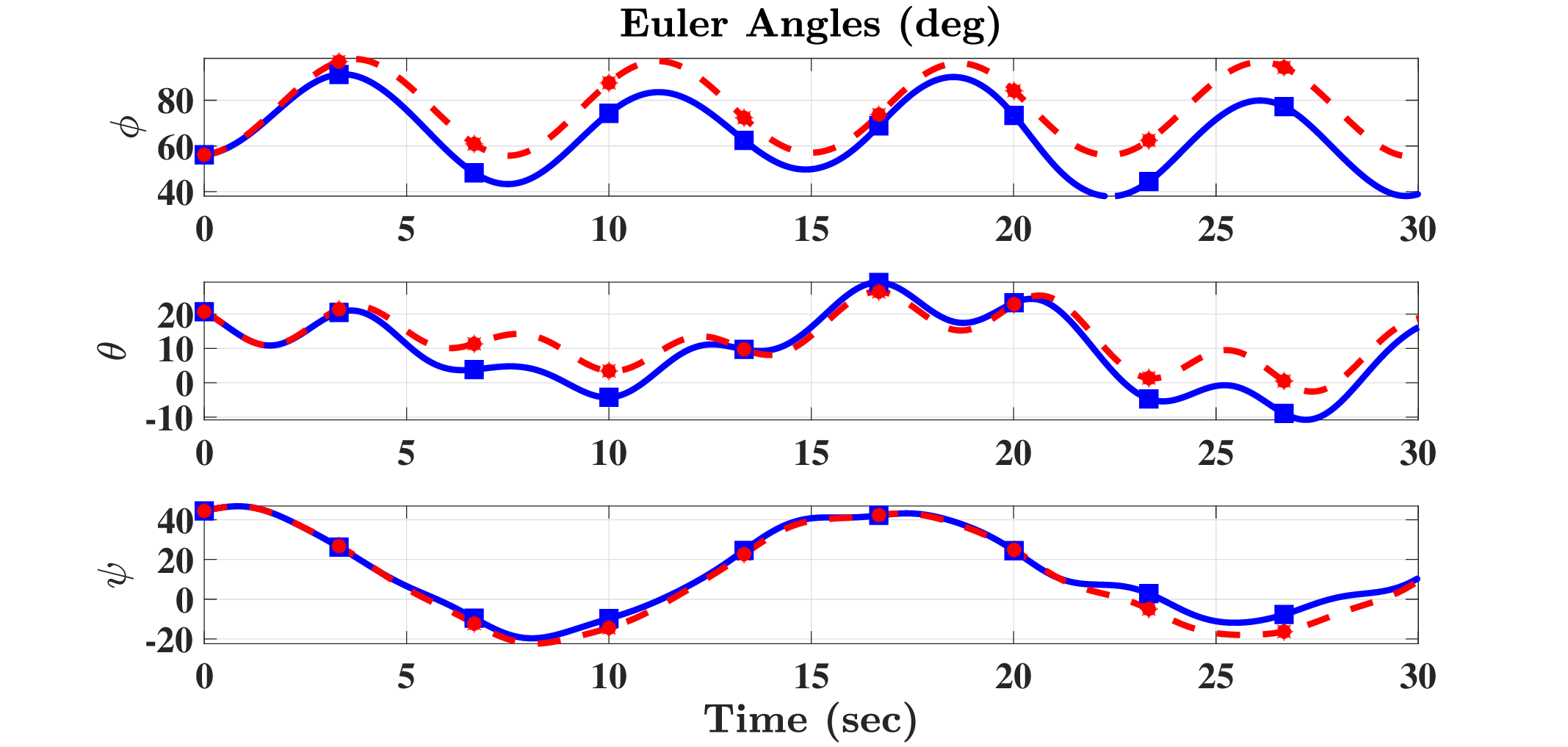}\caption{Euler Angles (\textcolor{blue}{$\dot{R}=R\left[\Omega\right]_{\times},\hspace{1em}R\rightarrow\xi$})
		vs (\textcolor{red}{$\dot{\xi}=\mathcal{J}\Omega$}) - Example 2}
	\label{fig:OVERVIEW_SO3_SIMU6} 
\end{figure}

\begin{figure}[h!]
	\centering{}\includegraphics[scale=0.43]{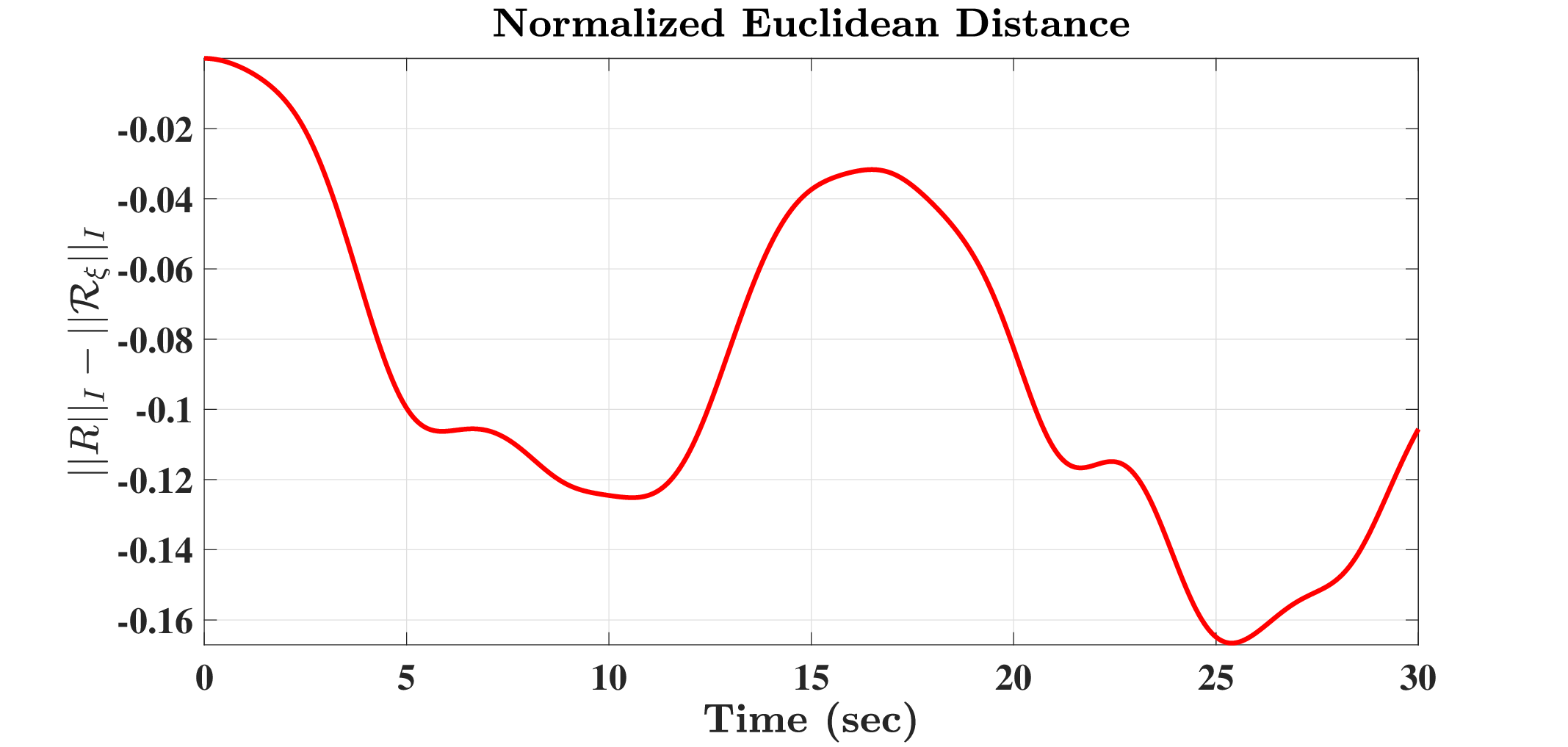}\caption{Error in normalized Euclidean distance associated with Euler Angles
		- Example 2}
	\label{fig:OVERVIEW_SO3_SIMU7} 
\end{figure}

\begin{figure}[h!]
	\centering{}\includegraphics[scale=0.43]{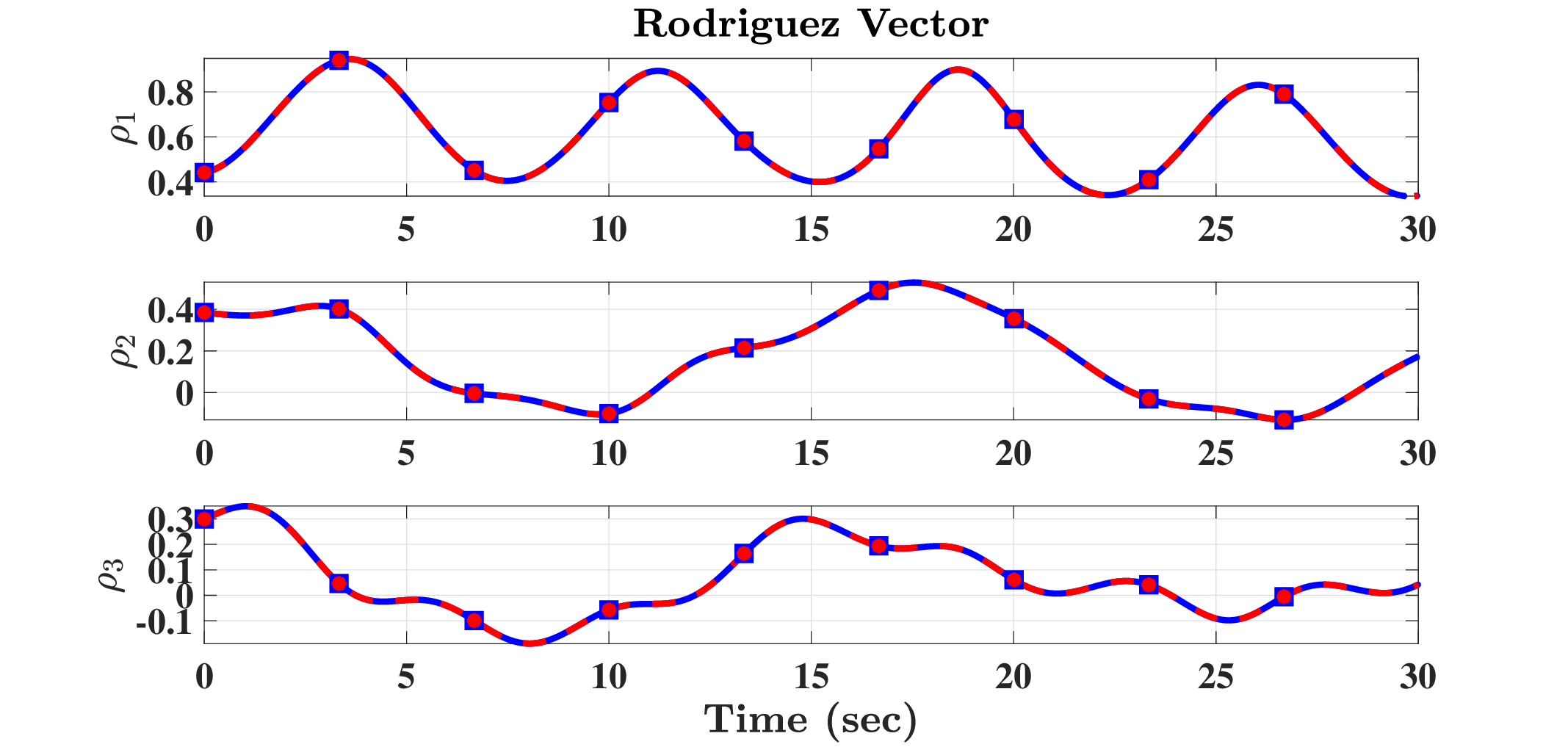}\caption{Rodriguez vector (\textcolor{blue}{$\dot{R}=R\left[\Omega\right]_{\times},\hspace{1em}R\rightarrow\rho$})
		vs (\textcolor{red}{$\dot{\rho}=\frac{1}{2}\left(\mathbf{I}_{3}+\left[\rho\right]_{\times}+\rho\rho^{\top}\right)\Omega$})
		- Example 2}
	\label{fig:OVERVIEW_SO3_SIMU8} 
\end{figure}

\begin{figure}[h!]
	\centering{}\includegraphics[scale=0.43]{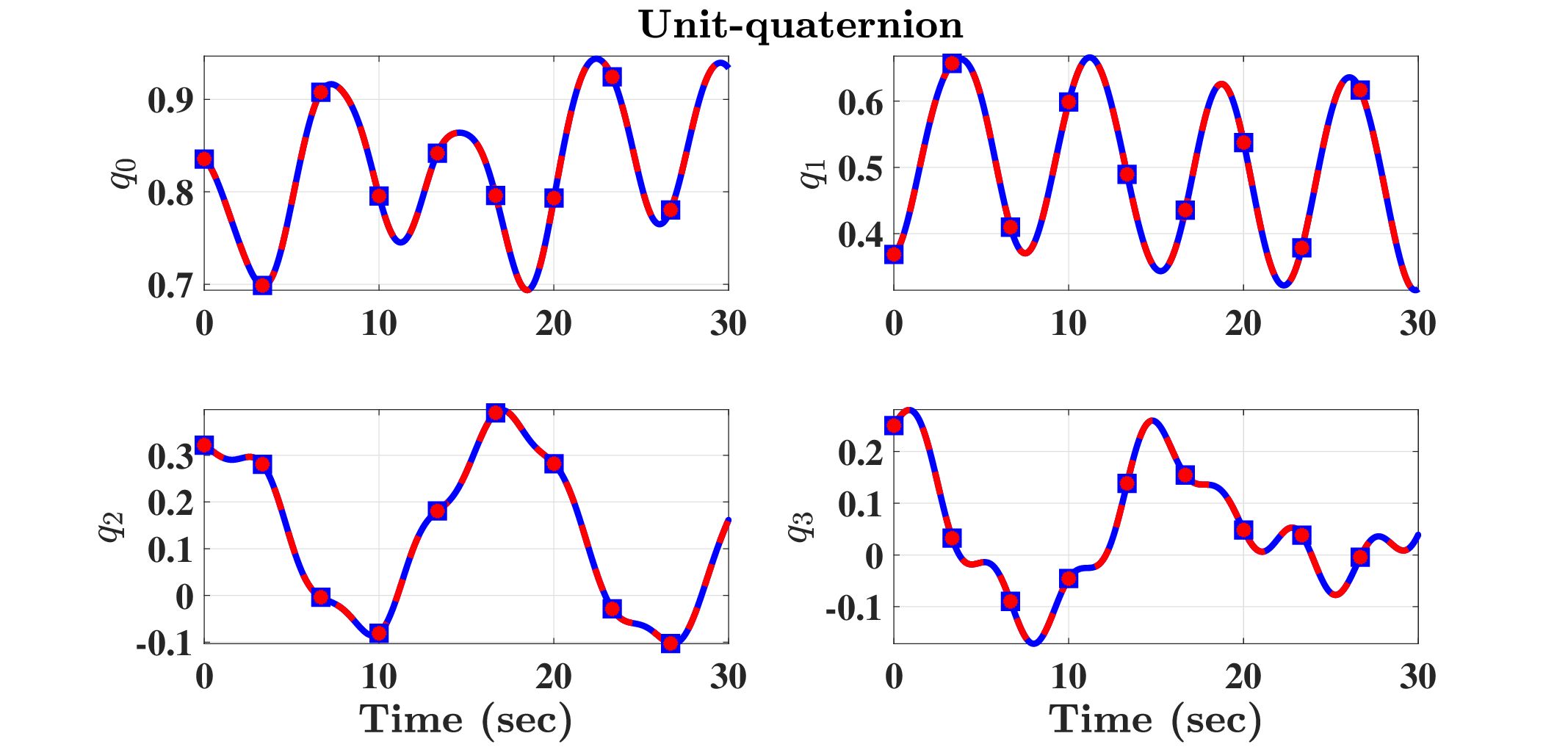}\caption{Unit-Quaternion (\textcolor{blue}{$\dot{R}=R\left[\Omega\right]_{\times},\hspace{1em}R\rightarrow Q$})
		vs (\textcolor{red}{$\dot{Q}=\frac{1}{2}\Gamma\left(\Omega\right)Q$})
		- Example 2}
	\label{fig:OVERVIEW_SO3_SIMU9} 
\end{figure}

The attitude $R$ has been obtained from the dynamics in \eqref{eq:OVERVIEW_SO3_dyn}
($\dot{R}=R\left[\Omega\right]_{\times}$) given $R\left(0\right)$
in \eqref{eq:R0-1} and the angular velocity in \eqref{eq:sim_ang-1}.
The mapping of the attitude $R$ from $\mathbb{SO}\left(3\right)$
to Euler angles, Rodriguez vector, and unit-quaternion is obtained
from \eqref{eq:SO3_EUL}, \eqref{eq:OVERVIEW_SO3_ROD}, and \eqref{eq:OVERVIEW_Q_comp_1},
respectively, and depicted in Figure \ref{fig:OVERVIEW_SO3_SIMU6},
\ref{fig:OVERVIEW_SO3_SIMU8}, and \ref{fig:OVERVIEW_SO3_SIMU9},
in blue colors. Euler angles obtained from Euler dynamics in \eqref{eq:OVERVIEW_EUL_dot}
are plotted in red color in Figure \ref{fig:OVERVIEW_SO3_SIMU6} against
Euler angles obtained through the mapping in \eqref{eq:SO3_EUL} from
$\mathbb{SO}\left(3\right)$ in blue colors. It can be clearly seen
in Figure \ref{fig:OVERVIEW_SO3_SIMU6} that Euler angles obtained
from Euler rate in \eqref{eq:SO3_EUL} failed to give the true Euler
angles obtained from the dynamics in \eqref{eq:OVERVIEW_SO3_dyn}.
Let Euler angles ($\xi$) obtained from Euler dynamics ($\dot{\xi}$)
in \eqref{eq:OVERVIEW_EUL_dot} be mapped such that $\xi:\mathbb{R}^{3}\rightarrow\mathcal{R}_{\xi}$.
In spite of $\mathcal{R}_{\xi}\in\mathbb{SO}\left(3\right)$ at every
time instant, $\mathcal{R}_{\xi}$ remains far from the true $R$.
Additionally, Figure \ref{fig:OVERVIEW_SO3_SIMU7} illustrates the
problem in Euler dynamics in terms of $\left\Vert R\right\Vert _{I}-\left\Vert \mathcal{R}_{\xi}\right\Vert _{I}$.
Figure \ref{fig:OVERVIEW_SO3_SIMU7} shows high error in the difference
between $\left\Vert R\right\Vert _{I}$ and $\left\Vert \mathcal{R}_{\xi}\right\Vert _{I}$.
On the other hand, the representation of Rodriguez vector and unit-quaternion
obtained from the dynamics in \eqref{eq:OVERVIEW_SO3_ROD}, and \eqref{eq:OVERVIEW_Q_comp_1},
respectively, is accurate and produces the same results as the mapping
of $R$ to Rodriquez vector and unit-quaternion obtained from \eqref{eq:OVERVIEW_SO3_dyn}.
Rodriguez vector obtained from Rodriguez vector dynamics in \eqref{eq:OVERVIEW_ROD_dot}
is plotted in red color in Figure \ref{fig:OVERVIEW_SO3_SIMU6} against
Rodriguez vector obtained through the mapping in \eqref{eq:OVERVIEW_SO3_ROD}
from $\mathbb{SO}\left(3\right)$ drawn in blue colors. Unit-quaternion
obtained from unit-quaternion dynamics in \eqref{eq:OVERVIEW_Q_dot}
is plotted in red color in Figure \ref{fig:OVERVIEW_SO3_SIMU8} against
unit-quaternion obtained through the mapping in \eqref{eq:OVERVIEW_Q_comp_1},
from $\mathbb{SO}\left(3\right)$ in blue colors. Table \ref{tab:colr_dynamics2}
summarizes the comparison between the three representations, while
the color notation are the same as those used in the above presented
plots. 
\begin{table}[H]
	\centering{}\caption{\label{tab:colr_dynamics2}Representation, color notation and related
		mapping - Example 2}
	\begin{tabular}{ccccc}
		\toprule 
		Representation  & Dynamics  & Color notation  & Mapping  & Figure \tabularnewline
		\midrule
		\midrule 
		\textcolor{blue}{$\mathbb{SO}\left(3\right)$} & \textcolor{blue}{$\dot{R}=R\left[\Omega\right]_{\times}$} & \textcolor{blue}{blue} & \textcolor{blue}{$R\left(t\right)\rightarrow\xi\left(t\right)$}  & \textcolor{blue}{9}\tabularnewline
		\midrule 
		&  &  & \textcolor{blue}{$R\left(t\right)\rightarrow\rho\left(t\right)$}  & \textcolor{blue}{11}\tabularnewline
		\midrule 
		&  &  & \textcolor{blue}{$R\left(t\right)\rightarrow Q\left(t\right)$}  & \textcolor{blue}{12}\tabularnewline
		\midrule 
		\textcolor{red}{Euler angles}  & \textcolor{red}{$\dot{\xi}=\mathcal{J}\Omega$}  & \textcolor{red}{Red}  & \textcolor{red}{$\int\dot{\xi}=\xi\left(t\right)$}  & \textcolor{red}{\ref{fig:OVERVIEW_SO3_SIMU6}}\tabularnewline
		\midrule 
		\textcolor{red}{Rodriguez vector}  & \textcolor{red}{$\dot{\rho}=\frac{1}{2}\left(\mathbf{I}_{3}+\left[\rho\right]_{\times}+\rho\rho^{\top}\right)\Omega$}  & \textcolor{red}{Red}  & \textcolor{red}{$\int\dot{\rho}=\rho\left(t\right)$}  & \textcolor{red}{\ref{fig:OVERVIEW_SO3_SIMU8}}\tabularnewline
		\midrule 
		\textcolor{red}{Unit-quaternion}  & \textcolor{red}{$\dot{Q}=\frac{1}{2}\Gamma\left(\Omega\right)Q$}  & \textcolor{red}{Red}  & \textcolor{red}{$\int\dot{Q}=Q\left(t\right)$}  & \textcolor{red}{\ref{fig:OVERVIEW_SO3_SIMU9} }\tabularnewline
		\bottomrule
	\end{tabular}
\end{table}

\begin{figure}[h!]
	\centering{}\includegraphics[scale=0.43]{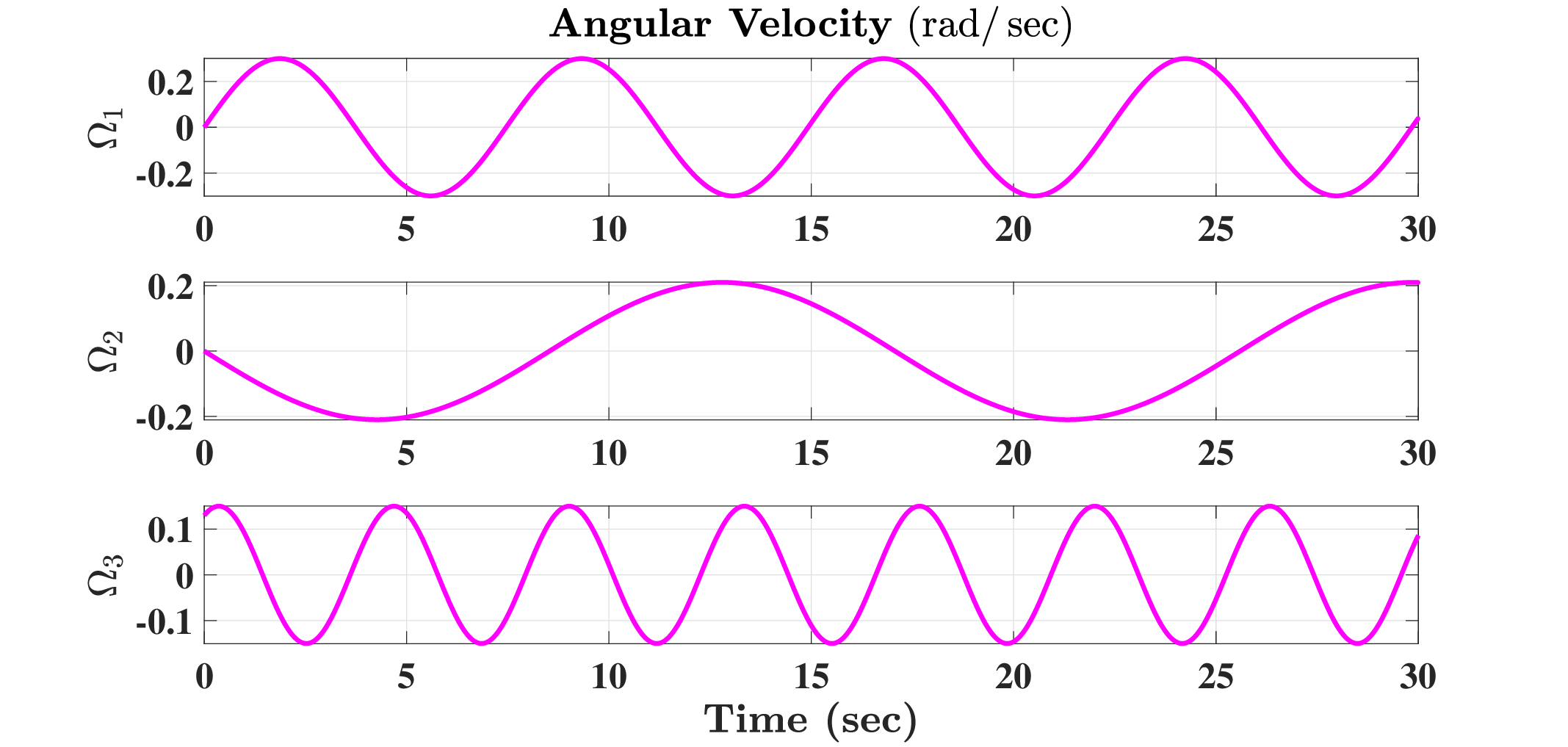}\caption{Angular velocity - Example 2}
	\label{fig:OVERVIEW_SO3_SIMU5} 
\end{figure}

\newpage{}

\newpage{}

\noindent\makebox[1\linewidth]{%
	\rule{0.6\textwidth}{1.4pt}%
}

\section{Conclusion \label{sec:Conclusion}}

This article presents the attitude configuration on the special orthogonal
group $\mathbb{SO}\left(3\right)$ and summarizes four popular methods
of attitude representation, namely Euler angles, angle-axis parameterization,
Rodriguez vector, and unit-quaternion. The $3\times3$ orthogonal
matrix gives a global and unique representation of the attitude, however
it does not have vector representation. Euler angles representation
is the most commonly used type of attitude representation, since it
can be easily visualized and understood. Nonetheless, it suffers from
singularity at certain configurations and resulting in unsuccessful
mapping from SO(3) to Euler angles. Rodriguez vector and angle-axis
parameterization are commonly used in the analysis of the attitude
filter and control design, and they have unique representation of
the attitude. Despite multiple advantages, both methods fail to represent
the attitude at quite a few configurations. In spite of the fact that
unit-quaternion does not have the singularity problem, it suffers
from non-uniqueness. The article gives a clear and detailed mapping
between different attitude parameterizations. A summary of the mapping
is provided in Table \ref{tab:SO3}, \ref{tab:Angle-axis}, \ref{tab:Rodriguez-vector},
and \ref{tab:Unit-Quaternion}. Also, the article provides some important
results which could be used to achieve successful estimation and/or
control process of a rigid-body rotating and/or moving in space. 

\noindent\makebox[1\linewidth]{%
	\rule{0.6\textwidth}{1.4pt}%
}
\begin{center}
	\begin{table}[H]
		\begin{centering}
			\caption{\label{tab:SO3}Mapping: $\mathbb{SO}\left(3\right)$ and other parameterizations}
			\par\end{centering}
		\centering{}%
		\begin{tabular}{l|l}
			\hline 
			\noalign{\vskip\doublerulesep}
			\textcolor{blue}{$R:\mathbb{R}^{3}\rightarrow\mathbb{SO}\left(3\right)$ } & \textcolor{blue}{$\xi:\mathbb{SO}\left(3\right)\rightarrow\mathbb{R}^{3}$ }\tabularnewline[\doublerulesep]
			\hline 
			\hline 
			\noalign{\vskip\doublerulesep}
			$\begin{array}{cc}
			& s=\sin\text{, }c=\cos\text{, }t=\tan\\
			R= & \left[\begin{array}{ccc}
			c\theta c\psi & -c\phi s\psi+s\phi s\theta c\psi & s\phi s\psi+c\phi s\theta c\psi\\
			c\theta s\psi & c\phi c\psi+s\phi s\theta s\psi & -s\phi c\psi+c\phi s\theta s\psi\\
			-s\theta & s\phi c\theta & c\phi c\theta
			\end{array}\right]
			\end{array}$ & $\left[\begin{array}{c}
			\phi\\
			\theta\\
			\psi
			\end{array}\right]=\left[\begin{array}{c}
			\arctan\left(\frac{R_{\left(3,2\right)}}{R_{\left(3,3\right)}}\right)\\
			\arctan\left(\frac{-R_{\left(3,1\right)}}{\sqrt{R_{\left(3,2\right)}^{2}+R_{\left(3,3\right)}^{2}}}\right)\\
			\arctan\left(\frac{R_{\left(2,1\right)}}{R_{\left(1,1\right)}}\right)
			\end{array}\right]$\tabularnewline[\doublerulesep]
			\hline 
			\noalign{\vskip\doublerulesep}
			\textcolor{blue}{$R:\mathbb{R}\times\mathbb{S}^{2}\rightarrow\mathbb{SO}\left(3\right)$ } & \textcolor{blue}{$\alpha,u:\mathbb{SO}\left(3\right)\rightarrow\mathbb{R}\times\mathbb{S}^{2}$ }\tabularnewline[\doublerulesep]
			\hline 
			\hline 
			\noalign{\vskip\doublerulesep}
			\multirow{2}{*}{$\mathcal{R}_{\alpha}\left(\alpha,u\right)=\mathbf{I}_{3}+\sin\left(\alpha\right)\left[u\right]_{\times}+\left(1-\cos\left(\theta\right)\right)\left[u\right]_{\times}^{2}$ } & $\alpha={\rm cos}^{-1}\left(\frac{1}{2}\left({\rm Tr}\left\{ R\right\} -1\right)\right)$\tabularnewline[\doublerulesep]
			\noalign{\vskip\doublerulesep}
			& $u=\frac{1}{{\rm sin}\left(\alpha\right)}\mathbf{vex}\left(\boldsymbol{\mathcal{P}}_{a}\left(R\right)\right)$\tabularnewline[\doublerulesep]
			\hline 
			\noalign{\vskip\doublerulesep}
			\textcolor{blue}{$R:\mathbb{R}^{3}\rightarrow\mathbb{SO}\left(3\right)$ } & \textcolor{blue}{$\rho:\mathbb{SO}\left(3\right)\rightarrow\mathbb{R}^{3}$ }\tabularnewline[\doublerulesep]
			\hline 
			\hline 
			\noalign{\vskip\doublerulesep}
			\multirow{2}{*}{$\mathcal{R}_{\rho}\left(\rho\right)=\frac{\left(1-\left\Vert \rho\right\Vert ^{2}\right)\mathbf{I}_{3}+2\rho\rho^{\top}+2\left[\rho\right]_{\times}}{1+\left\Vert \rho\right\Vert ^{2}}$} & $\rho=\mathbf{vex}\left(\left(R-\mathbf{I}_{3}\right)\left(R+\mathbf{I}_{3}\right)^{-1}\right)$\tabularnewline[\doublerulesep]
			\noalign{\vskip\doublerulesep}
			& $\rho=\frac{\mathbf{vex}\left(\boldsymbol{\mathcal{P}}_{a}\left(R\right)\right)}{2\left(1-\left\Vert R\right\Vert _{I}\right)}$ \tabularnewline[\doublerulesep]
			\hline 
			\noalign{\vskip\doublerulesep}
			\textcolor{blue}{$R:\mathbb{S}^{3}\rightarrow\mathbb{SO}\left(3\right)$ } & \textcolor{blue}{$Q:\mathbb{SO}\left(3\right)\rightarrow\mathbb{S}^{3}$ }\tabularnewline[\doublerulesep]
			\hline 
			\hline 
			\noalign{\vskip\doublerulesep}
			$\mathcal{R}_{Q}\left(Q\right)=\mathbf{I}_{3}+2q_{0}\left[q\right]_{\times}+2\left[q\right]_{\times}^{2}$ & \multirow{2}{*}{See \eqref{eq:OVERVIEW_Q_comp_1}, \eqref{eq:OVERVIEW_Q_comp_2},
				\eqref{eq:OVERVIEW_Q_comp_3}, and \eqref{eq:OVERVIEW_Q_comp_4}}\tabularnewline[\doublerulesep]
			\noalign{\vskip\doublerulesep}
			$\mathcal{R}_{Q}\left(Q\right)=\left(q_{0}^{2}-\left\Vert q\right\Vert ^{2}\right)\mathbf{I}_{3}+2q_{0}\left[q\right]_{\times}+2qq^{\top}$ & \tabularnewline[\doublerulesep]
			\hline 
		\end{tabular}
	\end{table}
	\par\end{center}

\begin{center}
	\begin{table}[H]
		\begin{centering}
			\caption{\label{tab:Angle-axis}Mapping: Angle-axis and other parameterizations}
			\par\end{centering}
		\centering{}%
		\begin{tabular}{l|l}
			\hline 
			\noalign{\vskip\doublerulesep}
			\textcolor{blue}{$\alpha,u:\mathbb{SO}\left(3\right)\rightarrow\mathbb{R}\times\mathbb{S}^{2}$ } & \textcolor{blue}{$R:\mathbb{R}\times\mathbb{S}^{2}\rightarrow\mathbb{SO}\left(3\right)$ }\tabularnewline[\doublerulesep]
			\hline 
			\hline 
			\noalign{\vskip\doublerulesep}
			$\alpha={\rm cos}^{-1}\left(\frac{1}{2}\left({\rm Tr}\left\{ R\right\} -1\right)\right)$ & \multirow{2}{*}{$\mathcal{R}_{\alpha}\left(\alpha,u\right)=\mathbf{I}_{3}+\sin\left(\alpha\right)\left[u\right]_{\times}+\left(1-\cos\left(\theta\right)\right)\left[u\right]_{\times}^{2}$ }\tabularnewline[\doublerulesep]
			\noalign{\vskip\doublerulesep}
			$u=\frac{1}{{\rm sin}\left(\alpha\right)}\mathbf{vex}\left(\boldsymbol{\mathcal{P}}_{a}\left(R\right)\right)$ & \tabularnewline[\doublerulesep]
			\hline 
			\noalign{\vskip\doublerulesep}
			\textcolor{blue}{$\alpha,u:\mathbb{R}^{3}\rightarrow\mathbb{R}\times\mathbb{S}^{2}$ } & \textcolor{blue}{$\rho:\mathbb{R}\times\mathbb{S}^{2}\rightarrow\mathbb{R}^{3}$ }\tabularnewline[\doublerulesep]
			\hline 
			\hline 
			\noalign{\vskip\doublerulesep}
			\multicolumn{1}{l|}{$\alpha=2{\rm tan}^{-1}\left(\left\Vert \rho\right\Vert \right)=2{\rm sin}^{-1}\left(\frac{\left\Vert \rho\right\Vert }{\sqrt{1+\left\Vert \rho\right\Vert ^{2}}}\right)$} & \multirow{2}{*}{$\rho={\rm tan}\left(\frac{\alpha}{2}\right)u$ }\tabularnewline[\doublerulesep]
			\noalign{\vskip\doublerulesep}
			$u={\rm cot}\left(\frac{\alpha}{2}\right)\rho$ & \tabularnewline[\doublerulesep]
			\hline 
			\noalign{\vskip\doublerulesep}
			\textcolor{blue}{$\alpha,u:\mathbb{S}^{3}\rightarrow\mathbb{R}\times\mathbb{S}^{2}$ } & \textcolor{blue}{$Q:\mathbb{R}\times\mathbb{S}^{2}\rightarrow\mathbb{S}^{3}$ }\tabularnewline[\doublerulesep]
			\hline 
			\hline 
			\noalign{\vskip\doublerulesep}
			\multicolumn{1}{l|}{$\alpha=2{\rm cos}^{-1}\left(q_{0}\right)$} & \multirow{2}{*}{$Q=\left[\begin{array}{c}
				\cos\left(\alpha/2\right)\\
				u\sin\left(\alpha/2\right)
				\end{array}\right]$ }\tabularnewline[\doublerulesep]
			\noalign{\vskip\doublerulesep}
			$u=\frac{1}{\sin\left(\alpha/2\right)}q$ & \tabularnewline[\doublerulesep]
			\hline 
		\end{tabular}
	\end{table}
	\par\end{center}

\begin{center}
	\begin{table}[H]
		\begin{centering}
			\caption{\label{tab:Rodriguez-vector}Mapping: Rodriguez vector and other parameterizations}
			\par\end{centering}
		\centering{}%
		\begin{tabular}{l|l}
			\hline 
			\noalign{\vskip\doublerulesep}
			\textcolor{blue}{$\rho:\mathbb{SO}\left(3\right)\rightarrow\mathbb{R}^{3}$ } & \textcolor{blue}{$R:\mathbb{R}^{3}\rightarrow\mathbb{SO}\left(3\right)$ }\tabularnewline
			\hline 
			\hline 
			\noalign{\vskip\doublerulesep}
			\multicolumn{1}{l|}{$\rho=\mathbf{vex}\left(\left(R-\mathbf{I}_{3}\right)\left(R+\mathbf{I}_{3}\right)^{-1}\right)$} & \multirow{2}{*}{$\mathcal{R}_{\rho}\left(\rho\right)=\frac{\left(1-\left\Vert \rho\right\Vert ^{2}\right)\mathbf{I}_{3}+2\rho\rho^{\top}+2\left[\rho\right]_{\times}}{1+\left\Vert \rho\right\Vert ^{2}}$}\tabularnewline[\doublerulesep]
			\noalign{\vskip\doublerulesep}
			$\rho=\frac{\mathbf{vex}\left(\boldsymbol{\mathcal{P}}_{a}\left(R\right)\right)}{2\left(1-\left\Vert R\right\Vert _{I}\right)}$  & \tabularnewline[\doublerulesep]
			\hline 
			\noalign{\vskip\doublerulesep}
			& \textcolor{blue}{$\xi:\mathbb{R}^{3}\rightarrow\mathbb{R}^{3}$}\tabularnewline[\doublerulesep]
			\hline 
			\hline 
			\noalign{\vskip\doublerulesep}
			& $\left[\begin{array}{c}
			\phi\\
			\theta\\
			\psi
			\end{array}\right]=\left[\begin{array}{c}
			\arctan\left(\frac{2\rho_{2}\rho_{3}+2\rho_{1}}{1+\rho_{3}^{2}-\rho_{1}^{2}-\rho_{2}^{2}}\right)\\
			\arctan\left(\frac{2\rho_{2}-2\rho_{1}\rho_{3}}{\sqrt{4\left(\rho_{2}\rho_{3}+\rho_{1}\right)^{2}+\left(1+\rho_{3}^{2}-\rho_{1}^{2}-\rho_{2}^{2}\right)^{2}}}\right)\\
			\arctan\left(\frac{2\rho_{1}\rho_{2}+2\rho_{3}}{1+\rho_{1}^{2}-\rho_{2}^{2}-\rho_{3}^{2}}\right)
			\end{array}\right]$\tabularnewline[\doublerulesep]
			\hline 
			\noalign{\vskip\doublerulesep}
			\textcolor{blue}{$\rho:\mathbb{R}\times\mathbb{S}^{2}\rightarrow\mathbb{R}^{3}$ } & \textcolor{blue}{$\alpha,u:\mathbb{R}^{3}\rightarrow\mathbb{R}\times\mathbb{S}^{2}$ }\tabularnewline[\doublerulesep]
			\hline 
			\hline 
			\noalign{\vskip\doublerulesep}
			\multirow{2}{*}{$\rho={\rm tan}\left(\frac{\alpha}{2}\right)u$ } & $\alpha=2{\rm tan}^{-1}\left(\left\Vert \rho\right\Vert \right)=2{\rm sin}^{-1}\left(\frac{\left\Vert \rho\right\Vert }{\sqrt{1+\left\Vert \rho\right\Vert ^{2}}}\right)$\tabularnewline[\doublerulesep]
			\noalign{\vskip\doublerulesep}
			& $u={\rm cot}\left(\frac{\alpha}{2}\right)\rho$\tabularnewline[\doublerulesep]
			\hline 
			\noalign{\vskip\doublerulesep}
			\textcolor{blue}{$\rho:\mathbb{S}^{3}\rightarrow\mathbb{R}^{3}$ } & \textcolor{blue}{$Q:\mathbb{R}^{3}\rightarrow\mathbb{S}^{3}$}\tabularnewline[\doublerulesep]
			\hline 
			\hline 
			\noalign{\vskip\doublerulesep}
			\multirow{2}{*}{$\rho=q/q_{0}$} & \multicolumn{1}{l}{$q_{0}=\pm\frac{1}{\sqrt{1+\left\Vert \rho\right\Vert ^{2}}}$}\tabularnewline[\doublerulesep]
			\noalign{\vskip\doublerulesep}
			& $q=q_{0}\rho$\tabularnewline[\doublerulesep]
			\hline 
		\end{tabular}
	\end{table}
	\par\end{center}

\begin{center}
	\begin{table}[H]
		\begin{centering}
			\caption{\label{tab:Unit-Quaternion}Mapping: Unit-quaternion and other parameterizations}
			\par\end{centering}
		\centering{}%
		\begin{tabular}{l|l}
			\hline 
			\noalign{\vskip\doublerulesep}
			\textcolor{blue}{$Q:\mathbb{SO}\left(3\right)\rightarrow\mathbb{S}^{3}$ } & \textcolor{blue}{$R:\mathbb{S}^{3}\rightarrow\mathbb{SO}\left(3\right)$ }\tabularnewline
			\hline 
			\hline 
			\noalign{\vskip\doublerulesep}
			\multirow{2}{*}{See \eqref{eq:OVERVIEW_Q_comp_1}, \eqref{eq:OVERVIEW_Q_comp_2},
				\eqref{eq:OVERVIEW_Q_comp_3}, and \eqref{eq:OVERVIEW_Q_comp_4}} & $\mathcal{R}_{Q}\left(Q\right)=\mathbf{I}_{3}+2q_{0}\left[q\right]_{\times}+2\left[q\right]_{\times}^{2}$\tabularnewline[\doublerulesep]
			\noalign{\vskip\doublerulesep}
			& $\mathcal{R}_{Q}\left(Q\right)=\left(q_{0}^{2}-\left\Vert q\right\Vert ^{2}\right)\mathbf{I}_{3}+2q_{0}\left[q\right]_{\times}+2qq^{\top}$\tabularnewline[\doublerulesep]
			\hline 
			\noalign{\vskip\doublerulesep}
			& \textcolor{blue}{$\xi:\mathbb{S}^{3}\rightarrow\mathbb{R}^{3}$ }\tabularnewline[\doublerulesep]
			\hline 
			\hline 
			\noalign{\vskip\doublerulesep}
			& $\left[\begin{array}{c}
			\phi\\
			\theta\\
			\psi
			\end{array}\right]=\left[\begin{array}{c}
			\arctan\left(\frac{2\left(q_{3}q_{2}+q_{0}q_{1}\right)}{1-2\left(q_{1}^{2}+q_{2}^{2}\right)}\right)\\
			\arcsin\left(2\left(q_{0}q_{2}-q_{3}q_{1}\right)\right)\\
			\arctan\left(\frac{2\left(q_{2}q_{1}+q_{0}q_{3}\right)}{1-2\left(q_{2}^{2}+q_{3}^{2}\right)}\right)
			\end{array}\right]$\tabularnewline[\doublerulesep]
			\hline 
			\noalign{\vskip\doublerulesep}
			\textcolor{blue}{$Q:\mathbb{R}\times\mathbb{S}^{2}\rightarrow\mathbb{S}^{3}$ } & \textcolor{blue}{$\alpha,u:\mathbb{S}^{3}\rightarrow\mathbb{R}\times\mathbb{S}^{2}$ }\tabularnewline[\doublerulesep]
			\hline 
			\hline 
			\noalign{\vskip\doublerulesep}
			\multirow{2}{*}{$Q=\left[\begin{array}{c}
				\cos\left(\alpha/2\right)\\
				u\sin\left(\alpha/2\right)
				\end{array}\right]$ } & $\alpha=2{\rm cos}^{-1}\left(q_{0}\right)$\tabularnewline[\doublerulesep]
			\noalign{\vskip\doublerulesep}
			& $u=\frac{1}{\sin\left(\alpha/2\right)}q$\tabularnewline[\doublerulesep]
			\hline 
			\noalign{\vskip\doublerulesep}
			\textcolor{blue}{$Q:\mathbb{R}^{3}\rightarrow\mathbb{S}^{3}$ } & \textcolor{blue}{$\rho:\mathbb{S}^{3}\rightarrow\mathbb{R}^{3}$}\tabularnewline[\doublerulesep]
			\hline 
			\hline 
			\noalign{\vskip\doublerulesep}
			$q_{0}=\pm\frac{1}{\sqrt{1+\left\Vert \rho\right\Vert ^{2}}}$ & \multirow{2}{*}{$\rho=q/q_{0}$}\tabularnewline[\doublerulesep]
			\noalign{\vskip\doublerulesep}
			$q=q_{0}\rho$ & \tabularnewline[\doublerulesep]
			\hline 
		\end{tabular}
	\end{table}
	\par\end{center}

\noindent\makebox[1\linewidth]{%
	\rule{0.6\textwidth}{1.4pt}%
}

\section*{Acknowledgment}

The author would
like to thank \textbf{Maria Shaposhnikova} for proofreading the article.

\bibliographystyle{IEEEtran}
\bibliography{bib_Overview_SO3}

\begin{thebibliography}{10}
\providecommand{\url}[1]{#1}
\csname url@samestyle\endcsname
\providecommand{\newblock}{\relax}
\providecommand{\bibinfo}[2]{#2}
\providecommand{\BIBentrySTDinterwordspacing}{\spaceskip=0pt\relax}
\providecommand{\BIBentryALTinterwordstretchfactor}{4}
\providecommand{\BIBentryALTinterwordspacing}{\spaceskip=\fontdimen2\font plus
\BIBentryALTinterwordstretchfactor\fontdimen3\font minus
  \fontdimen4\font\relax}
\providecommand{\BIBforeignlanguage}[2]{{%
\expandafter\ifx\csname l@#1\endcsname\relax
\typeout{** WARNING: IEEEtran.bst: No hyphenation pattern has been}%
\typeout{** loaded for the language `#1'. Using the pattern for}%
\typeout{** the default language instead.}%
\else
\language=\csname l@#1\endcsname
\fi
#2}}
\providecommand{\BIBdecl}{\relax}
\BIBdecl

\bibitem{hashim2020AtiitudeSurvey}
H.~A. Hashim, ``Attitude determination and estimation using vector
  observations: Review, challenges and comparative results,'' \emph{ar{X}iv
  preprint ar{X}iv:2001.03787}, 2020.

\bibitem{stuelpnagel1964parametrization}
J.~Stuelpnagel, ``On the parametrization of the three-dimensional rotation
  group,'' \emph{SIAM review}, vol.~6, no.~4, pp. 422--430, 1964.

\bibitem{wen1991attitude}
J.-Y. Wen and K.~Kreutz-Delgado, ``The attitude control problem,'' \emph{IEEE
  Transactions on Automatic control}, vol.~36, no.~10, pp. 1148--1162, 1991.

\bibitem{shuster1993survey}
M.~D. Shuster, ``A survey of attitude representations,'' \emph{Navigation},
  vol.~8, no.~9, pp. 439--517, 1993.

\bibitem{chaturvedi2011rigid}
N.~A. Chaturvedi, A.~K. Sanyal, and N.~H. McClamroch, ``Rigid-body attitude
  control,'' \emph{IEEE Control Systems}, vol.~31, no.~3, pp. 30--51, 2011.

\bibitem{kang2009linear}
Y.~Kang and J.~K. Hedrick, ``Linear tracking for a fixed-wing uav using
  nonlinear model predictive control,'' \emph{IEEE Transactions on Control
  Systems Technology}, vol.~17, no.~5, pp. 1202--1210, 2009.

\bibitem{mohamed2014improved}
H.~A.~H. Mohamed, ``Improved robust adaptive control of high-order nonlinear
  systems with guaranteed performance,'' \emph{M. Sc, King Fahd University Of
  Petroleum \& Minerals}, vol.~1, 2014.

\bibitem{stevens2015aircraft}
B.~L. Stevens, F.~L. Lewis, and E.~N. Johnson, \emph{Aircraft control and
  simulation: dynamics, controls design, and autonomous systems}.\hskip 1em
  plus 0.5em minus 0.4em\relax John Wiley \& Sons, 2015.

\bibitem{yeh2002nonlinear}
H.-H. Yeh, E.~Nelson, and A.~Sparks, ``Nonlinear tracking control for satellite
  formations,'' \emph{Journal of Guidance, Control, and Dynamics}, vol.~25,
  no.~2, pp. 376--386, 2002.

\bibitem{chwa2010tracking}
D.~Chwa, ``Tracking control of differential-drive wheeled mobile robots using a
  backstepping-like feedback linearization,'' \emph{IEEE Transactions on
  Systems, Man, and Cybernetics-Part A: Systems and Humans}, vol.~40, no.~6,
  pp. 1285--1295, 2010.

\bibitem{hashim2018SO3Stochastic}
H.~A. Hashim, L.~J. Brown, and K.~McIsaac, ``Nonlinear stochastic attitude
  filters on the special orthogonal group 3: Ito and stratonovich,'' \emph{IEEE
  Transactions on Systems, Man, and Cybernetics: Systems}, pp. 1--13, 2018.

\bibitem{markley2003attitude}
F.~L. Markley, ``Attitude error representations for kalman filtering,''
  \emph{Journal of guidance, control, and dynamics}, vol.~26, no.~2, pp.
  311--317, 2003.

\bibitem{crassidis2007survey}
J.~L. Crassidis, F.~L. Markley, and Y.~Cheng, ``Survey of nonlinear attitude
  estimation methods,'' \emph{Journal of guidance, control, and dynamics},
  vol.~30, no.~1, pp. 12--28, 2007.

\bibitem{mahony2008nonlinear}
R.~Mahony, T.~Hamel, and J.-M. Pflimlin, ``Nonlinear complementary filters on
  the special orthogonal group,'' \emph{IEEE Transactions on Automatic
  Control}, vol.~53, no.~5, pp. 1203--1218, 2008.

\bibitem{hashim2019SE3Det}
H.~A. Hashim, L.~J. Brown, and K.~McIsaac, ``Nonlinear pose filters on the
  special euclidean group se(3) with guaranteed transient and steady-state
  performance,'' \emph{IEEE Transactions on Systems, Man, and Cybernetics:
  Systems}, pp. 1--14, 2019.

\bibitem{mohamed2019filters}
H.~A.~H. Mohamed, ``Nonlinear attitude and pose filters with superior
  convergence properties,'' \emph{Ph. D, Western University}, 2019.

\bibitem{hashim2017neuro}
H.~A. Hashim, S.~El-Ferik, and F.~L. Lewis, ``Neuro-adaptive cooperative
  tracking control with prescribed performance of unknown higher-order
  nonlinear multi-agent systems,'' \emph{International Journal of Control}, pp.
  1--16, 2017.

\bibitem{hashim2017adaptive}
------, ``Adaptive synchronisation of unknown nonlinear networked systems with
  prescribed performance,'' \emph{International Journal of Systems Science},
  vol.~48, no.~4, pp. 885--898, 2017.

\bibitem{lefferts1982kalman}
E.~J. Lefferts, F.~L. Markley, and M.~D. Shuster, ``Kalman filtering for
  spacecraft attitude estimation,'' \emph{Journal of Guidance, Control, and
  Dynamics}, vol.~5, no.~5, pp. 417--429, 1982.

\bibitem{goddard1998pose}
J.~S. Goddard and M.~A. Abidi, ``Pose and motion estimation using dual
  quaternion-based extended kalman filtering,'' in \emph{Three-Dimensional
  Image Capture and Applications}, vol. 3313.\hskip 1em plus 0.5em minus
  0.4em\relax International Society for Optics and Photonics, 1998, pp.
  189--201.

\bibitem{choukroun2006novel}
D.~Choukroun, I.~Y. Bar-Itzhack, and Y.~Oshman, ``Novel quaternion kalman
  filter,'' \emph{IEEE Transactions on Aerospace and Electronic Systems},
  vol.~42, no.~1, pp. 174--190, 2006.

\bibitem{cheon2007unscented}
Y.-J. Cheon and J.-H. Kim, ``Unscented filtering in a unit quaternion space for
  spacecraft attitude estimation,'' in \emph{Industrial Electronics, 2007. ISIE
  2007. IEEE International Symposium on}.\hskip 1em plus 0.5em minus
  0.4em\relax IEEE, 2007, pp. 66--71.

\bibitem{barrau2015intrinsic}
A.~Barrau and S.~Bonnabel, ``Intrinsic filtering on lie groups with
  applications to attitude estimation,'' \emph{IEEE Transactions on Automatic
  Control}, vol.~60, no.~2, pp. 436--449, 2015.

\bibitem{hashim2019SO3Det}
H.~A. Hashim, L.~J. Brown, and K.~McIsaac, ``Guaranteed performance of
  nonlinear attitude filters on the special orthogonal group {SO}(3),''
  \emph{IEEE Access}, vol.~7, no.~1, pp. 3731--3745, 2019.

\bibitem{hashim2018Conf1}
------, ``Nonlinear explicit stochastic attitude filter on {SO}(3),'' in
  \emph{Proceedings of the 57th {IEEE} conference on Decision and Control
  ({CDC})}.\hskip 1em plus 0.5em minus 0.4em\relax IEEE, 2018, pp. 1210 --1216.

\bibitem{hashim2018SE3Stochastic}
------, ``Nonlinear stochastic position and attitude filter on the special
  euclidean group 3,'' \emph{Journal of the Franklin Institute}, pp. 1--27,
  2018.

\bibitem{hashim2019Conf1}
------, ``Guaranteed performance of nonlinear pose filter on {SE}(3),'' in
  \emph{Proceedings of the American Control Conference ({ACC})}, 2019, pp.
  1--6.

\bibitem{grood1983joint}
E.~S. Grood and W.~J. Suntay, ``A joint coordinate system for the clinical
  description of three-dimensional motions: application to the knee,''
  \emph{Journal of biomechanical engineering}, vol. 105, no.~2, pp. 136--144,
  1983.

\bibitem{mokhtari2004dynamic}
A.~Mokhtari and A.~Benallegue, ``Dynamic feedback controller of euler angles
  and wind parameters estimation for a quadrotor unmanned aerial vehicle,'' in
  \emph{Robotics and Automation, 2004. Proceedings. ICRA'04. 2004 IEEE
  International Conference on}, vol.~3.\hskip 1em plus 0.5em minus 0.4em\relax
  IEEE, 2004, pp. 2359--2366.

\bibitem{spong2008robot}
M.~W. Spong and M.~Vidyasagar, \emph{Robot dynamics and control}.\hskip 1em
  plus 0.5em minus 0.4em\relax John Wiley \& Sons, 2008.

\bibitem{murray1994mathematicall}
R.~M. Murray, \emph{A mathematical introduction to robotic manipulation}.\hskip
  1em plus 0.5em minus 0.4em\relax CRC press, 1994.

\bibitem{rodrigues1840lois}
O.~Rodrigues, \emph{Des lois g{\'e}om{\'e}triques qui r{\'e}gissent les
  d{\'e}placements d'un syst{\`e}me solide dans l'espace: et de la variation
  des cordonn{\'e}es provenant de ces d{\'e}placements consid{\'e}r{\'e}s
  ind{\'e}pendamment des causes qui peuvent les produire}, 1840.

\bibitem{gibbs1961scientific}
J.~W. Gibbs, \emph{The Scientific Papers, Vol II: Dover}.\hskip 1em plus 0.5em
  minus 0.4em\relax Dover Publications, 1961.

\bibitem{wilson1901vector}
E.~B. Wilson and J.~W. Gibbs, \emph{Vector analysis: a text-book for the use of
  students of mathematics \& physics: founded upon the lectures of JW
  Gibbs}.\hskip 1em plus 0.5em minus 0.4em\relax Scribner, 1901.

\bibitem{tsiotras1997higher}
P.~Tsiotras, J.~L. Junkins, and H.~Schaub, ``Higher-order cayley transforms
  with applications to attitude representations,'' \emph{Journal of Guidance,
  Control, and Dynamics}, vol.~20, no.~3, pp. 528--534, 1997.

\bibitem{cayley1845xiii}
A.~Cayley, ``Xiii. on certain results relating to quaternions: To the editors
  of the philosophical magazine and journal,'' \emph{Philosophical Magazine
  Series 3}, vol.~26, no. 171, pp. 141--145, 1845.

\bibitem{lizarralde1996attitude}
F.~Lizarralde and J.~T. Wen, ``Attitude control without angular velocity
  measurement: A passivity approach,'' \emph{IEEE transactions on Automatic
  Control}, vol.~41, no.~3, pp. 468--472, 1996.

\bibitem{mayhew2011quaternion}
C.~G. Mayhew, R.~G. Sanfelice, and A.~R. Teel, ``Quaternion-based hybrid
  control for robust global attitude tracking,'' \emph{IEEE Transactions on
  Automatic Control}, vol.~56, no.~11, pp. 2555--2566, 2011.

\bibitem{joshi1995robust}
S.~Joshi, A.~Kelkar, and J.-Y. Wen, ``Robust attitude stabilization of
  spacecraft using nonlinear quaternion feedback,'' \emph{IEEE Transactions on
  Automatic control}, vol.~40, no.~10, pp. 1800--1803, 1995.

\bibitem{euston2008complementary}
M.~Euston, P.~Coote, R.~Mahony, J.~Kim, and T.~Hamel, ``A complementary filter
  for attitude estimation of a fixed-wing uav,'' in \emph{Intelligent Robots
  and Systems, 2008. IROS 2008. IEEE/RSJ International Conference on}.\hskip
  1em plus 0.5em minus 0.4em\relax IEEE, 2008, pp. 340--345.

\bibitem{marins2001extended}
J.~L. Marins, X.~Yun, E.~R. Bachmann, R.~B. McGhee, and M.~J. Zyda, ``An
  extended kalman filter for quaternion-based orientation estimation using marg
  sensors,'' in \emph{Intelligent Robots and Systems, 2001. Proceedings. 2001
  IEEE/RSJ International Conference on}, vol.~4.\hskip 1em plus 0.5em minus
  0.4em\relax IEEE, 2001, pp. 2003--2011.

\bibitem{breckenridge1999quaternions}
W.~Breckenridge, ``Quaternions proposed standard conventions,'' \emph{Jet
  Propulsion Laboratory, Pasadena, CA, Interoffice Memorandum IOM}, pp.
  343--79, 1999.

\bibitem{bhat2000topological}
S.~P. Bhat and D.~S. Bernstein, ``A topological obstruction to continuous
  global stabilization of rotational motion and the unwinding phenomenon,''
  \emph{Systems \& Control Letters}, vol.~39, no.~1, pp. 63--70, 2000.

\end{thebibliography}
\end{document}